\documentclass[final]{myarticle}

\calclayout
\numberwithin{equation}{section}
\allowdisplaybreaks
\theoremstyle{plain}

\newtheorem{theorem}{Theorem}[section]
\newtheorem{lemma}[theorem]{Lemma}
\newtheorem{corollary}[theorem]{Corollary}
\newtheorem{proposition}[theorem]{Proposition}
\newtheorem{remark}[theorem]{Remark}
\newtheorem{definition}[theorem]{Definition}

\usepackage[pdfstartview=FitH]{hyperref}

\usepackage[notref,notcite]{showkeys}

\usepackage{enumerate}
\usepackage{amsmath}
\usepackage{amssymb}
\usepackage{mathtools}
\mathtoolsset{showonlyrefs}
\usepackage{mathrsfs}
\usepackage{comment}
\usepackage{xcolor}
\usepackage{extarrows}
\usepackage{graphicx}
\usepackage{subfigure}
\usepackage{doi}
\usepackage{cases}
\usepackage{tikz}
\usepackage[margin=1.25in]{geometry}



\def\i{\mathrm{i}}

\begin{document}
	
	\title{Ergodicity for Ginzburg-Landau equation with complex-valued space-time white noise on two-dimensional torus}
	
	\author{Huiping \textsc{Chen}}
	\address{Academy of Mathematics and Systems Science, Chinese Academy of Sciences, Beijing 100190, China}
	\email{chenhp@amss.ac.cn}
	
	\author{Yong \textsc{Chen}}
	\address{School of Big Data, Baoshan University, Baoshan, Yunnan 678000, China}
	\email{zhishi@pku.org.cn; chenyong77@gmail.com}
	
	\author{Yong \textsc{Liu}}
	\address{LMAM, School of Mathematical Sciences, Peking University, Beijing, 100871, China}
	\email{liuyong@math.pku.edu.cn}
	
	
	\keywords{Complex Ginzburg-Landau equation; ergodicity; complex space-time white noise; stochastic heat equation with dispersion; complex Wiener-It\^o integral. \\\indent 	
		We thank Prof. Deng Zhang's valuable disscusions and suggestions to pay attention to \cite{Matsuda20202,Matsuda2020, Trenberth2019}. H. Chen is supported by the China Postdoctoral Science Foundation under Grant Number 2024M763480. Y. Chen is supported by NSFC (No. 12461029). Y. Liu is supported by NSFC (No. 12231002) and Center for Statistical Science, PKU}


	\begin{abstract}
		We investigate the global well-posedness and ergodicity of the complex Ginzburg-Landau equation with a general nonlinear term on the two-dimensional torus, driven by complex-valued space-time white noise. Due to the roughness of noise, the solution to this singular equation is a distribution-valued stochastic process. As a result, the nonlinear term is ill-defined and requires renormalization. We establish global well-posedness by combining the fixed point theorem with an estimate that decays over time. Moreover, we prove ergodicity by applying the Krylov-Bogoliubov theorem along with an asymptotic coupling argument. A crucial tool in our proof is the theory of complex multiple Wiener-It\^o integrals, which enables direct estimates for random distributions themselves and provides a systematic framework for estimating complex Wick products.
	\end{abstract}

	\maketitle
	
	\tableofcontents

	\section{Introduction}
	
	We are interested in the stochastic complex Ginzburg-Landau equation on the two-dimensional torus defined as
	\begin{equation}\label{SCGL}
		\begin{cases}
			\partial_t u= ((\i+\mu) \Delta -1) u-\nu|u|^{2m} u+\tau u+\xi, & t>0, x \in \mathbb{T}^2, \\
			u(0, \cdot)=u_0\in \mathcal{C}^{-\alpha_0},
		\end{cases}
	\end{equation}
	where $\mathcal{C}^{-\alpha_0}$ is a Besov space with regularity $-\alpha_0<0$ (see Section \ref{Appendix A Besov space and Sobolev space} for the definition), $\mathbb{R}\ni\mu>0$, $\nu, \tau \in \mathbb{C}$, $\mathrm{Re}\, \nu>0$, $m\geq1$ is an integer and $\xi$ is a complex-valued space-time white noise as introduced in Section \ref{Sec 2.1 Complex multiple Wiener-Ito integral}. Here, we call $\i \Delta  u$ the dispersion term and $\mu \Delta  u$ the dissipation term. We aim to obtain the global well-posedness and ergodicity of \eqref{SCGL}.

	The complex Ginzburg-Landau equation is one of the most important nonlinear partial differential equations, describing a wide range of physical phenomena such as nonlinear waves, superconducting phase transitions and superfluidity (see \cite{RevModPhys.74.99}). Beyond its extensive applications in physics, it also holds significant mathematical interest as a deterministic or stochastic parabolic equation. The stochastic quantization approach, introduced by Parisi and Wu in \cite{ParisiWu1981}, provides a framework for constructing a quantum field theory via invariant measures of stochastic processes. We refer to \cite{Glimm1987,Simon1974} for more background and details. In this paper, we establish the global well-posedness and ergodicity of \eqref{SCGL}, which are essential to analyze the complex-valued $\Phi_2^{2(m+1)}$ measure in quantum field theory.

	For Ginzburg-Landau equation in the random setting, we refer to \cite{BartonSmith2004,BartonSmith20042,Kuksin2004,Odasso2006} for studies involving non-white or multiplicative noise. For the case of additive space-time white noise, Hairer in \cite{MH2002} established the exponential mixing property for the complex cubic Ginzburg-Landau equation driven by real-valued noise in one spatial dimension. When the spatial dimension $d\geq 2$, the complex Ginzburg-Landau equation driven by complex-valued space-time white noise becomes a singular stochastic partial differential equation (SPDE) due to the irregularity of the noise. For $d=2$ as we considered in this paper, one can use Da Prato-Debussche method (see \cite{DaPrato2003}) to handle the singularity. The theory of regularity structures introduced by Hairer in \cite{Hairer2014}, and the theory of paracontrolled distributions proposed by Gubinelli, Imkeller and Perkowski in \cite{Gubinelli2015} provide two strong tools for studying singular SPDEs. For $d=3$, the global well-posedness for the singular stochastic Ginzburg-Landau equation on $\mathbb{T}^3$ were derived in \cite{Hoshino2018,Hoshino2017} by using methods in \cite{Gubinelli2015,Hairer2014}.

		In contrast to its real counterpart, stochastic quantization equation, analyzing \eqref{SCGL} requires a more subtle approach. For instance, unlike the stochastic quantization equation where we can obtain $L^\infty$ estimates for the remainder, for \eqref{remainder in introduction} below, due to the presence of the dispersion term $\i \Delta  v$, we can only obtain $\left\|v\right\| _{L^{2p}}$ estimates for $p>1$ but close to $1$. This necessitates a slightly larger dissipation coefficient $\mu $ in the proofs of well-posedness and ergodicity. It is only when $m=1$ that the regularizing effect of the heat kernel can help us upgrade the $\left\|v\right\| _{L^{2p}}$ estimates to the required level, thereby removing the requirement for $\mu $. When $m>1$, this procedure of regularity enhancement is no longer applicable. Secondly, the presence of the nonlinear term $|u|^{2m} u$ induces a complex interaction between the real and imaginary parts of $u$, which prevents us from analyzing $u$ by separately considering its real and imaginary components. These points collectively demonstrate that \eqref{SCGL} possesses its own unique significance for study and cannot be simply regarded as a complex version of the stochastic quantization equation or as a real two-dimensional equation.
	
	\subsection{Main results}

	We utilize the Da Prato-Debussche method (see \cite{DaPrato2003}), decomposing $u$ into the solution of the stochastic heat equation with dispersion and a remainder respectively defined as
	\begin{align}
		&\begin{cases}
			\partial_t Z=\left[ (\i +\mu)\Delta-1\right]Z+ \xi, & t>0,x\in \mathbb{T}^2,\\
			Z(0,\cdot)=0,
		\end{cases}\label{heat equation in introduction}\\
		&\begin{cases}
			\partial_t v=\left[ (\i+\mu) \Delta-1\right]  v+	\Psi(v,\underline{Z}) , & t>0, x \in \mathbb{T}^2, \\
			v(0, \cdot)=u_0\in \mathcal{C}^{-\alpha_0},
		\end{cases}\label{remainder in introduction}
	\end{align} 
	where $	\Psi\left( v,\underline{Z}\right)$ is defined as
	\begin{equation*}\label{Psi in introduction}
		\Psi\left( v,\underline{Z}\right) :=-\nu\sum_{i=0}^{m+1}\sum_{j=0}^{m}\binom{m+1}{i}\binom{m}{j}v^i\overline{v}^jZ^{:m+1-i,m-j:} +\tau(v+Z),
	\end{equation*}
	and $Z^{:k,l:}$ for $0\leq k\leq m+1, 0\leq l\leq m$ are Wick products of \eqref{heat equation in introduction} defined via complex multiple Wiener-It\^o integrals in Section \ref{Stochastic heat equation with dispersion}. Then we defined $u:= v+Z$ as the solution to renormalized \eqref{SCGL}.
	
	The main results in this paper are as follows. 
	\begin{enumerate}
		\item Regularity theorem for general distribution-valued (namely $\mathscr{S}^{\prime}(\mathbb{T}^d)$-valued with $d\geq 1$) random processes with complex multiple Wiener-It\^o integral type (see \eqref{good process} for the definition), which include the Wick products of the stochastic heat equation defined on $\mathbb{T}^d$.
		\item Global well-posedness of the equation for the remainder \eqref{remainder in introduction}.
		\item Ergodicity for the renormalized Ginzburg-Landau equation \eqref{SCGL} by an asymptotic coupling argument.
	\end{enumerate}

	We next state the theorems specifically. Regardless of the spatial dimension $d\geq 1$ and the complexity of the structure of the target processes defined on $\mathbb{T}^d$, as long as they are of the complex multiple Wiener-It\^o integral type (see \eqref{good process}) and their Fourier coefficients satisfy conditions \eqref{eq: 2nd est for diff freq for general distribution}, \eqref{eq: 2nd est for general distribution} and \eqref{eq: time diff est for general distribution}, we can easily obtain their regularity. This simple criterion is establish by using isometry and hypercontractivity properties (see \cite[Theorem 7]{ito1952complex} and \cite[Proposition 3.9]{chenliu2019} respectively) of complex multiple Wiener-It\^o integrals. Applying the following theorem, we prove the regularity of $Z^{:k,l:}$ for $0\leq k\leq m+1, 0\leq l\leq m$ in Theorem \ref{regularity of Z}.
	\begin{theorem}\label{Prop: regularity for general distribution}
		Let $(X(t))_{t\geq 0}$ be a $\mathscr{S}^{\prime}(\mathbb{T}^d)$-valued random process with the form \eqref{good process}, where $d \geq 1$ is spatial dimension. If for any $t, t_1, t_2\geq 0$, $m, m_1, m_2\in \mathbb{Z}^d$ satisfying $m_1\neq m_2$, and any $0<\epsilon<1$,
		\begin{align}
			&\mathrm{E}\left[ \hat{X} (t_1,m_1)\overline{\hat{X} (t_2,m_2)}\right]=0,\label{eq: 2nd est for diff freq for general distribution}\\
			&\mathrm{E}\left[ \left| \hat{X} (t,m)\right| ^2\right]\lesssim\frac{1}{(1+|m|^2)^{1-\epsilon}},\label{eq: 2nd est for general distribution}\\
			&	\mathrm{E}\left[\left|  \hat{X} (t_2,m)- \hat{X} (t_1,m)\right| ^2 \right]\lesssim |t_2-t_1|^{\epsilon}\frac{1}{(1+|m|^2)^{1-\epsilon}}, \label{eq: time diff est for general distribution}
		\end{align}
		then for any $\alpha>0$, $X\in C\left( [0,\infty);\mathcal{C}^{1-d/2-\alpha}\right) $ almost surely. Moreover, for every $T>0$ and for any $p\geq1$, there exist a constant $C$ depending on $\alpha$, $T$, $p$, and $0<\lambda<1$ such that
		\begin{align}
			&\mathrm{E}\left[ 	\sup_{0\leq t\leq T} \left\|  X(t)\right\| _{{\mathcal{C}^{1-d/2-\alpha}}}^{p}\right]\leq C,\label{eq: regularity for general distribution}\\
			&	\mathrm{E}\left[\sup_{0\leq t_1<t_2\leq T}	\frac{ \left\| X(t_2)-X(t_1)\right\| _{{\mathcal{C}^{1-d/2-\alpha}}}^{p}}{(t_2-t_1)^{\lambda p}}\right]\leq C. \label{eq: time regularity for general distribution}
		\end{align}
	\end{theorem}
	
	Using complex Wiener-It\^o integrals for renormalizing complex-valued singular SPDEs offers four advantages. First, building on It\^o's work \cite{ito1952complex} and the authors' previous work on the theory of complex Wiener-It\^o integrals (see \cite{ChenHP2024,ccl22,ccl24b,ccl25,Chen2017,ChenLiu2014,chen2017fourth,chenliu2019}), it provides a systematic and concise framework for estimating Wick products. Second, it allows direct estimates for random distributions themselves. An alternative approach in \cite{RZ2025,Trenberth2019} utilizing generalized Laguerre polynomials defined only for the truncated distributions, lead to complicated proofs for regularity of a distribution through truncated sequences. Third, Theorem \ref{Prop: regularity for general distribution} leverages the definition of complex Besov spaces, improving the method in \cite[Proposition 5.1]{ccl22} by avoiding the exponential increase in dimensions from decomposing a complex Wiener-It\^o integral into real and imaginary parts. Fourth, we consider a general spatial dimension $d \geq 1$ and when $d=3$, the conditions in Theorem \ref{Prop: regularity for general distribution} are easier to be verified compared to \cite{Hoshino2018,Hoshino2017,Matsuda2020}. Specifically, conditions in Theorem \ref{Prop: regularity for general distribution} are for the Fourier coefficients of random distributions, which are finer structures, rather than for Littlewood-Paley blocks as considered in \cite[Proposition 5.3]{Hoshino2017} and \cite[Theorem 2.1]{Matsuda2020}.

	After understanding the properties of Wick products of \eqref{heat equation in introduction} well using Theorem \ref{Prop: regularity for general distribution}, we consider the global well-posedness of deterministic equation \eqref{remainder in introduction} for the remainder $v$, which also implies the global well-posedness of renormalized \eqref{SCGL}. Using the fixed point argument, in Theorem \ref{Local existence and uniqueness}, we obtain the local well-posedness of \eqref{remainder in introduction} in the solution space $C\left( \left( 0,T\right]; \mathcal{C}^{\beta}\right) $ with the norm $\sup_{0< t\leq T}t^{\gamma}\left\|v(t) \right\|_{\mathcal{C}^{\beta}} <\infty$, where $T>0$ and $\beta,\gamma$ satisfy some conditions (see \eqref{eq: para cond}). Combining local well-posedness with a priori estimates in Proposition \ref{priori estimate}, we derive the global well-posedness for \eqref{remainder in introduction} with $m>1$.
	
	\begin{theorem}[Global well-posedness for $m>1$]\label{Global existence and uniqueness}
		Let $m>1$, $\mu>\frac{2m-1}{2\sqrt{2m}}$ and $\alpha_0,\beta$, $\gamma$ satisfy 
		\begin{align}\label{eq: cond. for alpha0 beta gamma}
			\frac{1}{1+\mu(\mu+\sqrt{1+\mu^2})}\leq \alpha_0\leq \frac{2}{2m+1}, \quad \beta>0,\quad \frac{\alpha_0+\beta}{2}< \gamma<\frac{1}{2m+1}.
		\end{align}
		Then for every $u_0\in \mathcal{C}^{-\alpha_0}$, there exists a unique mild solution $v\in C\left( \left( 0,\infty\right); \mathcal{C}^{\beta}\right)$ of \eqref{remainder in introduction}. 
	\end{theorem}
	
	 The assumption on dissipation coefficient $\mu$ for $m>1$ is natural in the view of \cite[Theorem 4.1]{Doering1994}, where the global well-posedness for the deterministic complex Ginzburg-Landau equation was obtained. For $m=1$, using an improved estimate shown in Proposition \ref{prop: B1+ est}, we establish the global well-posedness for any $\mu>0$. Similar estimates can also be found in \cite[Theorem 3.3]{Matsuda2020}. However, it is worth noting that we have demonstrated that such an improved estimate exists only when $m = 1$. For a detailed explanation, please refer to Remark \ref{rmk: no improved est for m>1}.
	
	\begin{theorem}[Global well-posedness for $m=1$]\label{Global existence and uniqueness for m=1}
		Let $m=1$, $\mu>0$ and $\alpha_0,\beta$, $\gamma$ satisfy $$  \alpha_0> 0, \quad \beta>0,\quad \frac{\alpha_0+\beta}{2}< \gamma< \frac13.$$ Then for every $u_0\in \mathcal{C}^{-\alpha_0}$, there exists a unique mild solution $v\in C\left( \left( 0,\infty\right); \mathcal{C}^{\beta}\right)$ of \eqref{remainder in introduction}. 
	\end{theorem}

	We establish an estimate for $v$ that decays over time and is independent of the initial data in Proposition \ref{priori estimate} by using a coming down from infinity argument presented in \cite[Lemma 3.8]{TW2018}. This estimate not only can helps us prove Theorem \ref{Global existence and uniqueness}, but also is crucial for establishing the ergodicity of renormalized \eqref{SCGL}. While global well-posedness part was also addressed in \cite{Trenberth2019} and later extended to compact surfaces in \cite{RZ2025} following the argument in \cite{MourratWeber2017}, the estimates in \cite[Lemma 4.6]{RZ2025} and \cite[Proposition 6.3]{Trenberth2019}, depend on the initial data and is proportional to time, and thus are insufficient to prove ergodicity.

	Combining the estimates for $v$ in Proposition \ref{priori estimate} and  Proposition \ref{prop: B1+ est} with regularity for Wick products of \eqref{heat equation in introduction} in Theorem \ref{regularity of Z}, we show compact estimates for $u$ in Corollary \ref{priori extimate for u} and Corollary \ref{priori extimate for u with m=1}. These estimates are essential for proving the existence of invariant probability measures for $u$ by the Krylov-Bogoliubov theorem. For the uniqueness of invariant probability measure, we utilize an asymptotic coupling argument, which is introduced by Hairer, Mattingly and Scheutzow in \cite{HairerMattinglyScheutzow2011} and applied to stochastic quantization equations by R\"{o}ckner, Zhu and Zhu in \cite{RZZ2017}. Specifically, we derive a priori estimates of auxiliary system (see \eqref{auxiliary system}) by using the theory of complex Wiener-It\^o integrals. Then we construct a specific coupling for \eqref{SCGL} using this auxiliary system in the spirit of \cite{RZZ2017} and get the ergodicity. 
	\begin{theorem}\label{ergodicity of u}
		Let $m\geq1$ and $\alpha_0,\beta$, $\gamma$ satisfy \eqref{eq: cond. for alpha0 beta gamma}. There exists $\mu_0>0$ such that for any $\mu>\mu_0$, $u:=v+Z$ admits a unique invariant probability measure $\eta$ on $ \mathcal{C}^{-\alpha_0}$, where $Z$ and $v\in C\left( \left( 0,\infty\right); \mathcal{C}^{\beta}\right)$ are the solutions of \eqref{heat equation in introduction} and \eqref{remainder in introduction} respectively.		
	\end{theorem}

	\begin{remark}\label{rem: explaination for large mu}
		See \eqref{condition on mu} in Section \ref{sec: Proof of Ergodicity} for the definition of $\mu_0$. The technical assumption $\mu>\mu_0$ is necessary based on a priori estimates obtained in Proposition \ref{priori estimate} and Proposition \ref{priori estimate for auxiliary system}. However, for $m=1$, we can remove this assumption and establish the existence of invariant measures for any $\mu>0$ in Proposition \ref{existence of IM for m=1}, by applying an improved estimate from Proposition \ref{prop: B1+ est}.
	\end{remark}

		The asymptotic coupling method offers two advantages. First, it allows us to consider general nonlinearities ($m\geq1$). In \cite{Matsuda20202} and \cite{Matsuda2020}, motivated by \cite{TW2018}, Matsuda established the global well-posedness, strong Feller property, and support theorem for the case $m=1$. Regarding exponential ergodicity, he remarked that it might be possible to prove it by combining the ideas from \cite[Section 6]{TW2018} with the results in \cite{Matsuda20202} and \cite{Matsuda2020}. The method in \cite[Section 6]{TW2018} only yields a support theorem for the cubic nonlinearity and does not handle higher-order nonlinearities. Therefore, even if exponential ergodicity could be established, it would likely be limited to the case $m=1$. Second, its proof is simple. One can directly derive a priori estimates for the auxiliary system using arguments similar to those used for establishing the global well-posedness, which reduces overall length of the proof for ergodicity. In this work, we consider a general nonlinearity ($m\ge 1$) and concisely prove ergodicity using an asymptotic coupling approach.

	In \cite{RZZ17}, R\"{o}ckner, Zhu and Zhu proved that the $P(\Phi)_2$ quantum field is an invariant measure of the stochastic quantization equation. They achieved this by using solutions given by Dirichlet form theory in \cite{AR1991}. For further results on Dirichlet forms associated with stochastic quantization problems, we refer readers to \cite{ZZ182,ZZ18}. In future work, we will continue to investigate the explicit formulation of the invariant measure, namely complex-valued $\Phi_2^{2(m+1)}$ measure, and the exponential ergodicity of \eqref{SCGL}.

	\subsection{Existing results in literature}

	We provide a concise overview of both real and complex Wiener-It\^o integrals. For further details on complex Wiener-It\^o integrals, we refer to Section \ref{Sec 2.1 Complex multiple Wiener-Ito integral}. In 1951, It\^o published a seminal article \cite{ito1951real}, defining multiple Wiener integrals with respect to a normal random measure which was first introduced and termed polynomial chaos by Wiener in \cite{Wiener1938}. It\^o in \cite[Theorem 3.1]{ito1951real} revealed the close connection between real Hermite polynomials and multiple Wiener-It\^o integrals. Shortly thereafter, It\^o in \cite{ito1952complex} established the theory of complex multiple Wiener-It\^o integrals with respect to a complex normal random measure. It\^o showed isometry property, chaos expansion and a deep relation between complex Hermite polynomials (also named Hermite-Laguerre-It\^o polynomials in \cite{ChenLiu2014}) and complex multiple Wiener-It\^o integrals.

	Since then, theoretical interest in complex Gaussian systems, specifically complex Wiener-It\^o integrals and complex Hermite polynomials, has grown significantly. In \cite{ChenLiu2014,chenliu2019}, the last two authors of this paper explored properties of complex multiple Wiener-It\^o integrals, complex Ornstein-Uhlenbeck operators and semigroup. Especially, they proved the hypercontractivity property, namely the equivalence of the $L^r$ norm, where $r \geq 2$, for complex Wiener-It\^o integrals. In \cite{ccl22,ccl24b,Chen2017,chen2017fourth}, the relationship between real and complex Wiener-It\^o integrals, the asymptotic properties, product formula and independence of complex Wiener-It\^o integrals were established. Recent progress on complex Wiener-It\^o integrals was summarized in \cite{ccl25}. Further studies on complex Hermite polynomials have explored their analytical properties, connections with real Hermite polynomials, and related topics, see \cite{Agorram2016,Cotfas_2010,Ismail2016} for detailed discussions.

	Using complex Wiener-It\^o integrals for renormalizing is a natural choice inspired by the renormalization for real-valued singular SPDEs utilizing real Wiener-It\^o integral. In the real case, extensive research has employed Hermite polynomials and Wiener-It\^o integrals to renormalize singular SPDEs. In \cite{DaPrato2003}, Da Prato and Debussche defined Wick products through real Hermite polynomials for the approximate stochastic heat equation driven by real-valued space-time white noise, to investigate the global well-posedness of the $\Phi_2^{n+1}$ model ($n\geq3$) for almost all initial data with respect to the invariant measure. They proved the convergences of the approximate Wick products and defined the limits as the Wick products of the stochastic heat equation. In \cite{Hairer2014}, Hairer introduced the theory of regularity structures, which provides a powerful framework for the renormalization of singular SPDEs. As applications, he renormalized the generalized parabolic Anderson model and the $\Phi_3^4$ model, and established their convergence. The theory incorporates multiple Wiener-It\^o integrals as a method to represent the stochastic terms and handle the divergences arising from the nonlinearities. Using the method in \cite{Gubinelli2015}, where Gubinelli, Imkeller and Perkowski proposed an alternative framework for handling singular SPDEs based on the theory of paracontrolled distributions, Catellier and Chouk in \cite{CC2018AP} employed Hermite polynomials for renormalization and established local well-posedness for the $\Phi_3^{4}$ model. Based on this work, Mourrat and Weber in \cite{MW2017} obtained a priori estimates that ultimately lead to the global well-posedness of the $\Phi_3^{4}$ model. R\"{o}ckner, Zhu and Zhu in \cite{RZZ2017} also utilized real Hermite polynomials to renormalize the $\Phi_2^{4}$ model and established ergodicity through the generalized coupling method. Mourrat and Weber in \cite{MourratWeber2017}, and also Tsatsoulis and Weber in \cite{TW2018}, demonstrated that the Wick products of the stochastic heat equation correspond precisely to real multiple Wiener-It\^o integrals. Specifically, Mourrat and Weber in \cite{MourratWeber2017} established global well-posedness for the $\Phi_2^{4}$ model on the plane, while Tsatsoulis and Weber in \cite{TW2018} obtained exponential ergodicity for the $\Phi_2^{n+1}$ model on the torus by combining the strong Feller property with the support theorem. The convergences of the approximate Wick products involved in the $\Phi_3^{4}$ model and the multiple Wiener-It\^o integral representations of the limits can also be found in \cite{MWX2017}. The result that the limits coincide with multiple Wiener-It\^o integrals appears natural in light of \cite[Theorem 3.1]{ito1951real}.

	In \cite{ccl22}, we establish the regularity of \eqref{heat equation} by combining \cite[Proposition 5]{MWX2017}, which concerns the regularity of real-valued random processes, with the fact that the real and imaginary parts of a complex $(p,q)$-th Wiener-It\^o integral can be expressed as a real $(p+q)$-th Wiener-It\^o integral respectively. While this strategy proved effective for low-order Wick products, its applicability diminishes for high-order Wick products, because the kernels of the real and imaginary parts of a complex $(p,q)$-th Wiener-It\^o integral have dimensions of $2 ^{p+q}$. This exponential increase in dimensions complicates calculations when estimating high-order Wick products by this approach. However, when handling \eqref{SCGL} with a general nonlinearity ($m\geq1$), estimating high-order Wick products becomes unavoidable. To overcome the limitations caused by decomposing real and imaginary parts, we now adopt a new approach based on the definition of complex Besov spaces, developing novel techniques for estimating complex Wick products. This method resolves the previous challenges, and we anticipate that it will also be valuable for addressing related problems in future research.

	After our work \cite{ccl24c} appeared on arXiv, Robert and Zine in \cite{RZ2025} extended the global well-posedness part to the equations defined on compact surfaces. Their proof relies on a priori estimates that depend on the initial data and increase over time, as well as on generalized Laguerre polynomials for renormalization.

	This paper is organized as follows. In Section \ref{Appendix A Besov space and Sobolev space}, we introduce fundamental concepts and properties of Besov and Sobolev spaces. In Section \ref{Stochastic heat equation with dispersion}, we show the regularity of Wick products associated with the stochastic heat equation with dispersion \eqref{heat equation} using the theory of complex multiple Wiener-It\^o integrals. In Section \ref{Global well-posedness}, combining the Da Prato-Debussche method (see \cite{DaPrato2003}) with a priori estimates, we obtain the global well-posedness of the stochastic complex Ginzburg-Landau equation \eqref{SCGL} in the renormalized sense. In Section \ref{Ergodicity}, using the asymptotic (or generalized) coupling approach developed in \cite{HairerMattinglyScheutzow2011, KulikSch2018}, we demonstrate the existence and uniqueness of the invariant measure for the renormalized solution to \eqref{SCGL}. Some calculations and technical estimates are provided in Appendix \ref{Some technical estimation}.

	\section{Preliminary: Besov and Sobolev space}\label{Appendix A Besov space and Sobolev space}
	
	To characterize the regularity of the solution to the stochastic heat equation with dispersion \eqref{heat equation}, we introduce the Besov and Sobolev spaces, more details about which can be found in \cite[Chapter 1 and 2]{Bahouri2011}. Readers familiar with Besov and Sobolev spaces may skip this section and proceed directly to Section \ref{Stochastic heat equation with dispersion}, where the regularity of Wick products associated with \eqref{heat equation} is established using the theory of complex multiple Wiener-It\^o integrals.
	
	Let $d \geq 1$ denote the spatial dimension. Unless explicitly stated, all function spaces we considered throughout this paper are spaces of complex-valued functions with specific properties. For example, we denote by $C^{\infty}( \mathbb{T}^d) $ the space of all complex-valued smooth functions defined over $\mathbb{T}^d$. We denote the space of complex-valued Schwartz distributions by $\mathscr{S}^{\prime}(\mathbb{T}^d)$. For $p\geq1$, let $L^p(G) $ be the space of all complex-valued $p$-integrable functions defined on a domain $G$. Without ambiguity, we sometimes omit $G$ and use $L^p$ to denote this function space. For complex-valued $f,g\in L^2(G)$, we set $$\left\langle f,g\right\rangle=\int_{ G}f(x)g(x)\mathrm{d}x.$$
	Note that we do not take complex conjugate for $g$. Then the norm of $f\in L^2(G)$ is 
	\begin{equation*}
		\left\| f\right\| = \left\langle f,\overline{f}\right\rangle^{\frac12}.
	\end{equation*}
	We also use the notation $\left\langle f,g\right\rangle$ for a pairing of $f\in \mathscr{S}^{\prime}(\mathbb{T}^d)$ and $g\in C^{\infty}( \mathbb{T}^d) $. We write the Fourier coefficient of $f\in \mathscr{S}^{\prime}(\mathbb{T}^d)$ with frequency $k\in \mathbb{Z}^d$ as $$
	\hat{f}(k):=\left\langle f, \frac{1}{(2  \pi)^{d/2}} e^{-\i k \cdot x} \right\rangle =\frac{1}{(2  \pi)^{d/2}}\int_{\mathbb{T}^d} f(x) e^{-\i k \cdot x} \mathrm{d} x.$$ 
	Let $\left\lbrace \chi_j\right\rbrace _{j=-1}^{\infty}$ be a dyadic partition of unity satisfying the following properties:
	\begin{enumerate}
		\item $\chi_{-1}, \chi_0:\mathbb{R}^d\rightarrow [0,1]$ are radial and infinitely differentiable;
		\item $\mathrm{supp}(\chi_{-1})\subseteq B(0,\frac43)$ and $\mathrm{supp}(\chi_{0})\subseteq B(0,\frac83)\setminus B(0,\frac34)$, where $B(0,r):=\left\lbrace x\in \mathbb{R}^d: |x|\right. $ $\left. <r\right\rbrace $;
		\item For $j\geq0$, $\chi_j(\cdot)=\chi_{0}(\cdot/2^j)$;
		\item For any $x\in \mathbb{R}^d$, $\sum_{j=-1}^{\infty}\chi_{j}(x)=1$.
	\end{enumerate}
	The existence of dyadic partitions of unity can be found in \cite[Proposition 2.10]{Bahouri2011}. For any $f \in C^{\infty}(\mathbb{T}^d)$ and $j \geq-1$, we define the $j$-th Littlewood-Paley block as
	\begin{equation}\label{eq: def of Littlewood-Paley block}
		(\delta_j f)(x):=\sum_{k \in \mathbb{Z}^d}\chi_j(k) \hat{f}(k) \frac{1}{(2  \pi)^{d/2}}e^{ \i k \cdot x}=(\eta_{j}\ast f)(x)= \left\langle f, \eta_j(x-\cdot)\right\rangle ,
	\end{equation}
	where $\eta_j(x)=\sum_{k\in\mathbb{Z}^d} \frac{1}{(2  \pi)^{d/2}} \chi_j(k) \frac{1}{(2  \pi)^{d/2}}e^{ \i k \cdot x} \in C^{\infty}( \mathbb{T}^d) $. For $\alpha\in \mathbb{R}$ and $p,q\in [1,\infty]$, the Besov space $\mathcal{B}_{p,q}^{\alpha}\subset \mathscr{S}^{\prime}(\mathbb{T}^d)$ is defined as the completion of $C^{\infty}(\mathbb{T}^d)$ with respect to the norm
	$$
	\|f\|_{\mathcal{B}_{p,q}^{\alpha}}:=\Big( \sum _{j \geq-1} 2^{\alpha jq}\left\|\delta_j f\right\|_{L^{p}}^q \Big) ^{\frac{1}{q}}.
	$$
	When $p=q=\infty$, we denote by $\mathcal{C}^\alpha$ the space $\mathcal{B}_{\infty,\infty}^{\alpha}$ for $\alpha\in  \mathbb{R}$. That is,
	the space $\mathcal{C}^\alpha$ is defined as the completion of $C^{\infty}(\mathbb{T}^d)$ with respect to the norm
	$$
	\|f\|_{\mathcal{C}^\alpha}:=\sup _{j \geq-1} 2^{\alpha j}\left\|\delta_j f\right\|_{L^{\infty}}.
	$$
	
	We use $A \lesssim B$ to denote that there exists some constant $C>0$ such that $A \leq C B$. The constants $C$ we omitted may vary from line to line to simplify notations. We use $b\wedge c$ and $b \vee c$ to denote
	$\min\{b,c\}$ and $\max\{b,c\}$ respectively for some $b,c\in\mathbb{R}$.
	
	The following proposition (see \cite[Theorem 6.2.4]{Bergh1976} or \cite[Proposition 2.2 and 2.3]{Sawano2018}) explains the roles of parameters $\alpha$, $p$ and $q$ in Besov space $\mathcal{B}_{p,q}^{\alpha}$, where $\alpha\in \mathbb{R}$ and $p, q \in[1, \infty]$.
	\begin{proposition}\label{embedding} 
		Let $\alpha_1, \alpha_2 \in \mathbb{R}$, $p_1, p_2, q_1, q_2 \in[1, \infty]$. Then,
		\begin{align*}
			\|f\|_{\mathcal{B}_{p_1, q_1}^{\alpha_1}} & \lesssim \|f\|_{\mathcal{B}_{p_1, q_1}^{\alpha_2}}, \quad \text { whenever } \alpha_1 \leq \alpha_2, \\
			\|f\|_{\mathcal{B}_{p_1, q_1}^{\alpha_1}} & \lesssim\|f\|_{\mathcal{B}_{p_1, q_2}^{\alpha_1}}, \quad \text { whenever } q_1 \geq q_2, \\
			\|f\|_{\mathcal{B}_{p_1, q_1}^{\alpha_1}} & \lesssim\|f\|_{\mathcal{B}_{p_2, q_1}^{\alpha_1}}, \quad \text { whenever } p_1 \leq p_2, \\
			\|f\|_{\mathcal{B}_{p_1, q_1}^{\alpha_1}} & \lesssim\|f\|_{\mathcal{B}_{p_1, q_2}^{\alpha_2}}, \quad \text { whenever } \alpha_1<\alpha_2 .
		\end{align*}
	\end{proposition}

	From \cite[Theorem 6.2.4]{Bergh1976} or \cite[Proposition 2.1]{Sawano2018}, we get that $\mathcal{B}_{p,q}^{0}$ is close to $L^p$ for $p,q \in[1, \infty]$.

	\begin{proposition}\label{L^p embedding}
		Let $p \in[1, \infty]$. The space $\mathcal{B}_{p, 1}^0$ is continuously embedded in $L^p$ and $L^p$ is continuously embedded in $\mathcal{B}_{p, \infty}^0$. Namely,
		$$
		\|f\|_{\mathcal{B}_{p, \infty}^0} \lesssim \|f\|_{L^p} \lesssim\|f\|_{\mathcal{B}_{p, 1}^0} .
		$$
	\end{proposition}

	The following embedding proposition can be found in \cite[Proposition 2.71]{Bahouri2011}.
	\begin{proposition}\label{regularity embedding}
		Let $\alpha \in\mathbb{R}$, $r\in[1, \infty]$ and $1\leq q \leq p\leq \infty$. Then
		$$
		\|f\|_{\mathcal{B}_{p, r}^\alpha} \lesssim \|f\|_{\mathcal{B}_{q, r}^{\alpha+d\left(1/q-1/p\right)}} .
		$$
	\end{proposition}
	
	From \cite[Proposition 4]{MourratWeber2017}, we know the following interpolation proposition.
	\begin{proposition}\label{interpolation}
	 Let $\alpha_0, \alpha_1 \in \mathbb{R}$, $p_0, q_0, p_1, q_1 \in[1, \infty]$, $\theta \in[0,1]$, $\alpha:=(1- \theta) \alpha_0+ \theta \alpha_1$ and $p, q \in[1, \infty]$ such that
	\begin{align*}
		\frac{1}{p}=\frac{1-\theta}{p_0}+\frac{\theta}{p_1}, \quad \frac{1}{q}=\frac{1-\theta}{q_0}+\frac{\theta}{q_1} .
	\end{align*}
	Then, we have
	\begin{align*}
		\|f\|_{\mathcal{B}_{p, q}^{\alpha}} \leq \|f\|_{\mathcal{B}_{p_0, q_0}^{\alpha_0}}^{1-\theta} \|f\|_{\mathcal{B}_{p_1, q_1}^{\alpha_1}}^{\theta}.
	\end{align*}
	\end{proposition}

	The smoothing effect of the heat semigroup $\left(P_t=\mathrm{e}^{A  t}\right)_{t \geq 0}$ with the generator $A:=\left(\i+\mu \right) \Delta-1$ is shown in the following proposition (see \cite[Proposition 5]{MourratWeber2017}).
	\begin{proposition}\label{heat kernel smoothing}
		Let $f \in \mathcal{B}_{p, q}^\alpha$ for $\alpha\in\mathbb{R}$ and $p,q\in[1,\infty]$. For all $\beta \geq \alpha$ and every $t >0$,
		$$
		\left\|P_t f\right\|_{\mathcal{B}_{p, q}^\beta} \lesssim t^{\frac{\alpha-\beta}{2}}\|f\|_{\mathcal{B}_{p, q}^\alpha}.
		$$
	\end{proposition}
	The following proposition (see \cite[Proposition 6]{MourratWeber2017}) shows that, for $0 < \beta-\alpha \leq 2$ and $f \in \mathcal{B}_{p, q}^\beta $, the mapping $t \mapsto P_tf$ from $[0,\infty)$ to $\mathcal{B}_{p, q}^\alpha$
	is continuous. When $\beta=\alpha$, this is still valid. That is, for $\alpha\in \mathbb{R}$ and $f \in \mathcal{B}_{p, q}^\alpha $, the mapping $t \mapsto P_tf$ from $[0,\infty)$ to $\mathcal{B}_{p, q}^\alpha$ is continuous. 
	\begin{proposition}\label{1-heat kernel smoothing}
		Let $\alpha,\beta\in\mathbb{R}$ satisfy $0 \leq \beta-\alpha \leq 2$ and $p, q \in[1, \infty]$. Then for any $t \geq 0$,
		$$
		\left\|\left(1-P_t\right) f\right\|_{\mathcal{B}_{p, q}^\alpha} \lesssim t^{\frac{\beta-\alpha}{2}}\|f\|_{\mathcal{B}_{p, q}^\beta} .
		$$
	\end{proposition}
	
	The multiplicative structure (see \cite[Theorem 2.82 and 2.85]{Bahouri2011} or \cite[Corollary 1]{MourratWeber2017}) and duality property (see \cite[Proposition 2.76]{Bahouri2011} or \cite[Proposition 7]{MourratWeber2017}) of Besov space are presented in the following two propositions.

	\begin{proposition}\label{multiplicative structure}
		Let $\alpha,\beta\in\mathbb{R}$ satisfy $ \alpha+\beta>0$ and $p,p_1,p_2,q,q_1,q_2\in [1,\infty] $ be such that $\frac1p=\frac{1}{p_1}+\frac{1}{p_2}$ and $\frac1q=\frac{1}{q_1}+\frac{1}{q_2}$. The mapping $(f,g)\mapsto fg$ can be extended to a continuous linear map from $\mathcal{B}_{p_1,q_1}^{\alpha}\times \mathcal{B}_{p_2,q_2}^{\beta}$ to $\mathcal{B}_{p,q}^{\alpha\wedge \beta}$ and for $(f,g)\in\mathcal{B}_{p_1,q_1}^{\alpha}\times \mathcal{B}_{p_2,q_2}^{\beta}$,
		$$\left\|fg \right\|_{\mathcal{B}_{p,q}^{\alpha\wedge \beta}}\lesssim \left\|f \right\|_{\mathcal{B}_{p_1,q_1}^{\alpha}}  \left\|g \right\|_{\mathcal{B}_{p_2,q_2}^{\beta}}.$$
	\end{proposition}

	\begin{proposition}\label{duality property}
		Let $\alpha\in [0,1)$ and $p, q, p^{\prime}, q^{\prime}\in[ 1,\infty]$ be such that $1=\frac{1}{p}+\frac{1}{p^{\prime}}=\frac{1}{q}+\frac{1}{q^{\prime}}$. The mapping $(f,g)\mapsto \left\langle f,g\right\rangle $ can be extended to a continuous bilinear form on $\mathcal{B}_{p, q}^\alpha \times \mathcal{B}_{p^{\prime}, q^{\prime}}^{-\alpha}$ and	for all $(f, g) \in \mathcal{B}_{p, q}^\alpha \times \mathcal{B}_{p^{\prime}, q^{\prime}}^{-\alpha}$,
		$$
		|\langle f, g\rangle| \lesssim \|f\|_{\mathcal{B}_{p, q}^\alpha}\|g\|_{\mathcal{B}_{p^{\prime}, q^{\prime}}^{-\alpha}}.
		$$
	\end{proposition}
	
	For every $\alpha \in(0,1]$, by \cite[Proposition 8]{MourratWeber2017}, we know that the Besov norm $\left\| \cdot\right\| _{\mathcal{B}_{1,1}^\alpha}$ can be controlled by the norms $\|\cdot\|_{L^1}$ and $\|\nabla (\cdot)\|_{L^1}$, where $\nabla$ is the gradient operator.

	\begin{proposition}\label{gradient estimate}
		Let $\alpha \in(0,1]$ and $f \in \mathcal{B}_{1,1}^\alpha$. Then 
		$$
		\|f\|_{\mathcal{B}_{1,1}^\alpha} \lesssim \|f\|_{L^1}^{1-\alpha}\|\nabla f\|_{L^1}^\alpha+\|f\|_{L^1}.
		$$
	\end{proposition}

	Set $\Lambda=(-A)^{\frac12}$. For $s\geq0$ and $1\leq p\leq \infty$, we denote by $H_p^s$ the subspace of $L^p$ consisting of all $f$ that can be written as $f=\Lambda^{-s}g$, where $g\in L^p$. The $H_p^s$ norm of $f\in H_p^s$ is defined as the $L^p$ norm of $g$, that is, $\left\|f \right\| _{H_p^s}=\left\|\Lambda^s f \right\| _{L^p}$. Recall the following embedding results for the Sobolev space $H_p^s$ (see \cite[Theorem 6.4.4]{Bergh1976} and \cite[Theorem 2.13]{Sawano2018}).

	\begin{proposition}\label{Sobolev embedding}
		\begin{enumerate}
			\item Let $s\geq0$, $1<p<\infty$ and $\epsilon>0$. Then $H_p^{s+\epsilon}\subset\mathcal{B}_{p,1}^s\subset\mathcal{B}_{1,1}^s$ and the embedding is continuous and dense.
			\item Let $1\leq p_1\leq p_2\leq\infty$ and $\alpha\geq d(\frac{1}{p_1}-\frac{1}{p_2})$. Then $H_{p_1}^{\alpha}$ is continuously embedded in $H_{p_2}^{\alpha-d(\frac{1}{p_1}-\frac{1}{p_2})}$.
		\end{enumerate}
	\end{proposition}
	The following interpolation inequality and multiplicative inequality for the Sobolev space $H_p^s$ can be found in \cite[Theorem 6.4.5]{Bergh1976}, \cite[Lemma A.4]{Resnick1995} and \cite[Lemma 2.1]{RZZ15}.
	\begin{proposition}\label{interpolation and multiplicative for Sobolev}
		\begin{enumerate}
			\item Let  $s \in(0,1)$ and $p \in(1, \infty)$. Then for $u \in H_p^1$,
			$$
			\|u\|_{H_p^s} \lesssim\|u\|_{L^p}^{1-s}\|u\|_{H_p^1}^s.
			$$
			\item Let $s>0$ and $p \in(1, \infty)$. If $u \in L^{p_1}\cap H^s_{p_4}$ and $v\in L^{p_3}\cap H^s_{p_2}$, where $p_i \in(1, \infty]$ with $i=1, \ldots, 4$ satisfy
			$
			\frac{1}{p}=\frac{1}{p_1}+\frac{1}{p_2}=\frac{1}{p_3}+\frac{1}{p_4}
			$, then
			$$
			\left\|u v\right\|_{H^s_p} \lesssim\|u\|_{L^{p_1}}\left\| v\right\|_{H_{p_2}^s}+\|v\|_{L^{p_3}}\left\| u\right\|_{H^s_{p_4}}.
			$$
		\end{enumerate}
	\end{proposition}

	\section{Stochastic heat equation with dispersion}\label{Stochastic heat equation with dispersion}
	
	\subsection{Proof of Theorem \ref{Prop: regularity for general distribution} by complex multiple Wiener-It\^o integral}\label{Sec 2.1 Complex multiple Wiener-Ito integral}

	We introduce the complex-valued space-time white noise over $\mathbb{R}\times\mathbb{T}^d$, where $d \geq 1$ is spatial dimension, and the complex multiple Wiener-It\^o integral with respect to it. A complex-valued space-time white noise defined on some probability space $(\Omega,\mathcal{F},\mathrm{P})$ and over $\mathbb{R}\times\mathbb{T}^d$ is a family of centered complex Gaussian variables $\left\lbrace \xi(\varphi): \varphi\in L^2(\mathbb{R}\times\mathbb{T}^d )   \right\rbrace $ such that for any $\varphi,\psi\in L^2(\mathbb{R}\times\mathbb{T}^d) $,
	\begin{equation}\label{white noise}
		\mathrm{E}\left[ \xi(\varphi)\xi(\psi)\right]=0,\quad \mathrm{E}\left[ \xi(\varphi)\overline{\xi(\psi)}\right]=\left\langle\varphi,\overline{\psi} \right\rangle .
	\end{equation}
	To get an intuitive description, we construct the complex-valued space-time white noise as follows. Let $\left( B(\cdot, k)=\left( B_1(\cdot, k)+\i B_2(\cdot, k)\right) /\sqrt 2\right) _{k \in \mathbb{Z}^d}$ be a family of independent complex-valued Brownian motions over $\mathbb{R}$ and defined on $(\Omega,\mathcal{F},\mathrm{P})$, where $\left(B_1(t, k),B_2(t, k)\right)_{t\in\mathbb{R}}$ is a two-dimensional standard Brownian motion for a fixed $k\in\mathbb{Z}^d$. Then for $k \in \mathbb{Z}^d$, $\left( B(t, k)\right) _{t\in \mathbb{R}}$ satisfies 
	\begin{equation*}
		\mathrm{E}[B(t, k)] = \mathrm{E}[B^2(t, k)]=0,\quad \mathrm{E}[|B(t, k)|^2]=|t|, \quad t\in\mathbb{R}.
	\end{equation*}
	For every $\varphi\in L^2(\mathbb{R}\times\mathbb{T}^d ) $, we define
	\begin{equation}\label{construction of white noise}
		\xi(\varphi):=\sum_{k \in \mathbb{Z}^d}\int_{\mathbb{R}} \hat{\varphi}(t,-k) \mathrm{d}B(t, k) ,
	\end{equation}
	where the integral is interpreted in the sense of It\^o and $\hat{\varphi}(t,-k)$ denotes $\widehat{\varphi(t,\cdot)}(-k)$. Then $\xi(\varphi)$ with $\varphi\in L^2(\mathbb{R}\times\mathbb{T}^d ) $ is a centered complex Gaussian variable and by It\^o's isometry, for any $\varphi,\psi\in L^2(\mathbb{R}\times\mathbb{T}^d ) $,
	\begin{align}
		\mathrm{E}\left[ \xi(\varphi)\xi(\psi)\right]&=\sum_{k,k^{\prime}\in \mathbb{Z}^d}	\mathrm{E}\left[ \int_{\mathbb{R}} \hat{\varphi}(t,-k) \mathrm{d}B(t, k) \int_{\mathbb{R}} \hat{\psi}(t,-k^{\prime}) \mathrm{d}B(t, k^{\prime})\right] =0, \label{covariance of constructed white noise 1}\\
		\mathrm{E}\left[ \xi(\varphi)\overline{\xi(\psi)}\right]&=\sum_{k,k^{\prime}\in \mathbb{Z}^d}	\mathrm{E}\left[ \int_{\mathbb{R}} \hat{\varphi}(t,-k) \mathrm{d}B(t, k) \int_{\mathbb{R}} \overline{\hat{\psi}(t,-k^{\prime})} \mathrm{d}\overline{B(t, k^{\prime})}\right] \label{covariance of constructed white noise 2}\\
		&=\sum_{k\in \mathbb{Z}^d} \int_{\mathbb{R}} \hat{\varphi}(t,k) \overline{\hat{\psi}(t,k)}\mathrm{d} t =\left\langle\varphi,\overline{\psi} \right\rangle,\nonumber
	\end{align}
	which means that \eqref{white noise} is satisfied. Therefore, $\xi:=\left\lbrace \xi(\varphi), \varphi\in L^2(\mathbb{R}\times\mathbb{T}^d )\right\rbrace $ defined by \eqref{construction of white noise} is a complex-valued space-time white noise. From \eqref{construction of white noise}, \eqref{covariance of constructed white noise 1} and \eqref{covariance of constructed white noise 2}, we see that $\xi$ is formally given by 
	\begin{equation}\label{xi(t,x)}
		\xi(t,x)= \frac{\mathrm{d}}{{\mathrm{d}t}}W(t,x):=\sum_{k \in \mathbb{Z}^d} \frac{\mathrm{d} B(t, k) }{\mathrm{d}t} \frac{1}{2  \pi}e^{\i k \cdot x}, \quad t\in \mathbb{R},\, x\in\mathbb{T}^d ,
	\end{equation}
	with covariance
	\begin{equation*}
		\mathrm{E}\left[ \xi(t,x) \xi(s,y)\right] =0,\quad 	\mathrm{E}\left[ \xi(t,x) \overline{\xi(s,y)}\right]= \mathfrak{D}(t-s) \mathfrak{D}_{d}(x-y) ,\quad s,t\in \mathbb{R},\, x,y\in\mathbb{T}^d,
	\end{equation*}
	where $\mathfrak{D}$ and $\mathfrak{D}_{d}$ are Delta functions over $\mathbb{R}$ and $\mathbb{T}^d$ respectively. 
	
	Let $Z:=\mathbb{R}\times\mathbb{T}^d$ and $L^2(\Omega,\mathcal{F},\mathrm{P})$ denote the space of complex-valued square integrable random variables defined on $\Omega$. Following \cite[Chapter 2]{ito1952complex}, for $p,q\in\mathbb{N}$ with $p+q>0$, we define the complex $(p,q)$-th
	Wiener-It\^o integral of $f\in L^2(Z^{p}\times Z^{q})$ with respect to $\xi$ and denote it by
	\begin{align*}
		I_{p,q}(f)= \int_{(\mathbb{R}\times\mathbb{T}^d)^{p+q}} &f(s_1,y_1,\ldots, s_p,y_p; s_{p+1},y_{p+1},\ldots, s_{p+q}, y_{p+q}) \\
		&\xi(\mathrm{d}s_1,\mathrm{d}y_1)\cdots\xi(\mathrm{d}s_p,\mathrm{d}y_p)\overline{\xi}(\mathrm{d}s_{p+1},\mathrm{d}y_{p+1})\cdots\overline{\xi}(\mathrm{d}s_{p+q},\mathrm{d}y_{p+q}) .
	\end{align*}
	For $f\in L^2(\mathbb{R}\times\mathbb{T}^d)$, $I_{1,0}(f)$ coincides with $\xi(f)$ as defined in \eqref{construction of white noise}. For $p=q = 0$, we define $I_{0,0}$ as the identity map from $\mathbb{C}$ to $\mathbb{C}$. From \cite[Equation (2.15)]{chenliu2019} and \cite[Theorem 7]{ito1952complex}, we know that complex multiple Wiener-It\^o integrals satisfy the following properties. More results about complex multiple Wiener-It\^o integrals can be found in \cite[Theorem 4.1]{c24}, \cite[Theorem 3.1]{ccl24b}, \cite[Theorem 3.3, 3.9 and Corollary 4.14]{ccl22} and \cite[Theorem 1.1 and 3.3]{chen2017fourth}.
	
	\begin{lemma}\label{lem: property for WI}
		Let $p, q, \tilde{p},\tilde{q} \in \mathbb{N}$. For any $f\in L^2(Z^{p}\times Z^{q})$ and  $g\in L^2(Z^{\tilde{p}}\times Z^{\tilde{q}})$, the following statements hold:
		\begin{enumerate}
			\item The mapping $I_{p,q}: L^2(Z^{p}\times Z^{q})\rightarrow L^2(\Omega,\mathcal{F},\mathrm{P})$ is linear.
			\item $I_{p,q}( \mathrm{Sym}(f))=I_{p,q}(f)$, where $\mathrm{Sym}(f)$ is the symmetrization of $f$ defined as
			\begin{equation*}
				\mathrm{Sym}(f)(z_1,\ldots,z_p;z_{p+1},\ldots,z_{p+q})=\frac{1}{p!}\frac{1}{q!}\sum_{\sigma, \pi}f(z_{\sigma(1)},\ldots,z_{\sigma(p)};z_{\pi(p+1)},\ldots,z_{\pi(p+q)}),
			\end{equation*}
			where $z_j=(s_j,y_j)\in Z= \mathbb{R}\times\mathbb{T}^d$ for $1\leq j\leq p+q$, $\sigma$ and $\pi$ run over all permutations of $\{1,\ldots,p\}$ and $\{p+1,\ldots,p+q\}$ respectively.
			\item Complex Wiener-It\^o integrals satisfy the isometry property. That is,
			\begin{align}
				\mathrm{E}\left[ I_{p,q}(f)\overline{I_{\tilde{p},\tilde{q}}(g)}\right]
				&=\mathbf{1}_{\left\lbrace (p,q)=(\tilde{p},\tilde{q})\right\rbrace }	p!q!
				\left\langle \mathrm{Sym}(f), \overline{\mathrm{Sym}(g)}\right\rangle \label{eq: general isometry property} \\
				&= \mathbf{1}_{\left\lbrace (p,q)=(\tilde{p},\tilde{q})\right\rbrace }	p!q! \int_{Z^{p+q}}   \mathrm{Sym}(f) \overline{ \mathrm{Sym}(g)} \mathrm{d} z_1\cdots \mathrm{d}z_{p+q}, \nonumber
			\end{align}
			where $z_j=(s_j,y_j)\in Z= \mathbb{R}\times\mathbb{T}^d$ for $1\leq j\leq p+q$ and $\mathbf{1}_{A}$ denotes the indicator function of a set $A$.
			\item Complex Wiener-It\^o integrals satisfy that 
			\begin{equation}\label{reverse complex conjugate}
				I_{q,p}(h)=\overline{I_{p,q}(f)},
			\end{equation}
			where $h\in L^2(Z^{q}\times Z^{p})$ is the reverse complex conjugate of $f$ defined as 
			\begin{equation*}
				h(z_1,\ldots,z_q;z_{q+1},\ldots,z_{q+p})=\overline{f}(z_{q+1},\ldots,z_{q+p};z_1,\ldots,z_q).
			\end{equation*}
		\end{enumerate}
	\end{lemma}
	By the isometry property \eqref{eq: general isometry property} and Minkowski's inequality, we know that for $f\in L^2(Z^{p}\times Z^{q})$, where $p, q, \in \mathbb{N}$,
	\begin{equation*}
		\mathrm{E}\left[ |I_{p,q}(f)|^2\right]=	p!q!
		\|  \mathrm{Sym}(f)\|^2 _{L^2(Z^{p}\times Z^{q})}\leq 	p!q!   \|  f\|^2 _{L^2(Z^{p}\times Z^{q})}.
	\end{equation*}
	
	We now state the hypercontractivity property, namely the equivalence of the $L^r$ norm, where $r \geq 2$, for complex Wiener-It\^o integral, which was presented in \cite[Proposition 3.9]{chenliu2019}.
	\begin{lemma}[Hypercontractivity property]\label{lem: hypercontractivity}
		For any $f\in L^2(Z^{p}\times Z^{q})$ with some $p,q \in\mathbb{N}$ and any $r \geq 2$, we have
		\begin{equation*}
			\|	I_{p, q}(f)  \|_{L^r(\Omega)} \leq (r-1)^\frac{p+q}{2} \|I_{p, q}(f) \|_{L^2(\Omega)}.
		\end{equation*}
	\end{lemma}

	To prove the regularity of the Wick products of the stochastic heat equation with dispersion defined in \eqref{definition of wick product}, we first establish a straightforward criterion for deriving the regularity of a distribution-valued random process with complex multiple Wiener-It\^o integral type. This criterion simplifies the estimates in \cite[Proposition 5.3]{Hoshino2017} for $d=3$ and \cite[Theorem 2.1]{Matsuda2020} for $d=2$. Specifically, let $(X(t))_{t\geq 0}$ be a $\mathscr{S}^{\prime}(\mathbb{T}^d)$-valued random process, satisfying that for any fixed $t\geq0$ and any $\phi\in C^{\infty}(\mathbb{T}^d)$, there exists some function $f_{t,\phi}\in L^2(Z^{k}\times Z^{l}) $ with some $k,l\in\mathbb{N}$ such that
	\begin{equation}\label{good process}
		\left\langle X(t,\cdot), \phi \right\rangle = I_{k,l}(f_{t,\phi} ) .
	\end{equation}
	Using hypercontractivity property (Lemma \ref{lem: hypercontractivity}), we prove Theorem \ref{Prop: regularity for general distribution}. 
	
	\begin{proof}[Proof of Theorem \ref{Prop: regularity for general distribution}]
		By the definition of Littlewood-Paley block \eqref{eq: def of Littlewood-Paley block} and assumptions \eqref{eq: 2nd est for diff freq for general distribution}, \eqref{eq: 2nd est for general distribution}, for any $j\geq-1$, any $t\geq0$ and any $0<\epsilon<1$,
		\begin{align*}
			\mathrm{E}\left[ \left|  \left( \delta_j X\right) (t,x)\right| ^{2}\right]
			&=\sum _{m_1,m_2\in \mathbb{Z}^d} \chi_{j}(m_1)\chi_{j}(m_2)	\mathrm{E}\left[ \hat{X} (t,m_1)\overline{\hat{X} (t,m_2)}\right] \frac{1}{(2\pi)^d}e^{\i (m_1-m_2)\cdot x}\\
			&\lesssim \sum _{m\in \mathbb{Z}^d} |\chi_{j}(m)|^2 \frac{1}{(1+|m|^2)^{1-\epsilon}}\lesssim 2^{jd}2^{-2j(1-\epsilon)} =2^{2j(d/2-1+\epsilon)}.
		\end{align*}
		Since for fixed $t\geq0$ and $x\in\mathbb{T}^d$, $\left( \delta_j X\right) (t,x)= 	\left\langle X(t,\cdot), \eta_j(x-\cdot) \right\rangle= I_{k,l}(f_{t,\eta_j(x-\cdot)} )$, then for any $p\geq2$, by hypercontractivity property (see Lemma \ref{lem: hypercontractivity}), 
		\begin{equation*}
			\mathrm{E}\left[ \left\| \left( \delta_j X\right) (t,\cdot)\right\| _{L^p}^{p}\right]
			\leq \sup_{x\in\mathbb{T}^d} \mathrm{E}\left[ \left|  \left( \delta_j X\right) (t,x)\right| ^{p}\right]\lesssim \sup_{x\in\mathbb{T}^d} \left( \mathrm{E}\left[ \left|  \left( \delta_j X\right) (t,x)\right| ^{2}\right]\right) ^{\frac p2}\lesssim 2^{jp(d/2-1+\epsilon)}.
		\end{equation*}
		For any $\alpha>0$, we now take $p$ sufficiently large and $\epsilon$ sufficiently small such that $p\geq2$ and $-\alpha+\frac{d}{p}+\epsilon<0$, then by Proposition \ref{regularity embedding},
		\begin{align*}
			\mathrm{E}\left[ \left\|  X(t)\right\| _{{\mathcal{C}^{1-d/2-\alpha}}}^{p}\right]
			& \lesssim 	\mathrm{E}\left[ \left\|  X(t)\right\| _{{\mathcal{B}_{p,p}^{1-d/2-\alpha+d/p}}}^{p}\right]=	\sum_{j \geq-1} 2^{(1-\frac{d}{2}-\alpha+\frac dp)jp} \mathrm{E}\left[ \left\| \left( \delta_j  X\right) (t,\cdot)\right\| _{L^p}^{p}\right]\\
			&\lesssim\sum_{j \geq-1} 2^{(1-\frac{d}{2}-\alpha+\frac dp+\frac{d}{2}-1+\epsilon )jp} = \sum_{j \geq-1} 2^{(-\alpha+\frac dp+\epsilon )jp} \lesssim1.
		\end{align*}
		Then we get that for any fixed $t\in [0,\infty)$, $X(t)\in \mathcal{C}^{1-d/2-\alpha} $ almost surely for any $\alpha>0$.
		
		By \eqref{eq: 2nd est for diff freq for general distribution} and \eqref{eq: time diff est for general distribution}, we get that for any $j\geq-1$, any $ t_1,t_2\geq0$ and any $0<\epsilon<1$,
		\begin{align*}
			&\mathrm{E}\left[ \left| 	\left( \delta_j \left( X(t_2)-X(t_1)\right)\right) (x)\right| ^2  \right] \\
			=&\,\sum _{m_1,m_2\in \mathbb{Z}^d} \chi_{j}(m_1)\chi_{j}(m_2)\frac{1}{(2\pi)^d}e^{\i (m_1-m_2)\cdot x}\\ & \qquad \qquad \cdot	\mathrm{E}\Big[ \Big( \hat{X} (t_2,m_1)- \hat{X} (t_1,m_1)\Big) \Big( \overline{\hat{X} (t_2,m_2)}- \overline{\hat{X} (t_1,m_2)}\Big) \Big] \\
			\lesssim&\, \sum _{m\in \mathbb{Z}^d} |\chi_{j}(m)|^2	\mathrm{E}\left[\left|  \hat{X} (t_2,m)- \hat{X} (t_1,m)\right| ^2 \right] 
			\lesssim |t_2-t_1|^{\epsilon} \sum _{m\in \mathbb{Z}^d} |\chi_{j}(m)|^2 \frac{1}{(1+|m|^2)^{1-\epsilon}}\\
			\lesssim&\, |t_2-t_1|^{\epsilon} 2^{ jd}2^{-2 j(1-\epsilon)}= |t_2-t_1|^{\epsilon}2^{2j(d/2-1+\epsilon)}.
		\end{align*}
		Since for fixed $ t_1,t_2\geq0$ and $x\in\mathbb{T}^d$, $$	\left( \delta_j \left( X(t_2)-X(t_1)\right)\right) (x)=\left\langle X(t_2,\cdot)-X(t_1,\cdot), \eta_j(x-\cdot) \right\rangle= I_{k,l}(f_{t_2,\eta_j(x-\cdot)}- f_{t_1,\eta_j(x-\cdot)}),$$ then for any $p\geq2$, by hypercontractivity property (see Lemma \ref{lem: hypercontractivity}), 
		\begin{align*}
			\mathrm{E}\left[ \left\| \delta_j \left( X(t_2)-X(t_1)\right) \right\| _{L^p}^{p}\right]
			&\leq \sup_{x\in\mathbb{T}^d} \mathrm{E}\left[ \left|  	\left( \delta_j \left( X(t_2)-X(t_1)\right)\right) (x) \right| ^{p}\right]\\
			&\lesssim \sup_{x\in\mathbb{T}^d} \left( \mathrm{E}\left[ \left|  	\left( \delta_j \left( X(t_2)-X(t_1)\right)\right) (x)\right| ^{2}\right]\right) ^{\frac p2}\\
			&\lesssim |t_2-t_1|^{\frac{\epsilon p}{2}}2^{jp (d/2-1+\epsilon) }.
		\end{align*}
		For any $\alpha>0$, we now take $p$ sufficiently large and $\epsilon$ sufficiently small such that $p\geq2$, $\epsilon p>2$ and $-\alpha+\frac{d}{p}+\epsilon<0$, then by Proposition \ref{regularity embedding},
		\begin{align*}
			\mathrm{E}\left[ \left\| X(t_2)-X(t_1)\right\| _{{\mathcal{C}^{1-d/2-\alpha}}}^{p}\right]
			& \lesssim \mathrm{E}\left[ \left\| X(t_2)-X(t_1)\right\| _{{\mathcal{B}_{p,p}^{1-d/2-\alpha+d/p}}}^{p}\right]\\
			&=\sum_{j \geq-1} 2^{(1-d/2-\alpha+d/p)jp} \mathrm{E}\left[ \left\| \delta_j \left( X(t_2)-X(t_1)\right) \right\| _{L^p}^{p}\right]\\
			&\lesssim |t_2-t_1|^{\frac{\epsilon p}{2}}\sum_{j \geq-1} 2^{(-\alpha+\frac dp+\epsilon)jp}\lesssim |t_2-t_1|^{\frac{\epsilon p}{2}}.
		\end{align*}
		By Kolmogorov's continuity criterion (\cite[Theorem 3.3]{DPZ2014}), $X$ has a continuous version, which is still denoted by $X$, belongs to $C\left( [0,\infty);\mathcal{C}^{1-d/2-\alpha}\right) $ for any $\alpha>0$ and satisfies \eqref{eq: regularity for general distribution} and \eqref{eq: time regularity for general distribution}.
	\end{proof}

	\subsection{Regularity of Wick product}
	Utilizing the Da Prato-Debussche method (see \cite{DaPrato2003}), we first focus on the linearization of \eqref{SCGL}, namely the stochastic heat equation with dispersion, with zero initial data at time $a\in (-\infty, \infty)$,  
	\begin{equation}\label{heat equation}
		\begin{cases}
			\partial_t Z_a=\left[ (\i +\mu)\Delta-1\right]Z_a+ \xi, & t>a,x\in \mathbb{T}^2,\\
			Z_a(a,\cdot)=0,
		\end{cases}
	\end{equation}
	and derive the regularity properties of both the solution itself and its Wick products by using Theorem \ref{Prop: regularity for general distribution}.
	
	We set $A:=(\i+\mu)\Delta-1$ and let $\left\lbrace P_t=e^{tA}\right\rbrace _{t\geq0}$ be the semigroup with the generator $A$. Due to the roughness of $\xi$, we know that $\left( Z(t,\cdot)\right) _{t\geq a}$ is a distribution-valued stochastic process and by Duhamel's principle, for $\phi\in C^{\infty}(\mathbb{T}^2)$ and fixed $t\geq a$,
	\begin{align}
		\left\langle Z_a(t,\cdot),\phi\right\rangle &=\int_{a}^{t}\int_{\mathbb{T}^2}\left[ e^{(t-s)A}\xi(s,\cdot)\right] (x)\phi(x)\mathrm{d}x\mathrm{d}s \label{solution of heat equation}\\
		&=\int_{a}^{t}\int_{ \mathbb{T}^2}\left\langle H(t-s,\cdot-y),\phi\right\rangle  \xi(\mathrm{d}s,\mathrm{d}y),\nonumber
	\end{align}
	where $H(t,z)$ is the heat kernel defined as
	\begin{equation*}
		H(t,z)=\sum_{k \in \mathbb{Z}^2}\frac{1}{2\pi} e^{-\left( 1+(\i+\mu)|k|^2\right)t} \frac{1}{2\pi}e^{ \i k\cdot z}, \quad t\in\mathbb{R}\setminus\{0\},\,z\in \mathbb{T}^2.
	\end{equation*}
	For convenience, we extend the definition \eqref{solution of heat equation} to $a=-\infty$ and introduce the stationary solution $\left( Z_{-\infty}(t,\cdot)\right) _{t>-\infty}$ to \eqref{heat equation} as follows. For $\phi\in C^{\infty}(\mathbb{T}^2)$ and fixed $t>-\infty$,
	\begin{align}\label{stationary solution of heat equation}
		\left\langle Z_{-\infty}(t,\cdot),\phi\right\rangle
		=\int_{-\infty}^{t}\int_{ \mathbb{T}^2}\left\langle H(t-s,\cdot-y),\phi\right\rangle  \xi(\mathrm{d}s,\mathrm{d}y).
	\end{align}
	Now fix $t>-\infty$. We set $Z_{-\infty}^{:0,0:}(t,x)\equiv1$ for any $x\in\mathbb{T}^2$. For $k,l\in\mathbb{N}$ with $k+l>0$, we define $(k,l)$-th Wick product of $Z_{-\infty}$ as the complex $(k,l)$-th Wiener-It\^o integral with respect to $\xi$
	\begin{align}\label{definition of wick product}
		\left\langle Z_{-\infty}^{:k,l:}(t,\cdot),\phi\right\rangle :=I_{k,l}(f^{k,l}_{-\infty,t,\phi}),\quad \phi\in C^{\infty}(\mathbb{T}^2),
	\end{align}
	where $f^{k,l}_{-\infty,t,\phi}\in L^2((\mathbb{R}\times\mathbb{T}^2)^k \times (\mathbb{R}\times\mathbb{T}^2)^l)$ (see \eqref{variance of wick product}) is given by
	\begin{align}
		&f^{k,l}_{-\infty,t,\phi}(s_1,y_1,\ldots,s_k,y_k;s_{k+1},y_{k+1},\ldots,s_{k+l},y_{k+l}) \label{f k,l -infty}\\
		=&\,\left\langle \prod_{j=1}^{k}\mathbf{1}_{(-\infty,t]}(s_j) H(t-s_j,\cdot-y_j)\prod_{j=k+1}^{k+l}\mathbf{1}_{(-\infty,t]}(s_j) \overline{H}(t-s_j,\cdot-y_j),\phi\right\rangle.\nonumber
	\end{align}
	For $j=(j_1,j_2)\in \mathbb{Z}^2$, we define $|j|=\sqrt{j_1^2+j_2^2}$ and set $\rho_j  =1+\mu|j|^2$, $\theta_j=|j|^2$. From the calculation of \eqref{eq: 2nd est with test func} in the proof of Theorem \ref{regularity of stationary Z} below and \cite[Corollary C.3]{TW2018} (or \cite[Lemma 10.14]{Hairer2014}), we see that for any $0<\epsilon<1$,
	\begin{align}
		\left\| f^{k,l}_{-\infty,t,\phi}\right\|^2
		&= \frac{1}{(2\pi)^{2(k+l-1)}}\sum_{j\in \mathbb{Z}^2} |\hat{\phi}(-j)|^2 \left( \sum_{j_1+\cdots+ j_{k+l}=j}  \prod_{i=1}^{k+l}\frac{1}{2\rho_{j_i}}  \right) \label{variance of wick product}\\
		&\lesssim \sum_{j\in \mathbb{Z}^2} |\hat{\phi}(j)|^2 \frac{1}{(1+|j|^2)^{1-\epsilon}} \leq \left\|\phi \right\|^2  <\infty.\nonumber
	\end{align}
	Therefore, $f^{k,l}_{-\infty,t,\phi}\in L^2((\mathbb{R}\times\mathbb{T}^2)^k \times (\mathbb{R}\times\mathbb{T}^2)^l)$ and $\left\langle Z_{-\infty}^{:k,l:}(t,\cdot),\phi\right\rangle $ with $\phi\in C^{\infty}(\mathbb{T}^2)$ is a well-defined complex $(k,l)$-th Wiener-It\^o integral with respect to $\xi$. By \eqref{stationary solution of heat equation} and \eqref{definition of wick product}, we get that $Z_{-\infty}^{:1,0:}(t,x)=Z_{-\infty}(t,x)$ for $t> -\infty$ and $x\in\mathbb{T}^2$.

	Applying Theorem \ref{Prop: regularity for general distribution}, we now show that for every $k,l\in \mathbb{N}$, $t_0>-\infty$ and any $\alpha>0$, $\big( Z_{-\infty}^{:k,l:}(t_0+t,\cdot)\big) _{t\geq 0}$ is a $\mathcal{C}^{-\alpha}$-valued continuous random process.
	
	\begin{theorem}\label{regularity of stationary Z}
		Let $k,l\in \mathbb{N}$ and $t_0>-\infty$. Then for any $\alpha>0$, $\big( Z_{-\infty}^{:k,l:}(t_0+t,\cdot)\big) _{t\geq 0}\in C([0,\infty);\mathcal{C}^{-\alpha})$ almost surely. Moreover, for every $T>0$ and any $p\geq1$, there exist $0<\lambda<1$ and a constant $C$ depending on $\alpha$, $T$ and $p$ such that
		\begin{align}
			&\mathrm{E}\left[\sup_{0\leq t\leq T} \left\| Z_{-\infty}^{:k,l:}(t_0+t)\right\|^{p} _{\mathcal{C}^{-\alpha}} \right] \leq  C, \label{regularity}\\
			&\mathrm{E}\left[	\sup_{0\leq s<t\leq T}\frac{ \left\| Z_{-\infty}^{:k,l:}(t_0+t)-Z_{-\infty}^{:k,l:}(t_0+s)\right\|^{p} _{\mathcal{C}^{-\alpha}} }{(t-s)^{\lambda p}}\right] \leq  C. \nonumber
		\end{align} 
	\end{theorem}

	\begin{proof}
		The conclusions are obvious for $k=l=0$. Now we assume that $k,l\in\mathbb{N}$ with $k+l>0$. By Theorem \ref{Prop: regularity for general distribution}, it suffices to show that for any $s,t>-\infty$, $m,m_1,m_2\in \mathbb{Z}^2$ satisfying $m_1\neq m_2$ and any $0<\epsilon<1$,
		\begin{align}
			&\mathrm{E} \left[ \hat{Z}_{-\infty}^{:k,l:}(t,m_1)\overline{\hat{Z}_{-\infty}^{:k,l:}(s,m_2)}\right]=0,\label{eq: 2nd est for diff freq for Z}\\
			&\mathrm{E} \left[ \left| \hat{Z}_{-\infty}^{:k,l:}(t,m)\right| ^2\right]\lesssim \frac{1}{(1+|m|^2)^{1-\epsilon}},\label{eq: 2nd est for Z}\\
			&\mathrm{E}\left[\left|  \hat{Z}_{-\infty}^{:k,l:} (t,m)- \hat{Z}_{-\infty}^{:k,l:} (s,m)\right| ^2 \right]\lesssim |t-s|^{\epsilon}\frac{1}{(1+|m|^2)^{1-\epsilon}},\label{eq: time diff est for Z}
		\end{align}	
		which correspond to \eqref{eq: 2nd est for diff freq for general distribution}, \eqref{eq: 2nd est for general distribution} and \eqref{eq: time diff est for general distribution} respectively in Theorem \ref{Prop: regularity for general distribution}. For any $\phi_1,\phi_2\in C^{\infty}(\mathbb{T}^2)$, by isometry property \eqref{eq: general isometry property}, we obtain that (see \eqref{eq: 2nd est with test func in App} for the detail),
		\begin{align}
			&\mathrm{E} \left[ \left\langle Z_{-\infty}^{:k,l:}(t,\cdot),\phi_1\right\rangle \overline{\left\langle Z_{-\infty}^{:k,l:}(s,\cdot),\phi_2\right\rangle}\right]\label{eq: 2nd est with test func} \\
			=&\, \mathrm{E} \left[ I_{k,l}(f^{k,l}_{-\infty,t,\phi_1}) \overline{I_{k,l}(f^{k,l}_{-\infty,s,\phi_2}) } \right] = k!l!\left\langle f^{k,l}_{-\infty,t,\phi_1}, \overline{f^{k,l}_{-\infty,s,\phi_2} }\right\rangle \nonumber\\ 
			=&\, \frac{k!l!}{(2\pi)^{2(k+l-1)}}\sum_{j\in \mathbb{Z}^2} \hat{\phi_1}(-j) \overline{\hat{\phi_2}(-j)}  \Bigg( \sum_{j_1+\cdots+ j_{k+l}=j}\prod_{i=1}^{k+l}\frac{ e^{-\rho_{j_i}|t-s| -\i a_i\theta_{j_i} (t-s) }}{2\rho_{j_i}}  \Bigg) \nonumber,
		\end{align}
		where $a_i=1$ for $1\leq i\leq k$ and $a_i=-1$ for $k+1\leq i\leq k+l$.
		
		For $i=1,2$, we take $\phi_i(x)=\frac{1}{2\pi}e^{- \i m_i\cdot x}\in C^\infty(\mathbb{T}^2)$ with fixed $m_i\in\mathbb{Z}^2$, then
		\begin{align}
			&\mathrm{E} \left[ \hat{Z}_{-\infty}^{:k,l:}(t,m_1)\overline{\hat{Z}_{-\infty}^{:k,l:}(s,m_2)}\right]= \mathrm{E} \left[ \left\langle Z_{-\infty}^{:k,l:}(t,\cdot),\phi_1\right\rangle \overline{\left\langle Z_{-\infty}^{:k,l:}(s,\cdot),\phi_2\right\rangle}\right] \label{hatZ}\\
			=&\,\mathbf{1}_{\left\lbrace m_1=m_2\right\rbrace }\frac{k!l!}{(2\pi)^{2(k+l-1)}}\sum_{j_1+\cdots+j_{k+l}=m_1}\prod_{i=1}^{k+l}\frac{ e^{ -\rho_{j_i}|t-s|-\i a_i\theta_{j_i} (t-s) }}{2\rho_{j_i}} ,\nonumber
		\end{align}
		which implies \eqref{eq: 2nd est for diff freq for Z}. 
		
		Taking $s=t>-\infty$, $m_1=m_2=m\in \mathbb{Z}^2$ and using \cite[Corollary C.3]{TW2018}, we deduce that
		\begin{align}
			\mathrm{E} \left[ \left| \hat{Z}_{-\infty}^{:k,l:}(t,m)\right| ^2\right]
			&=\frac{k!l!}{(2\pi)^{2(k+l-1)}}	\sum_{j_1+\cdots+j_{k+l}=m}\prod_{i=1}^{k+l}\frac{ 1 }{2\rho_{j_i}} \lesssim \frac{1}{(1+|m|^2)^{1-\epsilon}} \label{hatZ for same time},
		\end{align}
		where $0< \epsilon<1$ is arbitrary. Then we obtain \eqref{eq: 2nd est for Z}. 
		
		Using \eqref{eq bound for A^{k,l}} and \cite[Corollary C.3]{TW2018}, we get that
		\begin{align}
			&\mathrm{E}\left[ \left|\hat{Z}_{-\infty}^{:k,l:}(t,m)-\hat{Z}_{-\infty}^{:k,l:}(s,m) \right| ^2\right] \label{hatZ diff time}
			\\=&\,\mathrm{E}\left[ \left|\hat{Z}_{-\infty}^{:k,l:}(t,m)\right| ^2+\left| \hat{Z}_{-\infty}^{:k,l:}(s,m) \right| ^2-\hat{Z}_{-\infty}^{:k,l:}(t,m)\overline{\hat{Z}_{-\infty}^{:k,l:}(s,m)}-\overline{\hat{Z}_{-\infty}^{:k,l:}(t,m)}\hat{Z}_{-\infty}^{:k,l:}(s,m)\right] \nonumber
			\\=&\, \frac{k!l!}{(2\pi)^{2(k+l-1)}}	\sum_{j_1+\cdots+j_{k+l}=m}  A^{k,l}_{s,t}(j_1,\ldots,j_{k+l})\nonumber \\\lesssim&\,  |t-s|^{\frac{\epsilon}{k+l}} \sum_{j_1+\cdots+j_{k+l}=m}\prod_{i=1}^{k+l}\frac{1}{(1+|j_i|^2)^{1-\frac{\epsilon}{k+l}}} \nonumber
			\\\lesssim &\, |t-s|^{\frac{\epsilon}{k+l}} \frac{1}{(1+|m|^2)^{1-\epsilon}}, \nonumber
		\end{align}
		where $0< \epsilon<1$ is arbitrary. Then we complete the proof of \eqref{eq: time diff est for Z}.
	\end{proof}

	By \eqref{solution of heat equation} and \eqref{stationary solution of heat equation}, we know that for any $t\geq a>-\infty$,
	\begin{equation*}
		Z_{a}(t)= Z_{-\infty}(t)- P_{t-a}Z_{-\infty}(a).
	\end{equation*}
	Then for $t> a>-\infty$, we set $Z_{a}^{:0,0:}(t,x)\equiv1$ for any $x\in\mathbb{T}^2$, and for $k,l\in\mathbb{N}$ with $k+l>0$, we define the $(k,l)$-th Wick product of $Z_{a}(t)$ as
	\begin{equation}\label{Z_{s,t}^{:k,l:}}
		Z_{a}^{:k,l:}(t)=\sum_{i=0}^{k}\sum_{j=0}^{l}\binom{k}{i}\binom{l}{j} \left(- P_{t-a}Z_{-\infty}(a)\right)  ^i \overline{\left(- P_{t-a}Z_{-\infty}(a)\right) }^j Z_{-\infty}^{:k-i,l-j:}(t).
	\end{equation}
	Similar definitions of Wick products for stochastic heat equation driven by real-valued space-time white noise can be found in \cite[Equation (5.42)]{MourratWeber2017} and \cite[Equation (2.7)]{TW2018}. Applying Theorem \ref{regularity of stationary Z}, we derive the regularity of $Z_{a}^{:k,l:}$ in the following theorem.
	\begin{theorem}\label{regularity of Z}
		Let $k,l\in\mathbb{N}$ satisfy $k+l>0$ and $a>-\infty$. For every $T>0$, any $p\geq1$ and $\alpha, \alpha^{\prime}\in (0,1)$, there exist $\theta>0$ depending on $\alpha, \alpha^{\prime},k,l$ and $C$ depending on $\alpha, \alpha^{\prime}, k,l,p,T$ such that 
		\begin{equation}\label{regularity with initial data}
			\mathrm{E}\left[ \sup_{0< t\leq T}t^{(k+l-1)\alpha^{\prime}p}\left\|Z_{a}^{:k,l:}(a+t)  \right\| _{\mathcal{C}^{-\alpha}}^{p}\right] \leq CT^{p\theta}.
		\end{equation}
	\end{theorem}
	\begin{proof}
		Since for any $s,t>-\infty$, $Z_{-\infty}(s)$ and $Z_{-\infty}(t)$ are Gaussian variables with the same variance (see \eqref{variance of wick product} for the case $(k,l)=(1,0)$), $Z_{-\infty}(s)$ and $Z_{-\infty}(t)$ have the same distribution. Combining this fact with \eqref{Z_{s,t}^{:k,l:}}, we know that for any $a,b>-\infty$ and $t>0$, $Z_{a}^{:k,l:}(a+t) $ and $Z_{b}^{:k,l:}(b+t) $ have the same distribution. Therefore, without loss of generality, we prove \eqref{regularity with initial data} for $a=0$. 
		
		Let $\bar{\alpha}=\alpha\wedge \frac{2}{3}\alpha^{\prime}$. We first consider the case $k+l=1$. By Proposition \ref{1-heat kernel smoothing},
		\begin{align*}
			\left\|Z_{0}^{:k,l:} (t) \right\| _{\mathcal{C}^{-\alpha}}&\leq \left\|Z_{-\infty}^{:k,l:}(t)  - Z_{-\infty}^{:k,l:}(0) \right\| _{\mathcal{C}^{-\alpha}}+ \left\|(1-P_t) Z_{-\infty}^{:k,l:} (0)\right\| _{\mathcal{C}^{-\alpha}}\\&\lesssim \left\|Z_{-\infty}^{:k,l:}(t)  - Z_{-\infty}^{:k,l:} (0)\right\| _{\mathcal{C}^{-\alpha}}+ t^{\frac{\alpha-\bar{\alpha}}{2}}\left\| Z_{-\infty}^{:k,l:}(0) \right\| _{\mathcal{C}^{-\bar{\alpha}}}.
		\end{align*}	
		Using Theorem \ref{regularity of stationary Z},
		\begin{align*}
			\mathrm{E}\left[ \sup_{0< t\leq T}\left\|Z_{0}^{:k,l:}(t)  \right\| _{\mathcal{C}^{-\alpha}}^{p}\right] &\lesssim  	\mathrm{E}\left[ \sup_{0< t\leq T}\left\|Z_{-\infty}^{:k,l:}(t)  - Z_{-\infty}^{:k,l:}(0) \right\|^p _{\mathcal{C}^{-\alpha}}\right] + T^{\frac{\alpha-\bar{\alpha}}{2}p}	\mathrm{E}\left[\left\| Z_{-\infty}^{:k,l:}(0) \right\|^p _{\mathcal{C}^{-\bar{\alpha}}}\right] \\& \lesssim T^{\lambda p}+T^{\frac{\alpha-\bar{\alpha}}{2}p}.
		\end{align*}
		Then there exists $\theta>0$ depending on $\alpha, \alpha^{\prime}$ such that \eqref{regularity with initial data} is valid for $k+l=1$.
		
		Now we consider the case $k+l>1$. By Proposition \ref{embedding}, \ref{heat kernel smoothing} and \ref{multiplicative structure}, we get that
		\begin{align*}
			&\left\|P_t Z_{-\infty}(0)^k\overline{P_t Z_{-\infty}(0)}^l\right\|_{\mathcal{C}^{-\alpha}} \lesssim	\left\|P_t Z_{-\infty}(0)^k\overline{P_t Z_{-\infty}(0)}^l\right\|_{\mathcal{C}^{-\bar{\alpha}}}\\\lesssim&\,\left\|P_t Z_{-\infty}(0)\right\|_{\mathcal{C}^{2 \bar{\alpha}}}^{k+l-1}\left\|P_t Z_{-\infty}(0)\right\|_{\mathcal{C}^{-\bar{\alpha}}} \lesssim t^{- \frac{3}{2} \bar{\alpha}(k+l-1)}\left\|Z_{-\infty}(0)\right\|_{\mathcal{C}^{-\bar{\alpha}}}^{k+l} .
		\end{align*}
		Similarly, for $(i,j)\neq (k,l)$, we have that
		\begin{align*}
			&\left\|P_t Z_{-\infty}(0)^i\overline{P_t Z_{-\infty}(0)}^jZ_{-\infty}^{:k-i,l-j:}(t)         \right\|_{\mathcal{C}^{-\alpha}} 
			\lesssim t^{- \frac{3}{2} \bar{\alpha}(i+j)}\left\|Z_{-\infty}(0)\right\|_{\mathcal{C}^{-\bar{\alpha}}}^{i+j} \left\|Z_{-\infty}^{:k-i,l-j:}(t)\right\|_{\mathcal{C}^{-\bar{\alpha}}}.
		\end{align*}
		Then by triangle inequality and Cauchy-Schwarz inequality, 
		\begin{align*}
			&\mathrm{E}\left[ \sup_{0< t\leq T}t^{(k+l-1)\alpha^{\prime}p}\left\|Z_{0}^{:k,l:} (t) \right\| _{\mathcal{C}^{-\alpha}}^{p}\right] \\
			\lesssim &\, T^{(\alpha^{\prime}- \frac{3}{2}\bar{\alpha})(k+l-1)p}	\mathrm{E}\left[\left\|Z_{-\infty}(0)\right\|_{\mathcal{C}^{-\bar{\alpha}}}^{(k+l)p} \right] +\sum_{0\leq i+j< k+l} \Bigg(  T^{\left((k+l-1)\alpha^{\prime}- \frac{3}{2} \bar{\alpha}(i+j)\right)p}	\\& \cdot \left( \mathrm{E}\left[\left\|Z_{-\infty}(0)\right\|_{\mathcal{C}^{-\bar{\alpha}}}^{2(i+j)p}\right] \right) ^{\frac12} \left( \mathrm{E}\left[\sup_{0< t\leq T}\left\|Z_{-\infty}^{:k-i,l-j:}(t)\right\|^{2p}_{\mathcal{C}^{-\bar{\alpha}}}\right]\right)^{\frac12}\Bigg) .
		\end{align*}
		Combining with Theorem \ref{regularity of stationary Z}, we derive \eqref{regularity with initial data} for $k+l>1$.
	\end{proof}
	
	In Section \ref{Global well-posedness}, we frequently use the properties of 
	\begin{align}\label{set of wick product}
		\left\lbrace Z_a^{:k,l:}: 0\leq k \leq m+1, 0\leq l\leq m ,k+l>0 \right\rbrace .
	\end{align}		
	From \eqref{reverse complex conjugate} and \eqref{definition of wick product}, we get that for any fixed $t>-\infty$ and any $k,l\in \mathbb{N}$,
	\begin{equation*}
		\left\langle \overline{Z_{-\infty}^{:k,l:}(t,\cdot)},\phi\right\rangle=\left\langle Z_{-\infty}^{:l,k:}(t,\cdot),\phi\right\rangle, \quad \forall\, \phi \in C^{\infty}(\mathbb{T}^2),
	\end{equation*}
	which implies that 
	\begin{equation}\label{conjugate of wick product}
		Z_{-\infty}^{:l,k:}(t,x)=\overline{	Z_{-\infty}^{:k,l:}(t,x)},\quad t>-\infty,\, x\in \mathbb{T}^2 , \, k,l\in \mathbb{N}.
	\end{equation}
	Combining \eqref{conjugate of wick product} with \eqref{Z_{s,t}^{:k,l:}}, we have 
	\begin{equation*}
		Z_{a}^{:l,k:}(t,x)=\overline{Z_{a}^{:k,l:}(t,x)}, \quad t> a >-\infty,\, x\in \mathbb{T}^2 , \, k,l\in \mathbb{N}.
	\end{equation*}
	Then for the set \eqref{set of wick product}, it suffices to consider the case $k\geq l$. Therefore, we define a set 
	\begin{equation}\label{set L}
		L:=\left\lbrace (i,j)\in \mathbb{N}^2: 0\leq i\leq m+1, 0\leq j\leq m, i\geq j, i+j>0\right\rbrace,
	\end{equation}
	with cardinality $M:=\frac{m(m+1)}{2}+2m+1$. By \eqref{Z_{s,t}^{:k,l:}}, Theorem \ref{regularity of stationary Z} and \ref{regularity of Z}, we know that for any $0<\alpha<1$,  
	\begin{equation}\label{space of Z}
		\underline{Z}_a:=\left\lbrace Z_a^{:i,j:}: (i,j)\in L\right\rbrace \in C([a,\infty);\mathcal{C}^{-\alpha})\times \left( C(( a,\infty);\mathcal{C}^{-\alpha})\right) ^{M-1},
	\end{equation}
	for almost surely $\omega\in\Omega$, and for every $T>0$, $0<\alpha^{\prime}<1$, there exists a constant $C$ depending on $\omega$, $\alpha, \alpha^{\prime}$ and $T$ such that 
	\begin{equation}\label{norm of Z}
		\left\| \underline{Z}_a\right\|_{\alpha,\alpha^{\prime},T}:=1\vee \sup_{(i,j)\in L} \sup_{0< t\leq T } t^{(i+j-1)\alpha^{\prime}} \left\|  Z_a^{:i,j:}(a+t)\right\| _{\mathcal{C}^{-\alpha}} \leq C  .
	\end{equation}
	As a random variable, $\mathrm{E}(C^p)<\infty$ for any $p\geq 1$. We sometimes omit the subscript $a$ in $Z_a$, $Z_a^{:k,l:}$, $\underline{Z}_a$ and $\left\| \underline{Z}_a\right\|_{\alpha,\alpha^{\prime},T}$ when $a=0$.

	For $\gamma>0$ and $K>0$, we define the set 
	\begin{equation}\label{E_{K,gamma}}
		E_{K,\gamma}:=\left\lbrace \left\| \underline{Z}\right\|_{\alpha,\alpha^{\prime},1}\leq K, \sum_{i=0}^{m+1}  \sum_{j=0}^{m}\int_{1}^T\left\| Z^{:i,j:}(s) \right\| _{\mathcal{C}^{-\alpha}}^\gamma \mathrm{d} s\leq K(1+T), \,\forall\, T\geq1\right\rbrace .
	\end{equation}
	In the following lemma, we estimate $ \underline{Z}$ and prove that for every $\gamma>0$, there exists $K>0$ such that $\mathrm{P}(E_{K,\gamma})>0$ in the spirit of \cite[Lemma 3.6]{RZZ2017}. We will show that the remainder \eqref{remainder} of \eqref{SCGL} converges to its auxiliary system \eqref{auxiliary system} on $E_{K,\gamma}$, and thus the generalized coupling we constructed (see \eqref{generalized coupling}) has positive measure on the diagonal at infinity $D$ (see \eqref{diagonal at infinity}) in the proof of Theorem \ref{ergodicity of u}. Consequently, the ergodicity of \eqref{SCGL} follows by using \cite[Corollary 2.2]{HairerMattinglyScheutzow2011} or Theorem \ref{HM Corollary 2.2}.
	\begin{lemma}\label{lem setE}
		For every $\gamma>0$ and $\varepsilon>0$, there exists $K>0$ such that $\mathrm{P}(E_{K,\gamma})\geq 1-\varepsilon$.
	\end{lemma}
	
	\begin{proof}
		By \cite[Theorem 3.3.1]{DZ1996}, we have that for every $q>0$, there exists a random variable $X\in L^2(\Omega,\mathcal{F},\mathrm{P})$ such that for $\mathrm{P}$-almost surely $\omega\in \Omega$,
		\begin{equation*}
			Y_{T}:=\frac{1}{T}\sum_{i=0}^{m+1}  \sum_{j=0}^{m}\int_{0}^{T}   \left\| Z_{-\infty}^{:i,j:}(s) \right\| _{\mathcal{C}^{-\alpha}}^q \mathrm{d}s\rightarrow X,\quad T\rightarrow\infty.
		\end{equation*}
		This implies that for every $\epsilon>0$, there exists $\Omega_1\subset\Omega$ such that $\mathrm{P}(\Omega_1)<\epsilon/2$ and  as $T\rightarrow\infty$, $\sup_{\omega\notin\Omega_1 }\left|	Y_{T}(\omega)-X(\omega)  \right| \rightarrow0$.
		Then there exists $T_0$ independent of $\omega$ such that for all $T \geq T_0$, $ Y_{T}(\omega)\leq X(\omega)+1$ for any $ \omega\notin \Omega_1$. Since $X \in L^2(\Omega,\mathcal{F},\mathrm{P})$, we get that for $\epsilon>0$, there exists $K_{\epsilon}=\sqrt{2\mathrm{E}X^2}  /\sqrt{\epsilon}> 0$ such that by Markov's inequality,
		\begin{equation*}
			\mathrm{P}(X>K_\epsilon)\leq \frac{\mathrm{E}X^2 }{K_\epsilon^2}=\frac{\epsilon}{2}.
		\end{equation*}
		On $\Omega_1^{c}\cap\left\lbrace X\leq K_{\epsilon} \right\rbrace $, we have that $	Y_{T}(\omega)\leq X(\omega)+1\leq K_{\epsilon}+1$ for any $T \geq T_0$. Therefore, let $K_1=K_\epsilon+1$, we get that for any $\gamma>0$,
		\begin{align*}
			&\mathrm{P}\Bigg(\sum_{i=0}^{m+1}  \sum_{j=0}^{m}\int_{0}^{T} \left\| Z_{-\infty}^{:i,j:} (s)\right\| _{\mathcal{C}^{-\alpha}}^{2\gamma}  \mathrm{d}s\leq K_1T, \,\forall\, T\geq T_0\Bigg)\\\geq&\, \mathrm{P}\left(\Omega_1^{c}\cap\left\lbrace X\leq K_{\epsilon} \right\rbrace\right) \geq \mathrm{P}\left(\Omega_1^{c}\right)+\mathrm{P}\left(\left\lbrace X\leq K_{\epsilon} \right\rbrace\right)-1\geq 1-\epsilon.
		\end{align*}
		Combining this with \eqref{regularity}, we know that there exists $K_2 > 0$ such that
		\begin{equation}\label{tech 1}
			\mathrm{P}\Bigg(\sum_{i=0}^{m+1}  \sum_{j=0}^{m}\int_{0}^{T} \left\| Z_{-\infty}^{:i,j:}(s) \right\| _{\mathcal{C}^{-\alpha}}^{2\gamma}  \mathrm{d}s\leq K_2(T+1), \forall\, T\geq 0\Bigg)\geq1-\epsilon.
		\end{equation}
		Using \eqref{Z_{s,t}^{:k,l:}}, Proposition \ref{heat kernel smoothing} and \ref{multiplicative structure}, and the fact that $s>1$, we have that for any $\gamma>0$,
		\begin{equation}\label{tech 2}
			\sum_{i=0}^{m+1}  \sum_{j=0}^{m} \left\| Z^{:i,j:}(s) \right\| _{\mathcal{C}^{-\alpha}}^\gamma \lesssim  \sum_{i=0}^{m+1}  \sum_{j=0}^{m} \left( \left\| Z_{-\infty}^{:i,j:}(0) \right\| _{\mathcal{C}^{-\alpha}}^{2\gamma(i+j)}+ \left\|Z_{-\infty}^{:i,j:}(s) \right\| _{\mathcal{C}^{-\alpha}}^{2\gamma}\right) .
		\end{equation}
		By \eqref{tech 1}, \eqref{tech 2} and \eqref{regularity}, we derive that there exists $K_3 > 0$ such that
		\begin{equation*}
			\mathrm{P}\Bigg(\sum_{i=0}^{m+1}  \sum_{j=0}^{m}\int_{1}^{T} \left\| Z^{:i,j:}(s) \right\| _{\mathcal{C}^{-\alpha}}^\gamma\mathrm{d}s\leq K_3(T+1), \forall\, T\geq 1\Bigg) \geq1-\epsilon.
		\end{equation*}
		By \eqref{norm of Z}, there exists $K_4$ such that $ \left\| \underline{Z}\right\|_{\alpha,\alpha^{\prime},1}\leq K_4$ almost surely. Then we finish the proof.
	\end{proof}

	\subsection{Approximation of Wick product}

	For any $\epsilon>0$, let $\phi_\epsilon(x)=\sum_{|m|<1/\epsilon}\frac{1}{2\pi}e^{\mathrm{i}m\cdot x}\in C^{\infty}(\mathbb{T}^2)$. We define a finite dimensional approximation of $Z_{a,\epsilon}(t,x)$ for $t\geq a>-\infty$ or $t> a=-\infty$, and $x\in\mathbb{T}^2$ as 
	\begin{equation}\label{Z_epsilon}
		Z_{a,\epsilon}(t,x):= \left\langle Z_a(t,\cdot),\phi_\epsilon(x-\cdot)\right\rangle=\int_{a}^{t}\int_{ \mathbb{T}^2}\left\langle H(t-s,\cdot-y),\phi_\epsilon(x-\cdot)\right\rangle  \xi(\mathrm{d}s,\mathrm{d}y).
	\end{equation}
	Following the calculation of \eqref{eq: 2nd est with test func}, we introduce the renormalization constant  
	\begin{equation*}
		c_\epsilon:=\mathrm{E}\left[ |Z_{-\infty,\epsilon}(t,x)|^2 \right]=\sum_{|j|<1/\epsilon}\frac{1}{2\rho_j}\sim \log(\epsilon^{-1}), \mbox{ as } \epsilon\rightarrow 0.
	\end{equation*}
	For $k,l\in\mathbb{N}$, we define the $(k,l)$-th Wick product of $Z_{a,\epsilon}$ as
	\begin{equation}\label{wick product of Z_epsilon}
		Z_{a,\epsilon}^{:k,l:}(t,x):=J_{k,l}(Z_{a,\epsilon}(t,x), c_\epsilon),
	\end{equation}
	where $J_{k,l}(z,\rho)$ is the $(k,l)$-th complex Hermite polynomial given by 
	\begin{equation}\label{complex Hermite polynomial}
		\exp\left\{ \lambda\overline{z}+\overline{\lambda}z-\rho|\lambda|^2\right\} =\sum_{k=0}^{\infty}\sum_{l=0}^{\infty}\frac{\overline{\lambda}^k\lambda^l}{k!l!}J_{k,l}(z,\rho),\quad z,\lambda\in\mathbb{C}, \, \rho>0.
	\end{equation}
	According to \cite[Equation (12.1)]{ito1952complex}, $J_{k,l}(z,\rho)$ has another equivalent definition as 
	\begin{equation*}
		J_{k,l}(z,\rho)= \sum_{r=0}^{k\wedge l} (-1)^r r!\binom{k}{r}\binom{l}{r} z^{k-r}\overline{z}^{l-r} \rho^r, \quad z\in\mathbb{C}, \, \rho>0.
	\end{equation*}

	To establish Proposition \ref{app of Wick power} below,  we show that the complex Hermite polynomial satisfies the following property (see Appendix \ref{Some technical estimation} for the proof). 
	\begin{lemma}\label{lem J_{k,l}(z_1+z_2)}
		For $k,l\in\mathbb{N}$, $z_1,z_2\in\mathbb{C}$ and $\rho>0$, 
		\begin{equation*}
			J_{k,l}(z_1+z_2,\rho)= \sum_{i=0}^{k} \sum_{j=0}^{l}\binom{k}{i}\binom{l}{j} J_{i,j}(z_1,\rho) z_2^{k-i}\overline{z_2}^{l-j}.
		\end{equation*}
	\end{lemma}
	By \eqref{Z_epsilon}, we know that for any $\epsilon>0$, $s\geq b>-\infty$ or  $s> b=-\infty$, and $r\geq0$,
	\begin{equation*}
		Z_{s,\epsilon}(s+r)= Z_{b,\epsilon}(s+r)- P_{r}Z_{b,\epsilon}(s).
	\end{equation*}
	Then combining Lemma \ref{lem J_{k,l}(z_1+z_2)} with \eqref{wick product of Z_epsilon}, we get that
	\begin{align}
		&Z_{s,\epsilon}^{:k,l:}(s+r) =J_{k,l}(Z_{s,\epsilon}(s+r), c_\epsilon)= J_{k,l}(Z_{b,\epsilon}(s+r)- P_{r}Z_{b,\epsilon}(s), c_\epsilon) \label{app time change of wick product}\\
		=&\,\sum_{i=0}^{k} \sum_{j=0}^{l}\binom{k}{i}\binom{l}{j} J_{k-i,l-j}(Z_{b,\epsilon}(s+r),c_\epsilon) (- P_{r}Z_{b,\epsilon}(s))^{i}(-\overline{ P_{r}Z_{b,\epsilon}(s)})^{j}\nonumber\\
		=&\,\sum_{i=0}^{k} \sum_{j=0}^{l}\binom{k}{i}\binom{l}{j}  (- P_{r}Z_{b,\epsilon}(s))^{i}(-\overline{ P_{r}Z_{b,\epsilon}(s)})^{j} Z^{:k-i,l-j:}_{b,\epsilon}(s+r).\nonumber
	\end{align}
	Utilizing \eqref{app time change of wick product}, we prove the following convergence results. 
	\begin{proposition}\label{app of Wick power}
		Let $k,l\in \mathbb{N}$, $T>0$, $p\geq1$ and $\alpha,\alpha^{\prime}\in(0,1)$. For every $a>-\infty$,
		\begin{align}
			&\lim\limits_{\epsilon\rightarrow 0}\mathrm{E}\left[\sup_{0\leq  t\leq T} \left\| Z_{-\infty,\epsilon}^{:k,l:}(a+t)- Z_{-\infty}^{:k,l:}(a+t)\right\|^{p} _{\mathcal{C}^{-\alpha}} \right] =0,\label{eq app of stationary wick power}\\
			&\lim\limits_{\epsilon\rightarrow 0}\mathrm{E}\left[\sup_{0< t\leq T} t^{(i+j-1)\alpha^{\prime}} \left\| Z_{a,\epsilon}^{:k,l:}(a+t)- Z_{a}^{:k,l:}(a+t)\right\|^{p} _{\mathcal{C}^{-\alpha}} \right] =0.\label{eq app of wick power}
		\end{align}
	\end{proposition}
	\begin{proof}
		Taking $s=a$, $r=t$ and $b=-\infty$ in \eqref{app time change of wick product} and recalling the definition of $Z_{a}^{:k,l:}(a+t)$ (see \eqref{Z_{s,t}^{:k,l:}}), it suffices to show \eqref{eq app of stationary wick power}. The conclusion is obvious for $k=l=0$. Now we assume that $k,l\in\mathbb{N}$ with $k+l>0$. We write $ Y_{-\infty,\epsilon}^{k,l}(t,x):= Z_{-\infty,\epsilon}^{:k,l:}(t,x)- Z_{-\infty}^{:k,l:}(t,x)$ for $t> -\infty$, and $x\in\mathbb{T}^2$. By \cite[Theorem 13.1]{ito1952complex}, we have that
		\begin{equation}\label{wick product of Z_epsilon as CMWII}
			Z_{-\infty,\epsilon}^{:k,l:}(t,x)=J_{k,l}(Z_{-\infty,\epsilon}(t,x), c_\epsilon)= I_{k,l}(h^{k,l}_{t,x,\epsilon}),
		\end{equation}
		where
		\begin{align*}
			&h^{k,l}_{t,x,\epsilon}(s_1,y_1,\ldots,s_k,y_k;s_{k+1},y_{k+1},\ldots,s_{k+l},y_{k+l})\\
			=&\,\prod_{j=1}^{k} \left\langle \mathbf{1}_{(-\infty,t]}(s_j) 	H(t-s_j,\cdot-y_j),\phi_\epsilon(x-\cdot)\right\rangle \\&\cdot \prod_{j=k+1}^{k+l} \left\langle \mathbf{1}_{(-\infty,t]}(s_j) \overline{H}(t-s_j,\cdot-y_j), \overline{\phi_\epsilon}(x-\cdot) \right\rangle .			
		\end{align*}
		Therefore, for fixed $t>-\infty$, and $x\in\mathbb{T}^2$, $Y_{-\infty,\epsilon}^{k,l}(t,x)$ is a complex $(k,l)$-th Wiener-It\^o integral with respect to $\xi$. Applying Theorem \ref{Prop: regularity for general distribution}, it suffices to show that for any $s,t> -\infty$, $m,m_1,m_2\in \mathbb{Z}^2$ satisfying $m_1\neq m_2$ and any $0<\eta<1$, there exists some $\gamma>0$ such that
		\begin{align}
			&\mathrm{E} \left[ \hat{Y}_{-\infty,\epsilon}^{k,l}(t,m_1)\overline{\hat{Y}_{-\infty,\epsilon}^{k,l}(s,m_2)}\right]=0,\label{eq: 2nd est for diff freq for Y}\\
			&\mathrm{E} \left[ \left| \hat{Y}_{-\infty,\epsilon}^{k,l}(t,m)\right| ^2\right]\lesssim \frac{\epsilon^{\gamma}}{(1+|m|^2)^{1-\eta}},\label{eq: 2nd est for Y}\\
			&\mathrm{E}\left[\left|  \hat{Y}_{-\infty,\epsilon}^{k,l} (t,m)- \hat{Y}_{-\infty,\epsilon}^{k,l} (s,m)\right| ^2 \right]\lesssim |t-s|^{\eta}\frac{\epsilon^{\gamma}}{(1+|m|^2)^{1-\eta}}.\label{eq: time diff est for Y}
		\end{align}	
		Using \eqref{definition of wick product}, \eqref{wick product of Z_epsilon as CMWII}, isometry property \eqref{eq: general isometry property} and the similar calculation in \eqref{eq: 2nd est with test func}, we get that
			\begin{align}
				&\mathrm{E}\left[ 	\hat{Z}_{-\infty,\epsilon}^{:k,l:}(t,m_1)\overline{\hat{Z}_{-\infty,\epsilon}^{:k,l:}(s,m_2)}\right]= (2\pi)^{k+l}\mathrm{E}\left[ 	\hat{Z}_{-\infty,\epsilon}^{:k,l:}(t,m_1)\overline{\hat{Z}_{-\infty}^{:k,l:}(s,m_2)}\right] \nonumber\\
				=&\, \mathbf{1}_{\left\lbrace m_1=m_2\right\rbrace }k!l!(2\pi)^{2}\sum_{j_1+\cdots+j_{k+l}=m_1, |j_i|<1/\epsilon}\prod_{i=1}^{k+l}\frac{ e^{-\rho_{j_i}|t-s| -\i a_i\theta_{j_i} (t-s) }}{2\rho_{j_i}}   \label{hatZepsilon}.
			\end{align}
			Combining with \eqref{hatZ}, we derive \eqref{eq: 2nd est for diff freq for Y}.
			
			Taking $s=t> -\infty$, $m_1=m_2=m\in \mathbb{Z}^2$ in \eqref{hatZepsilon}, using \eqref{hatZ for same time}, and applying \cite[Corollary C.3]{TW2018}, we deduce that for any $\eta\in(0,1)$,
			\begin{align*}
				\mathrm{E} \left[ \left| \hat{Y}_{-\infty,\epsilon}^{k,l}(t,m)\right| ^2\right]
				&\lesssim 	\Bigg( \sum_{j_1+\cdots+j_{k+l}=m}- \sum_{j_1+\cdots+j_{k+l}=m, |j_i|<1/\epsilon}\Bigg) \prod_{i=1}^{k+l}\frac{ 1}{2\rho_{j_i}}\\
				&\lesssim \frac{1}{(1+|m|^2)^{1-\eta}(1+1/\epsilon^2)^{\frac{\eta}{2}}} \lesssim \frac{\epsilon^{\eta}}{(1+|m|^2)^{1-\eta}}.
			\end{align*}
			Then we obtain \eqref{eq: 2nd est for Y}.
			
		By \eqref{hatZ}, \eqref{hatZepsilon}, a similar argument in \eqref{hatZ diff time}, Lemma \ref{lem bound for A^{k,l}} and \cite[Corollary C.3]{TW2018}, we get that for any $\eta\in(0,1)$,
			\begin{align*}
				&\mathrm{E}\left[\left|  \hat{Y}_{-\infty,\epsilon}^{k,l} (t,m)- \hat{Y}_{-\infty,\epsilon}^{k,l} (s,m)\right| ^2 \right]
				\\
				\lesssim&\,\Bigg( \sum_{j_1+\cdots+j_{k+l}=m}- \sum_{j_1+\cdots+j_{k+l}=m, |j_i|<1/\epsilon}\Bigg) A^{k,l}_{s,t}(j_1,\ldots,j_{k+l})\\
				\lesssim&\, |t-s|^{\frac{\eta}{2(k+l)}}\Bigg( \sum_{j_1+\cdots+j_{k+l}=m}- \sum_{j_1+\cdots+j_{k+l}=m, |j_i|<1/\epsilon}\Bigg)\prod_{i=1}^{k+l}\frac{1}{(1+|j_i|^2)^{1-\frac{\eta}{2(k+l)}}}\\
				\lesssim&\, |t-s|^{\frac{\eta}{2(k+l)}}\frac{1}{(1+|m|^2)^{1-\eta}(1+1/\epsilon^2)^{\frac{\eta}{2}}} \lesssim |t-s|^{\frac{\eta}{2(k+l)}}\frac{\epsilon^{\eta}}{(1+|m|^2)^{1-\eta}} .
			\end{align*}
		Then we finish the proof of \eqref{eq: time diff est for Y}.
		\end{proof}
		
		Taking $s=t$, $r=h$ and $b=a$ in \eqref{app time change of wick product}, letting $\epsilon\rightarrow 0$ and applying Proposition \ref{app of Wick power}, we obtain the following corollary.
		
		\begin{corollary}\label{Cor relation of wick product of diff initial moments}
			Let $k,l\in \mathbb{N}$, $t\geq a>-\infty$ or $t> a=-\infty$, and $h>0$. Then
			\begin{equation*}
				Z_t^{:k,l:}(t+h)= \sum_{i=0}^{k}\sum_{j=0}^{l}\binom{k}{i}\binom{l}{j}  	(- P_h Z_a(t))^{i}(-\overline{ P_h Z_a(t)})^{j} Z_a^{:k-i,l-j:}(t+h) .
			\end{equation*}
		\end{corollary}
		
		\begin{remark}
			In the definition of $Z_{a,\epsilon}^{:k,l:}(t)$ (see \eqref{wick product of Z_epsilon}), we choose the renormalization constant $c_{\epsilon}$, which is independent of $a$ and $t$. This ensures the validity of the corollary above, which is essential for establishing an estimate for $u=Z+v$ in Corollary \ref{priori extimate for u} by using the estimate  for $v$ in Proposition \ref{priori estimate}. This dissipative estimate for $u$ is in turn necessary for applying the Krylov-Bogoliubov theorem (see \cite[Corollary 3.1.2]{DZ1996}) to prove the existence of an invariant measure. In Proposition \ref{app of Wick power}, we show that  $Z_{a,\epsilon}^{:k,l:}(t)$ converges to $Z_{a}^{:k,l:}(t)$ as defined in \eqref{Z_{s,t}^{:k,l:}} for $t> a>-\infty$ in the sense of \eqref{eq app of wick power}.
			
			However, if we focus solely on the global well-posedness of \eqref{SCGL}, we can alternatively define the Wick products of $Z_a$ as follows. Let $t\geq a>-\infty$. We set $Z_a^{0,0}(t,x)\equiv1$ for any $x\in\mathbb{T}^2$. For $k,l\in\mathbb{N}$ with $k+l>0$, we define the $(k,l)$-th Wick product of $Z_a$ as the complex $(k,l)$-th Wiener-It\^o integral with respect to $\xi$
			\begin{align*}
				\left\langle Z_a^{k,l}(t,\cdot),\phi\right\rangle :=I_{k,l}(f^{k,l}_{a,t,\phi}),\quad \phi\in C^{\infty}(\mathbb{T}^2),
			\end{align*}
			where $f^{k,l}_{a,t,\phi}\in L^2((\mathbb{R}\times\mathbb{T}^2)^k \times (\mathbb{R}\times\mathbb{T}^d)^l)$ is defined by replacing $-\infty$ with $a$ in the definition of $f^{k,l}_{-\infty,t,\phi}$ (see \eqref{f k,l -infty}). We also define the approximate Wick product as
			\begin{align*}
				Z_{a,\epsilon}^{k,l}(t,x):=J_{k,l}(Z_{a,\epsilon}(t,x), c_{a,\epsilon}(t)),
			\end{align*}
			where 
			\begin{align*}
				c_{a,\epsilon}(t):= \mathrm{E}\left[ |Z_{a,\epsilon}^{k,l}(t,x)|^2 \right]=\sum_{|j|<1/\epsilon}\frac{1- e^{-2\rho_j(t-a)}}{2\rho_j}.
			\end{align*}
			In this formulation, both $Z_a^{k,l}$ and $Z_{a,\epsilon}^{k,l}$ are complex multiple Wiener-It\^o integrals. Analogously to Theorem \ref{regularity of stationary Z} and Proposition \ref{app of Wick power}, we can prove that for any $\alpha>0$,
			$Z_a^{:k,l:}(t,x)\in C([a,\infty);\mathcal{C}^{-\alpha})$ almost surely. Moreover, for every $T>a$ and any $p\geq1$, there exist $0<\lambda<1$ and a constant $C$ depending on $\alpha$, $T$ and $p$ such that
			\begin{align*}
				&\mathrm{E}\left[\sup_{a\leq t\leq T} \left\| Z_a^{:k,l:}(t)\right\|^{p} _{\mathcal{C}^{-\alpha}} \right] \leq  C, \\
				&\mathrm{E}\left[	\sup_{a\leq s<t\leq T}\frac{ \left\| Z_a^{:k,l:}(t)-Z_a^{:k,l:}(s)\right\|^{p} _{\mathcal{C}^{-\alpha}} }{(t-s)^{\lambda p}}\right] \leq  C,\\
				&\lim\limits_{\epsilon\rightarrow 0}\mathrm{E}\left[\sup_{a< t\leq T} \left\| Z_{a,\epsilon}^{:k,l:}(t)- Z_{a}^{:k,l:}(t)\right\|^{p} _{\mathcal{C}^{-\alpha}} \right] =0.
			\end{align*}
			Under this modified definition, the global well-posedness of \eqref{SCGL} still holds, as shown by arguments similar to those in Section \ref{Global well-posedness}, though the renormalization constant for \eqref{SCGL} now becomes time-dependent.
		\end{remark}

		\section{Global well-posedness}\label{Global well-posedness}
		
		\subsection{Local well-posedness for $m\geq1$}
		
		In this section, we let $T>a=0$, and $\alpha,\alpha^{\prime}>0$ sufficiently small (see \eqref{condition for alpha 1} and \eqref{condition for alpha 2} below). Then we fix $\omega\in\Omega$ such that \eqref{space of Z} and \eqref{norm of Z} hold for  $\alpha,\alpha^{\prime}$ and $T$, and consider a deterministic equation for the remainder. Inspired by the Da Prato-Debussche method (see \cite{DaPrato2003}), we say that $u$ solves the equation \eqref{SCGL} in the renormalized sense with initial data $u_0\in \mathcal{C}^{-\alpha_0}$ if $u=Z+v$, where $Z$, as defined in \eqref{solution of heat equation}, is the solution of the stochastic heat equation with dispersion \eqref{heat equation} with zero initial data at time $a=0$, and $v$ is a mild solution of \eqref{remainder} as defined below.
		\begin{definition}\label{def: mild solu v}
			Let $T>0$ and $\alpha_0$, $\beta$ and $\gamma$ satisfy 
			\begin{equation}\label{eq: para cond}
				\alpha_0>0,\quad \beta>0,\quad\frac{\alpha_0+\beta}{2}< \gamma<\frac{1}{2m+1}.
			\end{equation} 
			We say that a function $v$ is a mild solution of 
			\begin{equation}\label{remainder}
				\begin{cases}
					\partial_t v=\left[ (\i+\mu) \Delta-1\right]  v+	\Psi(v,\underline{Z}) , & t>0, x \in \mathbb{T}^2, \\
					v(0, \cdot)=u_0\in \mathcal{C}^{-\alpha_0},
				\end{cases}
			\end{equation}
			up to time $T$ if $v\in C\left( \left( 0,T\right]; \mathcal{C}^{\beta}\right) $ with the norm $\sup_{0< t\leq T}t^{\gamma}\left\|v(t) \right\|_{\mathcal{C}^{\beta}} <\infty$ and 
			\begin{equation*}
				v(t)=P_tu_0+\int_{0}^{t} P_{t-s} \Psi\left(v(s),\underline{Z}(s)\right)\mathrm{d}s,\quad 0\leq t\leq T.
			\end{equation*}
			Here, $	\Psi\left( v,\underline{Z}\right)$ is defined as
			\begin{equation}\label{Psi}
				\Psi\left( v,\underline{Z}\right) :=-\nu\sum_{i=0}^{m+1}\sum_{j=0}^{m}\binom{m+1}{i}\binom{m}{j}v^i\overline{v}^jZ^{:m+1-i,m-j:} +\tau(v+Z).
			\end{equation}
		\end{definition}
		
		Since $\underline{Z}\in C([0,\infty);\mathcal{C}^{-\alpha})\times \left( C((0,\infty);\mathcal{C}^{-\alpha})\right) ^{M-1}$ for any $\alpha>0$ (see \eqref{space of Z}), we assume that $\alpha<\beta$, then by Proposition \ref{multiplicative structure}, we know that $\Psi\left( v,\underline{Z}\right)$ is well-defined for $v\in C\left( \left( 0,T\right]; \mathcal{C}^{\beta}\right) $.

		\begin{theorem}[Local well-posedness]\label{Local existence and uniqueness}
			Let $T>0$ and $\alpha_0$, $\beta$, $\gamma$ satisfy \eqref{eq: para cond}. Then for any $u_0\in \mathcal{C}^{-\alpha_0}$, there exists $0<T^{*} \leq T $ depending on $ \left\| u_0\right\|_{\mathcal{C}^{-\alpha_0}} $ and $ \left\| \underline{Z}\right\|_{\alpha,\alpha^{\prime},T}$ such that \eqref{remainder} has a unique mild solution on $[0,T^{*}]$. 
		\end{theorem}
		
		\begin{proof}
			For any given $T > 0$, we define the solution map as
			\begin{equation*}
				\left(\mathcal{M}_Tv\right)(t)=P_tu_0+\int_{0}^{t} P_{t-s} \Psi\left(v(s),\underline{Z}(s)\right)\mathrm{d}s,\quad 0\leq t\leq T.
			\end{equation*} 
			We now prove that for a fixed $0 < T \leq 1$, if $R$ and $T^{*}$ satisfy that
			\begin{equation}\label{eq: cond for R and local T}
				R\geq C(\left\|u_0 \right\|_{\mathcal{C}^{-\alpha_0}}+ \left\| \underline{Z}\right\|_{\alpha,\alpha^{\prime},T}) \vee 1,\quad
				T^{*}\leq R^{-\frac{2m+1}{\gamma-\frac{\alpha_0+\beta}{2}}}\wedge T,
			\end{equation}
			then $\mathcal{M}_{T^{*}}$ is a contraction from $ \mathcal{B}_{T^{*}}(0,R):=\left\lbrace v\in C\left( \left( 0,T^{*}\right]; \mathcal{C}^{\beta}\right) : \sup_{0< t\leq T^{*}}t^{\gamma}\left\|v(t) \right\|_{\mathcal{C}^{\beta}} \right. $ $\left. \leq R \right\rbrace $ into itself. 
			
			Let $v\in \mathcal{B}_{T^{*}}(0,R)$. For fixed $s>0$ and any $s<t\leq T^{*}$, by Proposition \ref{heat kernel smoothing} and \ref{1-heat kernel smoothing},
			\begin{align}
				\left\|P_tu_0 -P_su_0\right\| _{\mathcal{C}^{\beta}} 
				&=	\left\|(P_{t-s}-1) P_su_0\right\| _{\mathcal{C}^{\beta}} \lesssim (t-s)^{\frac{\delta}{2}} \left\|P_su_0\right\| _{\mathcal{C}^{\beta+\delta}} \label{time cont of initial data} \\
				&\lesssim (t-s)^{\frac{\delta}{2}} s^{-\frac{\alpha_0+\beta+\delta}{2}} \left\|u_0 \right\|_{\mathcal{C}^{-\alpha_0}} \rightarrow0, \quad t\rightarrow s,\nonumber
			\end{align}
			where $0<\delta\leq2$. Therefore, $t\mapsto P_tu_0$
			is in $C\left( \left( 0,T^{*}\right]; \mathcal{C}^{\beta}\right)$. According to Proposition \ref{heat kernel smoothing} and the fact that $\frac{\alpha_0+\beta}{2}< \gamma$, we have that 
			\begin{equation*}
				\sup_{0< t\leq T^{*}} t^{\gamma}\left\|P_tu_0 \right\| _{\mathcal{C}^{\beta}}	\lesssim 	\sup_{0< t\leq T^{*}}  t^{\gamma-\frac{\alpha_0+\beta}{2}}	\left\|u_0 \right\|_{\mathcal{C}^{-\alpha_0}} 	\lesssim (T^{*})^{\gamma-\frac{\alpha_0+\beta}{2}}	\left\|u_0 \right\|_{\mathcal{C}^{-\alpha_0}} .
			\end{equation*}
			From a similar argument to \eqref{time cont of initial data}, we know $t\mapsto \int_{0}^{t} P_{t-s} \Psi\left(v(s),\underline{Z}(s)\right)\mathrm{d}s$ is in $C\left( \left( 0,T^{*}\right]; \mathcal{C}^{\beta}\right)$. We fix $\alpha, \alpha^{\prime}>0$ small enough such that 
			\begin{align}\label{condition for alpha 1}
				\alpha^{\prime}<\gamma,\quad \alpha<\beta\wedge \alpha_0,
			\end{align} 
			which implies that $\frac{\alpha+\beta}{2}+2m\gamma< \frac{\alpha_0+\beta}{2}+2m\gamma<(2m+1) \gamma<1$. For $0<t\leq T^{*}$, by Proposition \ref{heat kernel smoothing} and  \ref{multiplicative structure},
			\begin{align*}
				&\left\|\int_{0}^{t} P_{t-s} \Psi\left(v(s),\underline{Z}(s)\right)\mathrm{d}s\right\| _{\mathcal{C}^{\beta}}
				\lesssim\int_{0}^{t}\left\| P_{t-s} \Psi\left(v(s),\underline{Z}(s)\right)\right\| _{\mathcal{C}^{\beta}}\mathrm{d}s\\
				\lesssim&\,\int_{0}^{t } \left\|v(s)^{m+1}\overline{v(s)}^{m}\right\| _{\mathcal{C}^{\beta}} +\left\|v(s)\right\| _{\mathcal{C}^{\beta}}   \mathrm{d}s+ \int_{0}^{t } 
				(t-s)^{-\frac{\alpha+\beta}{2}}\left\| Z(s) \right\|_{\mathcal{C}^{-\alpha}}\mathrm{d}s\\&+ \sum_{0\leq i+j<2m+1}\int_{0}^{t } 
				(t-s)^{-\frac{\alpha+\beta}{2}}  \left\| v(s)^i\overline{v(s)}^jZ^{:m+1-i,m-j:}(s) \right\|_{\mathcal{C}^{-\alpha}}	
				\mathrm{d}s\\
				\lesssim&\,\left\| \underline{Z}\right\|_{\alpha,\alpha^{\prime},T}R^{2m+1}\Bigg(  \int_{0}^{t}s^{-(2m+1)\gamma}\mathrm{d}s \\& \qquad \qquad\qquad \qquad\quad+\sum_{0\leq i+j<2m+1} \int_{0}^{t}(t-s)^{-\frac{\alpha+\beta}{2}}s^{-(i+j)\gamma-(2m-(i+j))\alpha^{\prime}}\mathrm{d}s\Bigg)  \\
				\lesssim&\,  \left\| \underline{Z}\right\|_{\alpha,\alpha^{\prime},T} R^{2m+1}\left(t^{1-(2m+1)\gamma}+t^{1-\frac{\alpha+\beta}{2}-2m\gamma}\right)	\lesssim  \left\| \underline{Z}\right\|_{\alpha,\alpha^{\prime},T} R^{2m+1}t^{1-(2m+1)\gamma}.
			\end{align*}
			Using the fact that $\gamma-\frac{\alpha_0+\beta}{2}<1-2m\gamma$, we know that
			\begin{align*}
				\sup_{0< t\leq T^{*}} t^{\gamma} \left\|\int_{0}^{t} P_{t-s} \Psi\left(v(s),\underline{Z}(s)\right)\mathrm{d}s\right\| _{\mathcal{C}^{\beta}}&
				\lesssim \left\| \underline{Z}\right\|_{\alpha,\alpha^{\prime},T} R^{2m+1} (T^{*}) ^{1-2m\gamma} \\&\lesssim \left\| \underline{Z}\right\|_{\alpha,\alpha^{\prime},T} R^{2m+1} (T^{*})^{\gamma-\frac{\alpha_0+\beta}{2}}.
			\end{align*}
			Then we obtain that 
			\begin{align}\label{eq: invariant ball}
				\sup_{0< t\leq T^{*}} t^{\gamma}\left\|\left(\mathcal{M}_Tv\right)(t)\right\| _{\mathcal{C}^{\beta}}
				&\leq C( \left\|u_0 \right\|_{\mathcal{C}^{-\alpha_0}}+ \left\| \underline{Z}\right\|_{\alpha,\alpha^{\prime},T} )R^{2m+1}(T^{*})^{\gamma-\frac{\alpha_0+\beta}{2}}\leq R.
			\end{align}
			Therefore, $\mathcal{B}_{T^{*}}(0,R)$ is an invariant ball for the solution map $\mathcal{M}_{T^{*}}$. The contraction property of $\mathcal{M}_{T^{*}}$ can be proved by using the similar argument. By the fixed point theorem, there is a unique solution on $[0,T ^{*}]$ in $\mathcal{B}_{T^{*}}(0,R)$ if $R$ and $T^{*}$ satisfy \eqref{eq: cond for R and local T}.
			
			Next we show the uniqueness of solutions in $C\left( \left( 0,T^{*}\right]; \mathcal{C}^{\beta}\right)$ for fixed $T^{*}$ satisfying \eqref{eq: cond for R and local T}. Let $v_1,v_2\in C\left( \left( 0,T^{*}\right]; \mathcal{C}^{\beta}\right) $ be two solutions of \eqref{remainder}. By using the similar argument as above, we have that $v_1,v_2\in C\left( \left[ 0,T^{*}\right]; \mathcal{C}^{-\alpha_0}\right) $. We now take 
			\begin{equation*}
				R=C\left( \left( \left\| v_1\right\| _{C\left([0, T^{\ast}] ; \mathcal{C}^{-\alpha_0}\right)} \vee \left\| v_2\right\| _{C\left([0, T^{\ast}] ; \mathcal{C}^{-\alpha_0}\right)} \right) + \left\|\underline Z \right\| _{\alpha,\alpha^{\prime},T}\right)\vee 1<\infty,
			\end{equation*} 
			which satisfies \eqref{eq: cond for R and local T}. Following the similar calculations in \eqref{eq: invariant ball}, we can find $T^{\ast\ast}\leq T^{\ast}$ depending only on $\left\| v_1 \right\|_{C\left( \left( 0,T^{*}\right]; \mathcal{C}^{\beta}\right)}$ and $\left\| v_2 \right\|_{C\left( \left( 0,T^{*}\right]; \mathcal{C}^{\beta}\right)}$ such that $v_1,v_2\in B_{T^{\ast\ast}}(0,R)$. In fact, $T^{\ast\ast}\leq T^{\ast}$ is taken sufficiently small such that for $i=1,2$, 
			\begin{equation*}
				\left( \left\| v _i\right\| _{C\left( \left( 0,T^{*}\right]; \mathcal{C}^{\beta}\right)}\vee1\right)^{2m+1} (T^{**})^{\gamma-\frac{\alpha_0+\beta}{2}} \leq 1,
			\end{equation*}
			then
			\begin{align*}
				&\left\| v _i\right\| _{C\left( \left( 0,T^{**}\right]; \mathcal{C}^{\beta}\right)}=\left\| \mathcal{M}_{T^{\ast\ast}}v _i\right\| _{C\left( \left( 0,T^{**}\right]; \mathcal{C}^{\beta}\right)}\\
				\leq&\, C( \left\|u_0 \right\|_{\mathcal{C}^{-\alpha_0}}+ \left\| \underline{Z}\right\|_{\alpha,\alpha^{\prime},T} )\left( \left\| v _i\right\| _{C\left( \left( 0,T^{*}\right]; \mathcal{C}^{\beta}\right)}\vee1\right)^{2m+1} (T^{**})^{\gamma-\frac{\alpha_0+\beta}{2}}\leq R .
			\end{align*}
			Since $T^{\ast\ast}\leq T^{\ast}$, similar to \eqref{eq: invariant ball}, we get that $\mathcal{M}_{T^{\ast\ast}}$ is a contraction from $B_{T^{\ast\ast}}(0,R)$ into itself. By the fixed point theorem,  $v_1(t)=v_2(t)$ on $[0, T^{\ast\ast}]$. Then we consider $v_1(T^{\ast\ast}+t)$ and $v_2(T^{\ast\ast}+t)$ on $t\in [0, T^{\ast}-T^{\ast\ast}]$. Since
			\begin{equation*}
				\left\| v_1(T^{\ast\ast})\right\|_{\mathcal{C}^{-\alpha_0}}= \left\| v_2(T^{\ast\ast})\right\|_{\mathcal{C}^{-\alpha_0}} \leq \left\| v_1\right\| _{C\left([0, T^{\ast}] ; \mathcal{C}^{-\alpha_0}\right)} \vee \left\| v_2\right\| _{C\left([0, T^{\ast}] ; \mathcal{C}^{-\alpha_0}\right)} \leq \frac{R}{C}-\left\|\underline Z \right\| _{\alpha,\alpha^{\prime},T} ,
			\end{equation*}
			we can proceed the above step with the same time $T^{\ast\ast}$ and then derive that $v_1=v_2$ on $[0,2T^{\ast\ast}\wedge T^{\ast}]$. By iterating, we finally conclude that $v_1=v_2$ on $[0, T^{\ast}]$.
		\end{proof}

		\subsection{A priori estimates for $m\geq1$}\label{sec3.2 A priori estimates}
		Recall the definition of $\Psi\left( v,\underline{Z}\right)$ in \eqref{Psi}. For simplicity, we set $\Psi\left( v,\underline{Z}\right)= -\nu v^{m+1} \overline{v}^{m}+\Psi^{\prime}\left( v,\underline{Z}\right)$, where $\Psi^{\prime}\left( v,\underline{Z}\right)$ is defined as
		\begin{equation}\label{Psi^prime}
			\Psi^{\prime}\left( v,\underline{Z}\right):=-\nu\sum_{0\leq i+j<2m+1}\binom{m+1}{i}\binom{m}{j}v^i\overline{v}^jZ^{:m+1-i,m-j:} +\tau(v+Z).
		\end{equation}
		Testing the equation \eqref{remainder} with $|v|^{2p-2}\overline{v}$, we get the following expression of $	\|v(t)\|_{L^{2 p}}^{2 p}$ (see \cite[Proposition 5.2]{Hoshino2018}, \cite[Proposition 3.3]{Matsuda2020} and \cite[Proposition 6.2]{Trenberth2019} for the similar argument).
		\begin{proposition}\label{L^p of v(t)}
			Let $T>0$ and $v\in C\left( \left( 0,T\right]; \mathcal{C}^{\beta}\right) $ be a mild solution to \eqref{remainder}. Then for all $0<s_0\leq t\leq T$ and every $p\geq1$, 
			\begin{align*}
				&\frac{1}{2p}\left( \|v(t)\|_{L^{2 p}}^{2 p}-\|v(s_0)\|_{L^{2 p}}^{2 p}\right) \label{integral expression for L^p of v(t)} \\= &\, -\mathrm{Re}\left[(\i+\mu) \int_{s_0}^t \left\langle \nabla\left\{ |v(s)|^{2 p-2} \overline{v(s)}\right\}, \nabla v(s) \right\rangle \mathrm{d} s\right]  -\int_{s_0}^t\|v(s)\|_{L^{2 p}}^{2 p}\mathrm{d} s \nonumber\\&\,- \mathrm{Re} \nu \int_{s_0}^t\|v(s)\|_{L^{2 p+2m}}^{2 p+2m}\mathrm{d} s + \int_{s_0}^t\left\langle |v(s)|^{2 p-2},\mathrm{Re}\left(\overline{v(s)}\Psi^{\prime}\left( v(s),\underline{Z}(s)\right)\right) \right\rangle \mathrm{d} s \nonumber.
			\end{align*}
			In particular, differentiating with respect to $t$, we get that for every $t\in (0,T]$,
			\begin{align*}
				\frac{1}{2p}\partial_t\|v(t)\|_{L^{2 p}}^{2 p}
				= & -\mathrm{Re}\left[(\i+\mu) \left\langle  \nabla\left\{ |v(t)|^{2 p-2} \overline{v(t)}\right\} , \nabla v(t)\right\rangle  \right]-\|v(t)\|_{L^{2 p}}^{2 p} 
				\\& - \mathrm{Re} \nu \|v(t)\|_{L^{2 p+2m}}^{2 p+2m}+ \left\langle |v(t)|^{2 p-2},\mathrm{Re}\left(\overline{v(t)}\Psi^{\prime}\left( v(t),\underline{Z}(t)\right)\right) \right\rangle \nonumber .
			\end{align*}
		\end{proposition}	
		
		\begin{remark}
			Theorem \ref{Local existence and uniqueness} implies that $v\in C\left( \left( 0,T\right]; \mathcal{C}^{\beta}\right) $ for some $\beta< 1$ since the initial data $u_0$ of \eqref{remainder} belongs to $ \mathcal{C}^{-\alpha_0}$ with $\alpha_0>0$. However, the first order spatial derivatives and some time regularity of $v$ are required to show Proposition \ref{L^p of v(t)}. One can prove that for fixed $t > 0$, $v(t)\in \mathcal{C}^{2-\epsilon} $ for any $\epsilon>0$ if $v(0)=0$ and $v$ is a H\"older continuous function from $(0,T ]$ to $L^{\infty}(\mathbb{T}^2)$ with some positive exponent (see \cite[Proposition 3.1]{Matsuda2020} for the proof in the case $m=1$).
		\end{remark}

		\begin{proposition}\label{expression of v_2p}
			Let $T>0$, $v\in C\left( \left( 0,T\right]; \mathcal{C}^{\beta}\right) $ be a mild solution to \eqref{remainder} and $1\leq p\leq 1+\mu(\mu+\sqrt{1+\mu^2})$. For every $0\leq \delta<1$ such that
			\begin{equation}\label{delta condition}
				\frac{p-1}{\mu\left(\mu+\sqrt{1+\mu^2}\right)} \leq 1-\delta,
			\end{equation}
			we have that for every $t\in (0,T]$,
			\begin{align*}
				& \frac{1}{2 p}\partial_t\|v(t)\|_{L^{2 p}}^{2 p}+\delta \mu \left\||\nabla v(t)|^2|v(t)|^{2 p-2}\right\|_{L^1}+ \mathrm{Re}\nu \|v(t)\|_{L^{2 p+2m}}^{2 p+2m}+\|v(t)\|_{L^{2 p}}^{2 p} \\
				\leq&\,\left\langle |v(t)|^{2 p-2},\mathrm{Re}\left(\overline{v(t)}\Psi^{\prime}\left( v(t),\underline{Z}(t)\right)\right) \right\rangle.
			\end{align*}
			In particular, integrating with respect to $t$, we have that for all $0<s_0\leq t\leq T$,
			\begin{align*}
				& \frac{1}{2 p}\left(\|v(t)\|_{L^{2 p}}^{2 p}-\left\|v(s_0)\right\|_{L^{2 p}}^{2 p}\right)+\delta \mu \int_{s_0}^t\left\||\nabla v(s)|^2|v(s)|^{2 p-2}\right\|_{L^1}\mathrm{d} s \\+&\, \mathrm{Re}\nu \int_{s_0}^t\|v(s)\|_{L^{2 p+2m}}^{2 p+2m}\mathrm{d} s+\int_{s_0}^t\|v(s)\|_{L^{2 p}}^{2 p}\mathrm{d} s 
				\leq\int_{s_0}^t\left\langle |v(s)|^{2 p-2},\mathrm{Re}\left(\overline{v(s)}\Psi^{\prime}\left( v(s),\underline{Z}(s)\right)\right) \right\rangle \mathrm{d} s.
			\end{align*}
		\end{proposition}
		\begin{proof}
			By Proposition \ref{L^p of v(t)}, it suffices to show that for every $0\leq \delta<1$ satisfying \eqref{delta condition},																																																																																																								
			\begin{equation}\label{greater}
				\delta \mu\left\||\nabla v(t)|^2|v(t)|^{2 p-2}\right\|_{L^1}\leq \mathrm{Re}\left[(\i+\mu)  \left\langle \nabla\left\{|v(t)|^{2 p-2} \overline{v(t)}\right\}, \nabla v(t) \right\rangle \right], \quad t\in(0,T] .
			\end{equation}
			Note that
			\begin{align*}
				&2\mathrm{Re}\left[(\i+\mu) \left\langle \nabla\left\{|v(t)|^{2 p-2} \overline{v(t)}\right\}, \nabla v(t) \right\rangle \right]\\=&\,(\i+\mu)\left[ p\left\langle\left|v(t) \right| ^{2p-2}, \left| \nabla v(t)\right| ^2\right\rangle +(p-1)\left\langle \left|v(t) \right| ^{2p-4}, \overline{v(t)}^2\left( \nabla v(t)\right) ^2 \right\rangle  \right] 
				\\&\,+(-\i+\mu)\left[ p\left\langle\left|v(t) \right| ^{2p-2}, \left| \nabla v(t)\right| ^2\right\rangle +(p-1)\left\langle \left|v(t) \right| ^{2p-4}, {v(t)}^2\left( \nabla \overline{v(t)}\right) ^2 \right\rangle  \right] 
				\\=&\, 2p\mu \left\langle \left|v(t) \right| ^{2p-4}, x_t^2+y_t^2 \right\rangle +2(p-1) \left\langle \left|v(t) \right| ^{2p-4},  \mu x_t^2- \mu y_t^2 -2x_ty_t \right\rangle
				\\=&\, 2\left\langle  |v(t)|^{2 p-4},   \mu(2 p-1)x_t^2+ \mu y_t^2 -2(p-1)x_ty_t\right\rangle ,
			\end{align*}
			where $x_t=\frac12\nabla( v(t)\overline{v(t)} )=\frac12\nabla\left|v(t) \right| ^2$ and $y_t=\frac{\i}{2}(v(t)\nabla\overline{v(t)}-\overline{v(t)}\nabla v(t))$ are the real and imaginary parts of $\overline{v(t)}\nabla v(t)$ respectively. Then 
			\begin{align*}
				&2\mathrm{Re}\left[(\i+\mu) \left\langle \nabla\left\{|v(t)|^{2 p-2} \overline{v(t)}\right\}, \nabla v(t) \right\rangle \right]- 2	\delta \mu \left\||\nabla v(t)|^2|v(t)|^{2 p-2}\right\|_{L^1}
				\\=&\,2\mathrm{Re}\left[(\i+\mu) \left\langle \nabla\left\{|v(t)|^{2 p-2} \overline{v(t)}\right\}, \nabla v(t) \right\rangle \right]- 2	\delta \mu \left\langle |v(t)|^{2 p-4}, x_t^2+y_t^2 \right\rangle 
				\\=&\, 2\left\langle  |v(t)|^{2 p-4},   \mu(2 p-1-\delta)x_t^2+ \mu(1-\delta)y_t^2 -2(p-1)x_ty_t\right\rangle.
			\end{align*}
			To prove \eqref{greater}, it suffices to prove that the quadratic form 
			$$\mu(2 p-1-\delta)x_t^2+ \mu(1-\delta)y_t^2 -2(p-1)x_ty_t$$
			is non-negative. The matrix corresponding to this quadratic form is
			\begin{equation*}
				\begin{pmatrix}
					\mu(2 p-1-\delta) \quad & -(p-1)\\
					-(p-1)\quad &\mu(1-\delta)
				\end{pmatrix}.
			\end{equation*}
			Since $\mu(2 p-1-\delta)>0$ and when $\delta$ satisfies \eqref{delta condition}, the determinant of this matrix is non-negative, we complete the proof.
		\end{proof}
		Combining the dissipative effects of the Laplacian operator and the nonlinear term with Proposition \ref{expression of v_2p}, we derive a prioir estimates for the remainder $v$.
		\begin{proposition}\label{priori estimate}
			Let $T>0$, $v\in C\left( \left( 0,T\right]; \mathcal{C}^{\beta}\right) $ be a mild solution to \eqref{remainder}, $1\leq p\leq 1+\mu(\mu+\sqrt{1+\mu^2})$ and $0\leq \delta<1$ satisfy \eqref{delta condition}. Then there exist constant $C>0$ and $p_0>0$, such that for every $0<t_0<t\leq T$,  
			\begin{align}
				&\qquad\qquad\qquad\qquad\qquad\left\| v(t)\right\| _{L^{2p}}^{2p}  \leq C \left\| \underline{Z}\right\|_{\alpha,\alpha^{\prime},T}^{p_0}  t^{-\frac{p}{m}} ,\label{eq: Lp est}\\
				&\frac{\mathrm{Re}\nu}{2} \int_{t_0}^t\|v(s)\|_{L^{2 p+2m}}^{2 p+2m}\mathrm{d} s +\frac{\delta \mu}{2} \int_{t_0}^t\left\||\nabla v(s)|^2|v(s)|^{2 p-2}\right\|_{L^1}\mathrm{d} s + \int_{t_0}^t\|v(s)\|_{L^{2 p}}^{2 p}\mathrm{d} s\\
				\leq &\, C\left\| \underline{Z}\right\|_{\alpha,\alpha^{\prime},t_0}^{p_0}  t_0^{-\frac{p}{m}} + C \sum_{i=0}^{m+1}\sum_{j=0}^{m} \int_{t_0}^t \left\| Z^{:m+1-i,m-j:}(s) \right\| _{\mathcal{C}^{-\alpha}}^{p_0}   \mathrm{d} s \label{eq: Lp+2m est 2}
					\\
					\leq &\, C \left\| \underline{Z}\right\|_{\alpha,\alpha^{\prime},T}^{p_0}  t_0^{-\frac{p}{m}}\label{eq: Lp+2m est}.
				\end{align}
		\end{proposition}

		\begin{proof}
			By Proposition \ref{expression of v_2p}, to estimate $\left\| v(s)\right\| _{L^{2p}}^{2p}$ for $0<s\leq t\leq T$, we need to estimate $\left\langle |v(s)|^{2 p-2},\mathrm{Re}\left(\overline{v(s)}\Psi^{\prime}\left( v(s),\underline{Z}(s)\right)\right) \right\rangle$.	From the definition of $\Psi^{\prime}\left( v(s),\underline{Z}(s)\right)$ (see \eqref{Psi^prime}), we know that
			\begin{align}
				& \left\langle |v(s)|^{2 p-2},\mathrm{Re}\left(\overline{v(s)}\Psi^{\prime}\left( v(s),\underline{Z}(s)\right)\right) \right\rangle \label{total}\\=&\, \left\langle |v(s)|^{2 p-2},\mathrm{Re}\Bigg( -\nu\sum_{0\leq i+j<2m+1}\binom{m+1}{i}\binom{m}{j}v(s)^i\overline{v(s)}^{j+1}Z^{:m+1-i,m-j:} (s)  \Bigg) \right\rangle \nonumber
				\\&\,+ \left\langle |v(s)|^{2 p-2},\mathrm{Re}\left( \tau(|v(s)|^2+\overline{v(s)}Z(s))  \right) \right\rangle.\nonumber
			\end{align}
			Let
			\begin{equation*}
				A_s=\left\| v(s)\right\|_{L^{2 p+2m}}^{2 p+2m},\quad B_s=\left\||\nabla v(s)|^2|v(s)|^{2 p-2}\right\|_{L^1}.
			\end{equation*} 
			We show that for $0<s\leq t$, each term of $\left\langle |v(s)|^{2 p-2},\mathrm{Re}\left(\overline{v(s)}\Psi^{\prime}\left( v(s),\underline{Z}(s)\right)\right) \right\rangle$, namely 
			\begin{align*}
				&\left\| v(s)\right\|_{L^{2 p}}^{2 p}, \quad \left\langle |v(s)|^{2 p-2},\mathrm{Re}\left(\tau\overline{v(s)}Z(s)\right) \right\rangle,\\
				& 	\left\langle |v(s)|^{2 p-2},\mathrm{Re}\left(\nu v(s)^i\overline{v(s)}^{j+1}Z^{:m+1-i,m-j:}(s)\right) \right\rangle
			\end{align*} with $0\leq i\leq m+1, 0\leq j \leq m, (i,j)\neq(m+1,m)$, can be controlled by $A_s$ and $B_s$. 
			
			According to Proposition \ref{duality property}, for $\alpha\in(0,1)$,
			\begin{align}
				&\left|\left\langle |v(s)|^{2 p-2},\mathrm{Re}\left(\nu v(s)^i\overline{v(s)}^{j+1}Z^{:m+1-i,m-j:}(s)\right) \right\rangle \right| \label{ij estimate}\\ \lesssim &\,\left\langle |v(s)|^{2 p+i+j-1},\left| Z^{:m+1-i,m-j:}(s)\right|  \right\rangle\lesssim \left\| |v(s)|^{2 p+i+j-1} \right\| _{\mathcal{B}_{1,1}^{\alpha}}\left\| Z^{:m+1-i,m-j:}(s) \right\| _{\mathcal{C}^{-\alpha}}.\nonumber
			\end{align} 
			Using Proposition \ref{gradient estimate},
			\begin{equation}\label{b11}
				\left\| |v(s)|^{2 p+i+j-1} \right\| _{\mathcal{B}_{1,1}^{\alpha}}\lesssim \left\| |v(s)|^{2 p+i+j-1}\right\| _{L^1}^{1-\alpha}  \left\| \nabla |v(s)|^{2 p+i+j-1}  \right\| _{L^1}^{\alpha}  + \left\| |v(s)|^{2 p+i+j-1} \right\| _{L^1}.
			\end{equation}
			By Cauchy-Schwarz inequality,
			\begin{equation}\label{nabla}
				\left\| \nabla |v(s)|^{2 p+i+j-1}  \right\| _{L^1} \lesssim B_s^{\frac12} \left\| |v(s)|^{2 (p+i+j-1)}\right\| _{L^1}^{\frac12} . 
			\end{equation}
			Applying Sobolev inequality $\left\| f\right\| _{L^q}\lesssim \left(\left\| f\right\| _{L^2} ^2+ \left\| \nabla f\right\| _{L^2} ^2\right) ^{\frac12}$ for every $q<\infty$ (see \cite[Theorem 6.5.1]{Bergh1976} or \cite[Theorem 2.13]{Sawano2018}), taking $q=\frac{2(p+i+j-1)}{p}$ specifically and combining with Jensen's inequality, we have 
			\begin{equation}\label{2(p+i+j-1)}
				\left\| |v(s)|^{2 (p+i+j-1)}\right\| _{L^1}^{\frac12} \lesssim \left\|  |v(s)|^{2 p}  \right\| _{L^1}^{\frac{p+i+j-1}{2p}}+B_s^{\frac{p+i+j-1}{2p}}\lesssim A_s^{\frac{p+i+j-1}{2p+2m}}+B_s^{\frac{p+i+j-1}{2p}}.
			\end{equation}
			Moreover, by Jensen's inequality again, 
			\begin{equation}\label{L1}
				\left\| |v(s)|^{2 p+i+j-1} \right\| _{L^1}\lesssim A_s^{\frac{2p+i+j-1}{2p+2m}}.
			\end{equation}
			Therefore, combining with \eqref{b11}, \eqref{nabla}, \eqref{2(p+i+j-1)} and \eqref{L1}, we know that
			\begin{equation}\label{B11}
				\left\| |v(s)|^{2 p+i+j-1} \right\| _{\mathcal{B}_{1,1}^{\alpha}}\lesssim A_s^{\frac{2p+i+j-1-\alpha p}{2p+2m}} B_s^{\frac{\alpha}{2}}+  A_s^{\frac{2p+i+j-1}{2p+2m}(1-\alpha)} B_s^{\frac{2p+i+j-1}{2p}\alpha} +A_s^{\frac{2p+i+j-1}{2p+2m}}.
			\end{equation}
			Here we take  $0<\alpha<\frac{1}{m}\wedge\frac{p}{m(2m-1+2p)}= \frac{p}{m(2m-1+2p)}$ such that 
			\begin{equation*}
				\frac{2p+2m-1-\alpha p}{2p+2m} +\frac{\alpha}{2}<1,\quad \frac{2p+2m-1}{2p+2m}(1-\alpha) +\frac{2p+2m-1}{2p}\alpha<1.
			\end{equation*}
			Then we can find exponents $0<\gamma_{i,j}^k<1$, $1\leq k\leq 4$, such that 
			\begin{equation*}
				\frac{1}{\gamma_{i,j}^1}\frac{2p+2m-1-\alpha p}{2p+2m} +\frac{1}{\gamma_{i,j}^2}\frac{\alpha}{2}=1,
				\quad 
				\frac{1}{\gamma_{i,j}^3} \frac{2p+2m-1}{2p+2m}(1-\alpha) +\frac{1}{\gamma_{i,j}^4}\frac{2p+2m-1}{2p}\alpha=1.
			\end{equation*}
			According to \eqref{ij estimate}, \eqref{B11} and Young's inequality in the form $ab\leq \kappa a^{q_1}+C(\kappa)b^{q_2}$, where $\kappa>0$ is arbitrary, $a,b\geq0$, $q_1, q_2>1$ satisfy $\frac{1}{q_1}+\frac{1}{q_2}=1$ and $C(\kappa)=\frac{1}{q_2}(q_1\kappa)^{-\frac{q_2}{q_1}}$, we obtain that 
			\begin{align*}
				&\left|\left\langle |v(s)|^{2 p-2},\mathrm{Re}\left(\nu v(s)^i\overline{v(s)}^{j+1}Z^{:m+1-i,m-j:}(s)\right) \right\rangle \right|\\ \lesssim &\,\left( A_s^{\gamma_{i,j}^1}+B_s^{\gamma_{i,j}^2}+A_s^{\gamma_{i,j}^3}+B_s^{\gamma_{i,j}^4}+A_s^{\gamma_{i,j}^5}\right) \left\| Z^{:m+1-i,m-j:} (s) \right\| _{\mathcal{C}^{-\alpha}},
			\end{align*}
			where $0<\gamma_{i,j}^5=\frac{2p+i+j-1}{2p+2m}<1$.
			Using Young's inequality again, we get the bound 
			\begin{align}
				&\binom{m+1}{i}\binom{m}{j}\left|\left\langle |v(s)|^{2 p-2},\mathrm{Re}\left(\nu v(s)^i\overline{v(s)}^{j+1}Z^{:m+1-i,m-j:}(s)\right) \right\rangle \right| \label{estimate1}\\\leq&\,  \frac{1}{ (m+1)(m+2)} \left( \frac{\mathrm{Re}\nu}{4} A_s+\frac{\delta\mu}{2} B_s\right) +C\sum_{k=1}^{5}\left\| Z^{:m+1-i,m-j:}(s) \right\| _{\mathcal{C}^{-\alpha}}^{\frac{1}{1-\gamma_{i,j}^k}} \nonumber ,
			\end{align}
			where $0\leq \delta<1$ satisfy \eqref{delta condition} and $C$ is a constant depending on $m$, $\nu$, $\delta$, $\mu$, $\alpha$, $p$ and $\gamma_{i,j}^k$, $1\leq k\leq 5$.
			
			Now we estimate $\left\langle |v(s)|^{2 p-2},\mathrm{Re}\left(\tau\overline{v(s)}Z(s)\right) \right\rangle$. By the similar argument as above, we can find $0<\zeta_{m,m}^k<1$, $1\leq k\leq 5$, such that   
			\begin{align}
				&\left| \left\langle |v(s)|^{2 p-2},\mathrm{Re}\left(\tau\overline{v(s)}Z(s)\right) \right\rangle  \right| \label{estimate2}\\ \leq &\,\frac{1}{ (m+1)(m+2)} \left( \frac{\mathrm{Re}\nu}{4} A_s+\frac{\delta\mu}{2} B_s\right) +C\sum_{k=1}^{5}\left\| Z(s) \right\| _{\mathcal{C}^{-\alpha}}^{\frac{1}{1-\zeta_{m,m}^k}},\nonumber
			\end{align}
			where $C$ is a constant depending on $m$, $\nu$, $\delta$, $\mu$, $\alpha$, $p$, $\tau$ and $\zeta_{m,m}^k$, $1\leq k\leq 5$.

			Next we estimate $\left\| v(s)\right\|_{L^{2 p}}^{2 p}$. Applying Jensen's inequality and Young's inequality, we have that 
			\begin{equation}\label{estimate3}
				|\tau|\left\| v(s)\right\|_{L^{2 p}}^{2 p}\leq C	|\tau| A_s^{\frac{p}{p+m}}\leq \frac{\mathrm{Re}\nu}{4}A_s +C,
			\end{equation}
			where $C>0$ is a constant depending on $\nu,p,m$ and $\tau$.
			
			For $1\leq k\leq 5$, let $p_{m+1,m}^k=1$, $p_{m,m}^k=\frac{1}{1-\zeta_{m,m}^k}\vee \frac{1}{1-\gamma_{m,m}^k} $ if $\left\| Z(s) \right\| _{\mathcal{C}^{-\alpha}}\geq 1$ and $p_{m,m}^k=\frac{1}{1-\zeta_{m,m}^k}\wedge \frac{1}{1-\gamma_{m,m}^k} $ if $\left\| Z(s) \right\| _{\mathcal{C}^{-\alpha}}< 1$, and  $p_{i,j}^k=\frac{1}{1-\gamma_{i,j}^k}$, where $0\leq i\leq m+1$, $0\leq j\leq m$ and $(i,j)\neq (m+1,m), (m,m)$. Then combining with \eqref{total}, \eqref{estimate1}, \eqref{estimate2} and \eqref{estimate3}, we get that 
			\begin{align}
				&\left| 	\left\langle |v(s)|^{2 p-2},\mathrm{Re}\left(\overline{v(s)}\Psi^{\prime}\left( v(s),\underline{Z}(s)\right)\right) \right\rangle \right| \label{vRev}\\ \leq&\,   \frac{\mathrm{Re}\nu}{2} A_s+\frac{\delta\mu}{2} B_s+C\sum_{i=0}^{m+1}\sum_{j=0}^{m} \sum_{k=1}^{5}  \left\| Z^{:m+1-i,m-j:}(s) \right\| _{\mathcal{C}^{-\alpha}}^{p_{i,j}^k} \nonumber \\\leq&\, \frac{\mathrm{Re}\nu}{2} A_s+ \frac{\delta\mu}{2} B_s+ C \sum_{i=0}^{m+1}\sum_{j=0}^{m} \left\| Z^{:m+1-i,m-j:}(s) \right\| _{\mathcal{C}^{-\alpha}}^{\tilde{p}},\nonumber
			\end{align}
			where $C>0$ is a constant depending on $m$, $\nu$, $\delta$, $\mu$, $\alpha$, $p$, $\tau$ and all $p_{i,j}^k$, and
			\begin{equation}\label{tildep}
				\tilde{p}=\max\left\lbrace p_{i,j}^k: 1\leq k\leq 5, 0\leq i\leq m+1,0\leq j\leq m \right\rbrace .
			\end{equation}
			By Proposition \ref{expression of v_2p}, we have that
			\begin{align}
				&\frac{1}{2 p}\partial_s\|v(s)\|_{L^{2 p}}^{2 p}+\frac{\delta\mu}{2} \left\||\nabla v(s)|^2|v(s)|^{2 p-2}\right\|_{L^1} +\frac{\mathrm{Re}\nu}{2} \|v(s)\|_{L^{2 p+2m}}^{2 p+2m} +\|v(s)\|_{L^{2 p}}^{2 p}  \nonumber \\
				\leq &\, C \sum_{i=0}^{m+1}\sum_{j=0}^{m} \left\| Z^{:m+1-i,m-j:}(s) \right\| _{\mathcal{C}^{-\alpha}}^{\tilde{p}}  \label{total estimate}.
			\end{align}
			Let $t>s$. For $r\in[s,t]$,
			\begin{align*}
				&\sum_{i=0}^{m+1}\sum_{j=0}^{m} \left\| Z^{:m+1-i,m-j:}(r) \right\| _{\mathcal{C}^{-\alpha}}^{\tilde{p}}\\\leq&\, 1+ \sum_{0<i+j\leq 2m+1}  s^{-\alpha^{\prime}(i+j-1)\tilde{p}} \sup_{s\leq r\leq t }r^{\alpha^{\prime}(i+j-1)\tilde{p}}\left\| Z^{:i,j:}(r) \right\| _{\mathcal{C}^{-\alpha}}^{\tilde{p}} \\
				\lesssim &\,	\left\| \underline{Z}\right\|^{\tilde{p}}_{\alpha,\alpha^{\prime},t}\sum_{j=0}^{2m}  s^{-\alpha^{\prime}j\tilde{p}}.
			\end{align*}
			Using Jensen's inequality, we get that for $0<s\leq r\leq t\leq T$,
			\begin{align*}
				&\partial_r\|v(r)\|_{L^{2 p}}^{2 p}+  \left( \left\| v(r)\right\| _{L^{2p}}^{2p}\right)^{\frac{p+m}{p}} \\
				\lesssim&\,  \frac{1}{2p}	\partial_r\|v(r)\|_{L^{2 p}}^{2 p}+  \frac{\mathrm{Re}\nu}{2} \|v(r)\|_{L^{2 p+2m}}^{2 p+2m}+\frac{\delta\mu}{2} \left\||\nabla v(r)|^2|v(r)|^{2 p-2}\right\|_{L^1}\\
				\leq&\, C 	\left\| \underline{Z}\right\|^{\tilde{p}}_{\alpha,\alpha^{\prime},t}\sum_{j=0}^{2m}  s^{-\alpha^{\prime}j\tilde{p}} .
			\end{align*}
 			Let $f(r)=\left\| v(r)\right\| _{L^{2p}}^{2p}$ for $r\in[s,t]$. Then according to \cite[Lemma 3.8]{TW2018}, for $s<r\leq t$,
			\begin{align*}
				f(r)&\leq \frac{f(s)}{\left( 1+\frac{m}{2p}(r-s)f(s)^{\frac{m}{p}}\right) ^{\frac{p}{m}}}\vee \Bigg[ 2C  	\left\| \underline{Z}\right\|^{\tilde{p}}_{\alpha,\alpha^{\prime},t}\sum_{j=0}^{2m}  s^{-\alpha^{\prime}j\tilde{p}}\Bigg] ^{\frac{p}{p+m}}\\
				&\leq C \Bigg((r-s)^{-\frac{p}{m}}\vee \Bigg[ 	\left\| \underline{Z}\right\|^{\tilde{p}}_{\alpha,\alpha^{\prime},t}\sum_{j=0}^{2m}  s^{-\alpha^{\prime}j\tilde{p}}\Bigg] ^{\frac{p}{p+m}}\Bigg).
			\end{align*}
			Taking $s=t/2$, $r=t$ and $\alpha^{\prime}>0$ small enough such that
			\begin{align}\label{condition for alpha 2}
				\frac{2m\alpha^{\prime}\tilde{p}p}{p+m}\leq \frac{p}{m},\quad 2m\alpha^{\prime}\tilde{p}<1,
			\end{align}
			then
			\begin{align}
				\left\| v(t)\right\| _{L^{2p}}^{2p}  &\leq C \left\| \underline{Z}\right\|_{\alpha,\alpha^{\prime},T}^{\frac{p\tilde{p}}{p+m}} \left( t^{-\frac{p}{m}}\vee \Bigg(\sum_{j=0}^{2m}  t^{-\alpha^{\prime}j\tilde{p}}\Bigg)^{\frac{p}{p+m}}\right) \nonumber\\&\leq C \left\| \underline{Z}\right\|_{\alpha,\alpha^{\prime},T}^{\frac{p\tilde{p}}{p+m}}  t^{-\frac{p}{m}},\quad  0< t\leq T. \label{eq: Lp est2}
			\end{align}
			Integrating both sides of \eqref{total estimate} from $t_0$ to $t$ and applying \eqref{eq: Lp est2}, we get that 
			\begin{align*}
				&\frac{\mathrm{Re}\nu}{2} \int_{t_0}^t\|v(s)\|_{L^{2 p+2m}}^{2 p+2m}\mathrm{d} s +\frac{\delta \mu}{2} \int_{t_0}^t\left\||\nabla v(s)|^2|v(s)|^{2 p-2}\right\|_{L^1}\mathrm{d} s +\int_{t_0}^t\|v(s)\|_{L^{2 p}}^{2 p}\mathrm{d} s   \\ 
				\lesssim &\, \left\| \underline{Z}\right\|_{\alpha,\alpha^{\prime},t_0}^{\frac{p\tilde{p}}{p+m}} t_0^{-\frac{p}{m}} +\sum_{i=0}^{m+1}\sum_{j=0}^{m} \int_{t_0}^t \left\| Z^{:m+1-i,m-j:}(s) \right\| _{\mathcal{C}^{-\alpha}}^{\tilde{p}}  \mathrm{d} s\\
				\lesssim  &\, \left\| \underline{Z}\right\|_{\alpha,\alpha^{\prime},T}^{\tilde{p}} t_0^{-\frac{p}{m}} +	\left\| \underline{Z}\right\|^{\tilde{p}}_{\alpha,\alpha^{\prime},T} \sum_{j=0}^{2m} \int_{t_0}^t s^{-\alpha^{\prime}j\tilde{p}} \mathrm{d} s
				\\ 
				\lesssim &\, \left\| \underline{Z}\right\|_{\alpha,\alpha^{\prime},T}^{\tilde{p}}  (t_0^{-\frac{p}{m}} \vee 1 ) 	\lesssim \left\| \underline{Z}\right\|_{\alpha,\alpha^{\prime},T}^{\frac{p\tilde{p}}{p+m}}  t_0^{-\frac{p}{m}}.
			\end{align*}
			Then we finish the proof. 
		\end{proof}
		
	\subsection{Improved a priori estimates for $m=1$}\label{sec3.3 Better A priori estimates}
	
	We establish an estimate for $\|v\|_{L^{2p}}$ with some $p>1$ in Proposition \ref{priori estimate}. Then by Proposition \ref{L^p embedding}, we actually obtain an estimate for $\|v\|_{\mathcal{B}_{2p,\infty}^0}$. We next show that this estimate can be improved when $m=1$.

		\begin{proposition}\label{prop: B-epsilon-6p est}
		Let $T>0$ and $v\in C\left( \left( 0,T\right]; \mathcal{C}^{\beta}\right) $ be a mild solution to \eqref{remainder} with $m=1$. Then there exist $\epsilon_0>0$, $p>1$, $p_0>0$ and $\kappa>0$, such that for any $0<\epsilon<\epsilon_0$ and any $0<t_0\leq T$,   
		\begin{equation*}
			\int_{t_0}^{T}\left\| v(t)\right\|^3 _{\mathcal{B}_{6p,\infty}^{\epsilon}} \mathrm{d}t \lesssim \left\| \underline{Z}\right\|_{\alpha,\alpha^{\prime},T}^{p_0}  t_0^{-\kappa}.
		\end{equation*}
	\end{proposition}
	\begin{proof}
		A solution $v\in C\left( \left( 0,T\right]; \mathcal{C}^{\beta}\right)$ satisfies
		\begin{align*}
			v(t)= P_{t-\frac{t_0}{2}}v\Big(\frac{t_0}{2}\Big) +\int_{\frac{t_0}{2}}^{t} P_{t-s} \Psi\left(v(s),\underline{Z}(s)\right)\mathrm{d}s, \quad t_0\leq t\leq T.
		\end{align*}
		By Proposition \ref{priori estimate}, we first take $1<p<6/5$ (but close to 1 since $\mu>0$) such that \eqref{eq: Lp est} and \eqref{eq: Lp+2m est} hold. Then we let $0<\epsilon_0=\frac{(p-1)(8p-1)}{3p(p+1)}<  \frac{1}{3p}$. For any $0<\epsilon<\epsilon_0$, we estimate the $L^{3}([t_0,T]; \mathcal{B}_{6p,\infty}^{\epsilon})$-norm of each term on the right hand side.
		
		For the initial data part, by Proposition \ref{regularity embedding}, \ref{heat kernel smoothing} and \ref{L^p embedding}, and estimate \eqref{eq: Lp est} with $m=1$, there exists $p_0>0$ such that 
		\begin{align*}
			\left\| P_{t-\frac{t_0}{2}}v\Big(\frac{t_0}{2}\Big)\right\| _{\mathcal{B}_{6p,\infty}^{\epsilon}}
			&\lesssim 	\left\| P_{t-\frac{t_0}{2}}v\Big(\frac{t_0}{2}\Big)\right\| _{\mathcal{B}_{2p,\infty}^{\epsilon+\frac{2}{3p}}} \lesssim \Big(t-\frac{t_0}{2}\Big)^{-\frac{\epsilon}{2}- \frac{1}{3p}}\left\| v\Big(\frac{t_0}{2}\Big)\right\| _{\mathcal{B}_{2p,\infty}^{0}}\\
			&\lesssim \Big(t-\frac{t_0}{2}\Big)^{-\frac{\epsilon}{2}- \frac{1}{3p}} \left\| v\Big(\frac{t_0}{2}\Big)\right\| _{L^{2p}} \lesssim \left\| \underline{Z}\right\|_{\alpha,\alpha^{\prime},T}^{p_0} t_0^{-\frac{1}{2}} \Big(t-\frac{t_0}{2}\Big)^{-\frac{\epsilon}{2}- \frac{1}{3p}}.
		\end{align*}
		Note that the exponent $p_0>0$ will vary from line to line in this proof. Then
		\begin{align*}
			\int_{t_0}^{T} 	 \left\| P_{t-\frac{t_0}{2}}v\Big(\frac{t_0}{2}\Big)\right\|^3 _{\mathcal{B}_{6p,\infty}^{\epsilon}} \mathrm{d}t
			& \lesssim \left\| \underline{Z}\right\|_{\alpha,\alpha^{\prime},T}^{p_0} t_0^{-\frac{3}{2}-\frac{3\epsilon}{2}- \frac{1}{p}}.
		\end{align*}
		
		Recall that $\Psi(v,\underline{Z})=-\nu v^{2} \overline{v}-\nu\sum_{0\leq i+j<3}\binom{2}{i}\binom{1}{j}v^i\overline{v}^jZ^{:2-i,1-j:} +\tau(v+Z)$ in \eqref{Psi}. Let $p_1=\frac{2 m+2p}{2m+1}=\frac{2+2p}{3}$. Using Proposition \ref{regularity embedding}, \ref{heat kernel smoothing} and \ref{L^p embedding}, 
		\begin{align*}
			\left\| P_{t-s} v(s)^{2}\overline{v(s)}\right\| _{\mathcal{B}_{6p,\infty}^{\epsilon}}
			&\lesssim \left\| P_{t-s} v(s)^{2}\overline{v(s)}\right\|_{\mathcal{B}_{p_1,\infty}^{\epsilon+2(\frac{1}{p_1}-\frac{1}{6p})}}\lesssim (t-s)^{-(\frac{\epsilon}{2}+\frac{1}{p_1} -\frac{1}{6p} )} \left\|  v(s)^{2}\overline{v(s)}\right\| _{\mathcal{B}_{p_1,\infty}^{0}}\\
			& \lesssim (t-s)^{-(\frac{\epsilon}{2}+\frac{1}{p_1} -\frac{1}{6p} )}  \left\|  {v(s)}^{3}\right\| _{L^{p_1}} 
			=(t-s)^{-(\frac{\epsilon}{2}+\frac{1}{p_1} -\frac{1}{6p} )}  \left\|  {v(s)}^{2+2p}\right\| _{L^{1}}^{\frac{1}{p_1}}.
		\end{align*}
		Take $q= \frac{6+6p}{8p-1}$ such that $\frac{1}{p_1}+ \frac{1}{q}=1+\frac13$. Applying Young's convolution inequality and \eqref{eq: Lp+2m est} with $m=1$, 
		\begin{align}
			\int_{t_0}^{T} \left(\int_{\frac{t_0}{2}}^{t} 	\left\| P_{t-s} v(s)^{2}\overline{v(s)}\right\| _{\mathcal{B}_{6p,\infty}^{\epsilon}} \mathrm{d} s \right)^3 \mathrm{d} t \nonumber
			&\lesssim \int_{t_0}^{T} \left( \int_{\frac{t_0}{2}}^{t} (t-s)^{-(\frac{\epsilon}{2}+\frac{1}{p_1} -\frac{1}{6p} )}  \left\|  {v(s)}^{2+2p}\right\| _{L^{1}}^{\frac{1}{p_1}} \mathrm{d} s\right)^3 \mathrm{d} t \nonumber\\
			&\lesssim \left( \int_{0}^{T} s^{-q(\frac{\epsilon}{2}+\frac{1}{p_1} -\frac{1}{6p} )} \mathrm{d} s\right) ^{\frac{3}{q}} \left( \int_{\frac{t_0}{2}}^{T}   \left\|  {v(s)}^{2+2p}\right\| _{L^{1}} \mathrm{d} s\right) ^{\frac{3}{p_1}} \nonumber\\
			&\lesssim  \left\| \underline{Z}\right\|_{\alpha,\alpha^{\prime},T}^{p_0} t_0 ^{-\frac{3p}{ p_1}} \label{eq: reason for m=1},
		\end{align}
		where we use the fact that $0< q(\frac{\epsilon}{2}+\frac{1}{p_1} -\frac{1}{6p} )<1 $ since 
		\begin{align*}
			\frac{\epsilon}{2}< \frac{\epsilon_0}{2}  = \frac{(p-1)(8p-1)}{6p(p+1)}  = \frac{1}{q} -\frac{1}{p_1} +\frac{1}{6p}.
		\end{align*} 
		Using the similar argument and \eqref{eq: Lp est} with $m=1$, we know that 
		\begin{align*}
			\int_{t_0}^{T} \left( \int_{\frac{t_0}{2}}^{t} 	\left\| P_{t-s} v(s)\right\| _{\mathcal{B}_{6p,\infty}^{\epsilon}} \mathrm{d} s \right)^3 \mathrm{d} t 
			&\lesssim  \left(\int_{0}^{T} s^{-\frac{\epsilon}{2}-\frac{1}{3p}} \mathrm{d} s\right)^3  \int_{\frac{t_0}{2}}^{T}  \left\|  {v(t)}\right\| _{L^{2p}}^3 \mathrm{d} t  \\
			&\lesssim   \left\| \underline{Z}\right\|^{p_0}_{\alpha,\alpha^{\prime},T} t_0^{-\frac{3}{2}},
		\end{align*}
		where we use the fact that $\frac{\epsilon}{2}<\epsilon_0<\frac{1}{3p} < 1-\frac{1}{3p}$. By \eqref{norm of Z} and taking $\alpha<1/3$ such that $\alpha+\epsilon<2$, we have
		\begin{align*}
			\int_{t_0}^{T} \left( \int_{\frac{t_0}{2}}^{t} \left\| P_{t-s} Z(s)\right\| _{\mathcal{B}_{6p,\infty}^{\epsilon}} \mathrm{d} s \right)^3 \mathrm{d} t  
			&\lesssim 	\int_{t_0}^{T} \left( \int_{\frac{t_0}{2}}^{t} (t-s)^{-\frac{\alpha+\epsilon}{2}}\left\| Z(s)\right\| _{\mathcal{C}^{-\alpha}} \mathrm{d} s \right)^3 \mathrm{d} t \\
			&\lesssim  \left(\int_{0}^{T} s^{-\frac{\alpha+\epsilon}{2}} \mathrm{d} s\right)^3 \left\| \underline{Z}\right\|_{\alpha,\alpha^{\prime},T}^{p_0} \lesssim  \left\| \underline{Z}\right\|_{\alpha,\alpha^{\prime},T}^{p_0}  .
		\end{align*}
		For $i+j=0$, from Proposition \ref{embedding} and \ref{heat kernel smoothing}, Young's convolution inequality and \eqref{norm of Z},
		\begin{align*}
			\int_{t_0}^{T} \left(  \int_{\frac{t_0}{2}}^{t}	\left\| P_{t-s} Z^{:2,1:}(s) \right\| _{\mathcal{B}_{6p,\infty}^{\epsilon}} \mathrm{d} s \right)^3 \mathrm{d} t  
			&\lesssim \left\| \underline{Z}\right\|_{\alpha,\alpha^{\prime},T}^{p_0} \int_{t_0}^{T} \left( \int_{\frac{t_0}{2}}^{t} (t-s)^{-\frac{\alpha+\epsilon}{2}} s^{-2\alpha^{'}} \mathrm{d} s \right)^3 \mathrm{d} t    \\
			&\lesssim \left\| \underline{Z}\right\|_{\alpha,\alpha^{\prime},T}^{p_0}  \left(\int_{0}^{T} s^{-\frac{\alpha+\epsilon}{2}} \mathrm{d} s\right)^3 \int_{\frac{t_0}{2}}^{T} s^{-6\alpha^{'}} \mathrm{d} s \\
			&\lesssim \left\| \underline{Z}\right\|_{\alpha,\alpha^{\prime},T}^{p_0} {t_0}^{-6 \alpha^{'}} .
		\end{align*}
		For $i+j=1$, according to Proposition \ref{embedding}, \ref{regularity embedding}, \ref{heat kernel smoothing}, \ref{multiplicative structure}, the fact that $\alpha<1/3$ and \eqref{norm of Z},
		\begin{align*}
			\left\| P_{t-s} v(s)^i\overline{v(s)}^jZ^{:2-i,1-j:}(s) \right\| _{\mathcal{B}_{6p,\infty}^{\epsilon}} 
			&\lesssim 	\left\| P_{t-s} v(s)^i\overline{v(s)}^jZ^{:2-i,1-j:}(s) \right\| _{\mathcal{B}_{1,1}^{\epsilon+2(1-\frac{1}{6p})}} \\
			&\lesssim (t-s)^{-(\frac{\alpha+\epsilon}{2}+1-\frac{1}{6p})}\left\| v(s)^i\overline{v(s)}^jZ^{:2-i,1-j:}(s) \right\| _{\mathcal{B}_{1,1}^{-\alpha}}
			\\
			&\lesssim (t-s)^{-(\frac{\alpha+\epsilon}{2}+1-\frac{1}{6p})}\left\| v(s)^i\overline{v(s)}^j \right\| _{\mathcal{B}_{1,1}^{1/3}} \left\| Z^{:2-i,1-j:}(s) \right\| _{\mathcal{C}^{-\alpha}}\\
			&\lesssim (t-s)^{-(\frac{\alpha+\epsilon}{2}+1-\frac{1}{6p})} s^{-\alpha^{\prime}} \left\| v(s) \right\| _{\mathcal{B}_{1,1}^{1/3}} \left\|\underline{Z}\right\|_{\alpha,\alpha^{\prime},T}^{p_0}.
		\end{align*} 
		By Proposition \ref{gradient estimate}, H\"older inequality and \eqref{eq: Lp est},
		\begin{align*}
			\left\| v(s) \right\| _{\mathcal{B}_{1,1}^{1/3}}
			&\lesssim \left\| v(s)\right\| _{L^1}^{2/3}  \left\| \nabla v(s)  \right\| _{L^1}^{1/3}  + \left\| v(s) \right\| _{L^1}\\
			&\lesssim \left\| v(s)\right\| _{L^{2}}^{2/3}  \left\| \nabla v(s)  \right\| _{L^2}^{1/3}  + \left\| v(s) \right\| _{L^{2}}\\
			&\lesssim s^{-\frac{1}{3}} \left\| \nabla v(s)  \right\| _{L^2}^{1/3} + s^{-\frac{1}{2}}.
		\end{align*}
		Then applying Young's convolution inequality and \eqref{eq: Lp+2m est},
		\begin{align*}
			&	\int_{t_0}^{T} \left( \int_{\frac{t_0}{2}}^{t} 	\left\| P_{t-s} v(s)^i\overline{v(s)}^jZ^{:2-i,1-j:}(s) \right\| _{\mathcal{B}_{6p,\infty}^{\epsilon}}  \mathrm{d} s \right)^3 \mathrm{d} t \\
			\lesssim &\, \left\|\underline{Z}\right\|_{\alpha,\alpha^{\prime},T}^{p_0} \left( \int_{0}^{T} s^{-(\frac{\alpha+\epsilon}{2}+1-\frac{1}{6p})} \mathrm{d} s\right) ^{3}  \left( {t_0}^{-3(\alpha^{\prime}+ \frac{1}{3})} \int_{\frac{t_0}{2}}^{T}   \left\|  \nabla v(s) \right\| _{L^{2}} \mathrm{d} s + {t_0}^{-3(\alpha^{\prime}+ \frac{1}{2})}   \right)  \\
			\lesssim &\, \left\|\underline{Z}\right\|_{\alpha,\alpha^{\prime},T}^{p_0} {t_0}^{-3(\alpha^{\prime}+ \frac{1}{2})}  \left( \int_{\frac{t_0}{2}}^{T}   \left\|  \nabla v(s) \right\| _{L^{2}}^2 \mathrm{d} s \right) ^{\frac12}  \\
			\lesssim &\,\left\|\underline{Z}\right\|_{\alpha,\alpha^{\prime},T}^{p_0} {t_0}^{-(3\alpha^{\prime}+2)}   ,
		\end{align*}
		where we use the fact that $0<\frac{\alpha+\epsilon}{2}+1-\frac{1}{6p}<1$ since $\alpha+\epsilon<\frac{1}{3p}$ for sufficiently small $\alpha$ and $\epsilon< \epsilon_0 = \frac{(p-1)(8p-1)}{3p(p+1)}  < \frac{1}{3p}$. 
		For $i+j=2$,
		\begin{align*}
			\left\| P_{t-s} v(s)^i\overline{v(s)}^jZ^{:2-i,1-j:}(s) \right\| _{\mathcal{B}_{6p,\infty}^{\epsilon}} 
			&\lesssim 	\left\| P_{t-s} v(s)^i\overline{v(s)}^jZ^{:2-i,1-j:}(s) \right\| _{\mathcal{B}_{1,1}^{\epsilon+2(1-\frac{1}{6p})}} \\
			&\lesssim (t-s)^{-(\frac{\alpha+\epsilon}{2}+1-\frac{1}{6p})}\left\| v(s)^i\overline{v(s)}^jZ^{:2-i,1-j:}(s) \right\| _{\mathcal{B}_{1,1}^{-\alpha}}
			\\
			&\lesssim (t-s)^{-(\frac{\alpha+\epsilon}{2}+1-\frac{1}{6p})}\left\| v(s)^i\overline{v(s)}^j \right\| _{\mathcal{B}_{1,1}^{2/3}} \left\| Z^{:2-i,1-j:}(s) \right\| _{\mathcal{C}^{-\alpha}}\\
			&\lesssim (t-s)^{-(\frac{\alpha+\epsilon}{2}+1-\frac{1}{6p})}  \left\| |v(s)|^2 \right\| _{\mathcal{B}_{1,1}^{2/3}} \left\|\underline{Z}\right\|_{\alpha,\alpha^{\prime},T}^{p_0}.
		\end{align*} 
		By Proposition \ref{gradient estimate}, Cauchy-Schwarz inequality and \eqref{eq: Lp est}, 
		\begin{align*}
			\left\| |v(s)|^2 \right\| _{\mathcal{B}_{1,1}^{2/3}}
			&\lesssim \left\| |v(s)|^2\right\| _{L^1}^{1/3}  \left\| \overline{v(s)}\nabla v(s)  \right\| _{L^1}^{2/3}  + \left\| |v(s)|^2 \right\| _{L^1}\\
			&\lesssim    \left\| v(s)\right\| _{L^{2}}^{4/3} \left\| \nabla v(s)  \right\| _{L^2}^{2/3}  + \left\| v(s) \right\| _{L^{2}}^{2}\\
			&\lesssim s^{-2/3} \left\| \nabla v(s)  \right\| _{L^2}^{2/3} + s^{-1}.
		\end{align*}
		Then from Young's convolution inequality and \eqref{eq: Lp+2m est},
		\begin{align*}
			&	\int_{t_0}^{T} \left( \int_{\frac{t_0}{2}}^{t} 	\left\| P_{t-s} v(s)^i\overline{v(s)}^jZ^{:m+1-i,m-j:}(s) \right\| _{\mathcal{B}_{6p,\infty}^{\epsilon}}  \mathrm{d} s \right)^3 \mathrm{d} t \\
			\lesssim &\, \left\|\underline{Z}\right\|_{\alpha,\alpha^{\prime},T}^{p_0}  \left( \int_{0}^{T} s^{-(\frac{\alpha+\epsilon}{2}+1-\frac{1}{6p})} \mathrm{d} s\right) ^3  \left( {t_0}^{-2}\int_{\frac{t_0}{2}}^{T}   \left\|  \nabla v(s) \right\| _{L^{2}}^2 \mathrm{d} s + {t_0}^{-3} \right)  \\
			\lesssim &\,\left\|\underline{Z}\right\|_{\alpha,\alpha^{\prime},T}^{p_0} {t_0}^{- 3}   ,
		\end{align*}
		where we still use the fact that $0<\frac{\alpha+\epsilon}{2}+1-\frac{1}{6p}<1$.
	\end{proof}
	
	\begin{remark}\label{rmk: no improved est for m>1}
		For $m>1$, when handle $ v(s)^{2m}\overline{v(s)}^m $ using the similar argument in \eqref{eq: reason for m=1}, the condition
		$ q(\frac{\epsilon}{2}+\frac{1}{p_1} -\frac{1}{2p(2m+1)} )<1 $ is necessary, where $\epsilon>0$, $p_1=\frac{2 m+2p}{2m+1}$ and $\frac1q= 1+\frac{1}{2m+1}- \frac{1}{p_1}$. This condition is equivalent to 
		\begin{align*}
			0<\frac{\epsilon}{2}< \frac{1}{q} -\frac{1}{p_1} +\frac{1}{2p(2m+1)}  =\frac{(p-m)(4(m +1)p-1)}{2p(2 m+1)(p+m)}  .
		\end{align*} 
		However, only for $m=1$, we can find $1<p<2$ such that $(p-m)(4(m +1)p-1)>0$ and then $\frac{(p-m)(4(m +1)p-1)}{2p(2 m+1)(p+m)}>0$. This is why we can only make such an improved estimation for the case where $m=1$.
	\end{remark}

	\begin{proposition}\label{prop: B-epsilon est}
		Let $T>0$ and $v\in C\left( \left( 0,T\right]; \mathcal{C}^{\beta}\right) $ be a mild solution to \eqref{remainder} with $m=1$. Then there exist $\epsilon_0>0$, $p>1$, $p_0>0$ and $\kappa>0$, such that for any $0<\epsilon<\epsilon_0$ and any $0<t_0\leq T$,   
		\begin{equation*}
			\sup_{t_0\leq t\leq T}\left\| v(t)\right\| _{\mathcal{B}_{2p,\infty}^{\epsilon}} \lesssim \left\| \underline{Z}\right\|_{\alpha,\alpha^{\prime},T}^{p_0}  t_0^{-\kappa}.
		\end{equation*}
	\end{proposition}
	\begin{proof}
		A solution $v\in C\left( \left( 0,T\right]; \mathcal{C}^{\beta}\right)$ satisfies
		\begin{align*}
			v(t)= P_{t-\frac{t_0}{2}}v\Big(\frac{t_0}{2}\Big) +\int_{\frac{t_0}{2}}^{t} P_{t-s} \Psi\left(v(s),\underline{Z}(s)\right)\mathrm{d}s, \quad t_0\leq t\leq T.
		\end{align*}
		By Proposition \ref{priori estimate}, we first take $1<p<3/2$ (but close to 1 since $\mu>0$) such that \eqref{eq: Lp est} and \eqref{eq: Lp+2m est} hold. Then we let $0<\epsilon_0=\frac{(p-1)(2p-1)}{p(p+1)}<\frac12$. For any $0<\epsilon<\epsilon_0$, we estimate the $L^{\infty}([t_0,T]; \mathcal{B}_{2p,\infty}^{\epsilon})$-norm of each term on the right hand side.
		
		For the initial data part, by Proposition \ref{heat kernel smoothing} and \ref{L^p embedding}, and estimate \eqref{eq: Lp est} with $m=1$, there exists $p_0>0$ such that 
		\begin{align*}
			\left\| P_{t-\frac{t_0}{2}}v\Big(\frac{t_0}{2}\Big)\right\| _{\mathcal{B}_{2p,\infty}^{\epsilon}}
			&\lesssim  \Big(t-\frac{t_0}{2}\Big)^{-\frac{\epsilon}{2}}\left\| v\Big(\frac{t_0}{2}\Big)\right\| _{\mathcal{B}_{2p,\infty}^{0}}\lesssim \Big(t-\frac{t_0}{2}\Big)^{-\frac{\epsilon}{2}}\left\| v\Big(\frac{t_0}{2}\Big)\right\| _{L^{2p}}\\
			& \lesssim \left\| \underline{Z}\right\|_{\alpha,\alpha^{\prime},T}^{p_0} t_0^{-\frac{1}{2}} \Big(t-\frac{t_0}{2}\Big)^{-\frac{\epsilon}{2}}.
		\end{align*}
		Note that the exponent $p_0>0$ will vary from line to line in this proof. Then
		\begin{align*}
			\sup_{t_0\leq t\leq T} 	 \left\| P_{t-\frac{t_0}{2}}v\Big(\frac{t_0}{2}\Big)\right\| _{\mathcal{B}_{2p,\infty}^{\epsilon}}
			& \lesssim \left\| \underline{Z}\right\|_{\alpha,\alpha^{\prime},T}^{p_0} t_0^{-\frac{1}{2}-\frac{\epsilon}{2}}.
		\end{align*}
		
		Recall that $\Psi(v,\underline{Z})=-\nu v^{2} \overline{v}-\nu\sum_{0\leq i+j<3}\binom{2}{i}\binom{1}{j}v^i\overline{v}^jZ^{:2-i,1-j:} +\tau(v+Z)$ in \eqref{Psi}. Let $p_1=\frac{2 m+2p}{2m+1}=\frac{2+2p}{3}$. By Proposition \ref{heat kernel smoothing} and \ref{multiplicative structure}, 
		\begin{align*}
			\left\| P_{t-s} v(s)^{2}\overline{v(s)}\right\| _{\mathcal{B}_{2p,\infty}^{\epsilon}}
			&\lesssim 	\left\| v(s)^{2}\overline{v(s)}\right\| _{\mathcal{B}_{2p,\infty}^{\epsilon}}\lesssim 	\left\| v(s)\right\| _{\mathcal{B}_{6p,\infty}^{\epsilon}}^3.
		\end{align*}
		Combining with Proposition \ref{prop: B-epsilon-6p est}, there exists $\kappa>0$ such that
		\begin{align}
			\sup_{t_0\leq t\leq T} \int_{\frac{t_0}{2}}^{t} 	\left\| P_{t-s} v(s)^{2}\overline{v(s)}\right\| _{\mathcal{B}_{2p,\infty}^{\epsilon}} \mathrm{d} s \nonumber
			\lesssim  \int_{\frac{t_0}{2}}^{T} 	\left\| v(s)\right\| _{\mathcal{B}_{6p,\infty}^{\epsilon}}^3 \mathrm{d} s \lesssim  \left\| \underline{Z}\right\|_{\alpha,\alpha^{\prime},T}^{p_0} t_0 ^{-\kappa} .
		\end{align}
		Using Proposition \ref{heat kernel smoothing} and \ref{L^p embedding}, Young's convolution inequality and \eqref{eq: Lp est} with $m=1$, we know that 
		\begin{align*}
			\sup_{t_0\leq t\leq T} \int_{\frac{t_0}{2}}^{t} 	\left\| P_{t-s} v(s)\right\| _{\mathcal{B}_{2p,\infty}^{\epsilon}} \mathrm{d} s
			&\lesssim  \int_{0}^{T} s^{-\frac{\epsilon}{2}} \mathrm{d} s \sup_{\frac{t_0}{2}\leq t\leq T}   \left\|  {v(t)}\right\| _{L^{2p}}  \\
			&\lesssim   \left\| \underline{Z}\right\|^{p_0}_{\alpha,\alpha^{\prime},T} t_0^{-\frac{1}{2}},
		\end{align*}
		where we use the fact that $0<\epsilon<\epsilon_0< 1/2$. Similarly, by \eqref{norm of Z} and taking $\alpha$ small enough such that $\alpha+\epsilon<2$, we have
		\begin{align*}
			\sup_{t_0\leq t\leq T} \int_{\frac{t_0}{2}}^{t} \left\| P_{t-s} Z(s)\right\| _{\mathcal{B}_{2p,\infty}^{\epsilon}} \mathrm{d} s \lesssim \sup_{t_0\leq t\leq T} \int_{\frac{t_0}{2}}^{t} (t-s)^{-\frac{\alpha+\epsilon}{2}}\left\| Z(s)\right\| _{\mathcal{C}^{-\alpha}} \mathrm{d} s \lesssim \left\| \underline{Z}\right\|_{\alpha,\alpha^{\prime},T}^{p_0} .
		\end{align*}
		For $i+j=0$, applying Proposition \ref{heat kernel smoothing}, \eqref{norm of Z} and Young's convolution inequality,
		\begin{align*}
			\sup_{t_0\leq t\leq T} \int_{\frac{t_0}{2}}^{t}	\left\| P_{t-s} Z^{:2,1:}(s) \right\| _{\mathcal{B}_{2p,\infty}^{\epsilon}} \mathrm{d} s 
			&\lesssim \left\| \underline{Z}\right\|_{\alpha,\alpha^{\prime},T}^{p_0}\sup_{t_0\leq t\leq T} \int_{\frac{t_0}{2}}^{t} (t-s)^{-\frac{\alpha+\epsilon}{2}} s^{-2 \alpha^{'}} \mathrm{d} s  \\
			&\lesssim \left\| \underline{Z}\right\|_{\alpha,\alpha^{\prime},T}^{p_0} {t_0}^{-2 \alpha^{'}} .
		\end{align*}
		For $i+j=1$, from Proposition \ref{embedding}, \ref{regularity embedding}, \ref{heat kernel smoothing}, \ref{multiplicative structure} and \eqref{norm of Z}, 
		\begin{align*}
			\left\| P_{t-s} v(s)^i\overline{v(s)}^jZ^{:2-i,1-j:}(s) \right\| _{\mathcal{B}_{2p,\infty}^{\epsilon}} 
			&\lesssim 	\left\| P_{t-s} v(s)^i\overline{v(s)}^jZ^{:2-i,1-j:}(s) \right\| _{\mathcal{B}_{1,1}^{\epsilon+2(1-\frac{1}{2p})}} \\
			&\lesssim (t-s)^{-(\frac{\alpha+\epsilon}{2}+1-\frac{1}{2p})}\left\| v(s)^i\overline{v(s)}^jZ^{:2-i,1-j:}(s) \right\| _{\mathcal{B}_{1,1}^{-\alpha}}
			\\
			&\lesssim (t-s)^{-(\frac{\alpha+\epsilon}{2}+1-\frac{1}{2p})}\left\| v(s)^i\overline{v(s)}^j \right\| _{\mathcal{B}_{1,1}^{2\alpha}} \left\| Z^{:2-i,1-j:}(s) \right\| _{\mathcal{C}^{-\alpha}}\\
			&\lesssim (t-s)^{-(\frac{\alpha+\epsilon}{2}+1-\frac{1}{2p})} s^{-\alpha^{\prime}} \left\| v(s) \right\| _{\mathcal{B}_{1,1}^{2\alpha}} \left\|\underline{Z}\right\|_{\alpha,\alpha^{\prime},T}^{p_0}.
		\end{align*} 
		By Proposition \ref{gradient estimate}, H\"older inequality and \eqref{eq: Lp est} with $m=1$,
		\begin{align*}
			\left\| v(s) \right\| _{\mathcal{B}_{1,1}^{2\alpha}}
			&\lesssim \left\| v(s)\right\| _{L^1}^{1-2\alpha}  \left\| \nabla v(s)  \right\| _{L^1}^{2\alpha}  + \left\| v(s) \right\| _{L^1}\\
			&\lesssim \left\| v(s)\right\| _{L^{2}}^{1-2\alpha}  \left\| \nabla v(s)  \right\| _{L^2}^{2\alpha}  + \left\| v(s) \right\| _{L^{2}}\\
			&\lesssim s^{-\frac{1-2\alpha}{2}} \left\| \nabla v(s)  \right\| _{L^2}^{2\alpha} + s^{-\frac{1}{2}}.
		\end{align*}
		Then according to Young's convolution inequality and \eqref{eq: Lp+2m est} with $m=1$,
		\begin{align*}
			&\sup_{t_0\leq t\leq T} \int_{\frac{t_0}{2}}^{t} 	\left\| P_{t-s} v(s)^i\overline{v(s)}^jZ^{:2-i,1-j:}(s) \right\| _{\mathcal{B}_{2p,\infty}^{\epsilon}}  \mathrm{d} s\\
			\lesssim &\, \left\|\underline{Z}\right\|_{\alpha,\alpha^{\prime},T}^{p_0} {t_0}^{-(\alpha^{\prime}+ \frac{1-2\alpha}{2})} \left( \int_{0}^{T} s^{-\frac{1}{1-\alpha}(\frac{\alpha+\epsilon}{2}+1-\frac{1}{2p})} \mathrm{d} s\right) ^{1-\alpha} \left( \int_{\frac{t_0}{2}}^{T}   \left\|  \nabla v(s) \right\| _{L^{2}}^2 \mathrm{d} s\right) ^{\alpha} \\
			&+ \left\|\underline{Z}\right\|_{\alpha,\alpha^{\prime},T}^{p_0} {t_0}^{-(\alpha^{\prime}+ \frac{1}{2})} \int_{0}^{T} s^{-(\frac{\alpha+\epsilon}{2}+1-\frac{1}{2p})} \mathrm{d} s  \\
			\lesssim &\,\left\|\underline{Z}\right\|_{\alpha,\alpha^{\prime},T}^{p_0} {t_0}^{-(\alpha^{\prime}+ \frac{1}{2})}   ,
		\end{align*}
		where we use the fact that $0<\frac{1}{1-\alpha}(\frac{\alpha+\epsilon}{2}+1-\frac{1}{2p})<1$ since $\frac{\alpha+\epsilon}{2}+1-\frac{1}{2p}<1-\alpha$ for sufficiently small $\alpha$ and $\epsilon< \epsilon_0 = \frac{(p-1)(2p-1)}{p(p+1)} < \frac{1}{p}$. (This is from $(p-1)(2p-1)< p+1\Leftrightarrow p(p-2)< 0 \Leftarrow 1<p<2$.) 
		For $i+j=2$, by Proposition \ref{embedding}, \ref{regularity embedding}, \ref{heat kernel smoothing}, \ref{multiplicative structure} and \eqref{norm of Z}, 
		\begin{align*}
			\left\| P_{t-s} v(s)^i\overline{v(s)}^jZ^{:2-i,1-j:}(s) \right\| _{\mathcal{B}_{2p,\infty}^{\epsilon}} 
			&\lesssim 	\left\| P_{t-s} v(s)^i\overline{v(s)}^jZ^{:2-i,1-j:}(s) \right\| _{\mathcal{B}_{1,1}^{\epsilon+2(1-\frac{1}{2p})}} \\
			&\lesssim (t-s)^{-(\frac{\alpha+\epsilon}{2}+1-\frac{1}{2p})}\left\| v(s)^i\overline{v(s)}^jZ^{:2-i,1-j:}(s) \right\| _{\mathcal{B}_{1,1}^{-\alpha}}
			\\
			&\lesssim (t-s)^{-(\frac{\alpha+\epsilon}{2}+1-\frac{1}{2p})}\left\| v(s)^i\overline{v(s)}^j \right\| _{\mathcal{B}_{1,1}^{2\alpha}} \left\| Z^{:2-i,1-j:}(s) \right\| _{\mathcal{C}^{-\alpha}}\\
			&\lesssim (t-s)^{-(\frac{\alpha+\epsilon}{2}+1-\frac{1}{2p})}  \left\| |v(s)|^2 \right\| _{\mathcal{B}_{1,1}^{2\alpha}} \left\|\underline{Z}\right\|_{\alpha,\alpha^{\prime},T}^{p_0}.
		\end{align*} 
		Using Proposition \ref{gradient estimate}, Cauchy-Schwarz inequality and \eqref{eq: Lp est} with $m=1$, 
		\begin{align*}
			\left\| |v(s)|^2 \right\| _{\mathcal{B}_{1,1}^{2\alpha}}
			&\lesssim \left\| |v(s)|^2\right\| _{L^1}^{1-2\alpha}  \left\| \overline{v(s)}\nabla v(s)  \right\| _{L^1}^{2\alpha}  + \left\| |v(s)|^2 \right\| _{L^1}\\
			&\lesssim \left\| v(s)\right\| _{L^{2}}^{2(1-2\alpha)}   \left\| v(s)\right\| _{L^{2}}^{2\alpha} \left\| \nabla v(s)  \right\| _{L^2}^{2\alpha}  + \left\| v(s) \right\| _{L^{2}}^{2}\\
			&\lesssim s^{-(1-\alpha)} \left\| \nabla v(s)  \right\| _{L^2}^{2\alpha} + s^{-1}.
		\end{align*}
		Then by Young's convolution inequality,
		\begin{align*}
			&\sup_{t_0\leq t\leq T} \int_{\frac{t_0}{2}}^{t} 	\left\| P_{t-s} v(s)^i\overline{v(s)}^jZ^{:m+1-i,m-j:}(s) \right\| _{\mathcal{B}_{2p,\infty}^{\epsilon}}  \mathrm{d} s\\
			\lesssim &\, \left\|\underline{Z}\right\|_{\alpha,\alpha^{\prime},T}^{p_0} {t_0}^{-(1-\alpha)} \left( \int_{0}^{T} s^{-\frac{1}{1-\alpha}(\frac{\alpha+\epsilon}{2}+1-\frac{1}{2p})} \mathrm{d} s\right) ^{1-\alpha} \left( \int_{\frac{t_0}{2}}^{T}   \left\|  \nabla v(s) \right\| _{L^{2}}^2 \mathrm{d} s\right) ^{\alpha} \\
			&+ \left\|\underline{Z}\right\|_{\alpha,\alpha^{\prime},T}^{p_0} {t_0}^{-1} \int_{0}^{T} s^{-(\frac{\alpha+\epsilon}{2}+1-\frac{1}{2p})} \mathrm{d} s  \\
			\lesssim &\,\left\|\underline{Z}\right\|_{\alpha,\alpha^{\prime},T}^{p_0} {t_0}^{- 1}   ,
		\end{align*}
		where we still use the fact that $0<\frac{1}{1-\alpha}(\frac{\alpha+\epsilon}{2}+1-\frac{1}{2p})<1$. 
	\end{proof}

	\begin{proposition}\label{prop: B-2epsilon est}
		Let $T>0$ and $v\in C\left( \left( 0,T\right]; \mathcal{C}^{\beta}\right) $ be a mild solution to \eqref{remainder} with $m=1$. Assume there exist $p>1$, $0<\epsilon_0\leq 2$, $p_0>0$ and $\kappa>0$, such that for any $0<\epsilon<\epsilon_0$ and any $0<t_0\leq T$,   
		\begin{equation}\label{eq: bootstrap}
			\int_{\frac{t_0}{2}}^{T}\left\| v(t)\right\|^3 _{\mathcal{B}_{6p,\infty}^{\epsilon}} \mathrm{d}t + \sup_{\frac{t_0}{2}\leq t\leq T}\left\| v(t)\right\| _{\mathcal{B}_{2p,\infty}^{\epsilon}} \lesssim \left\| \underline{Z}\right\|_{\alpha,\alpha^{\prime},T}^{p_0}  t_0^{-\kappa}.
		\end{equation}
		Then for $\epsilon_1$ satisfying $0<\epsilon_1 < \frac53(1-\frac{1}{p})$ and $\epsilon+\epsilon_1<\frac53-\frac{2}{3p}$, there exists $\kappa_1>\kappa$ such that 
		\begin{equation*}
			\int_{t_0}^{T}\left\| v(t)\right\|^3 _{\mathcal{B}_{6p,\infty}^{\epsilon+ \epsilon_1}} \mathrm{d}t + \sup_{t_0\leq t\leq T}\left\| v(t)\right\| _{\mathcal{B}_{2p,\infty}^{\epsilon+ \epsilon_1}} \lesssim \left\| \underline{Z}\right\|_{\alpha,\alpha^{\prime},T}^{p_0}  t_0^{-\kappa_1}.
		\end{equation*}
	\end{proposition}
	
	\begin{proof}
		A solution $v\in C\left( \left( 0,T\right]; \mathcal{C}^{\beta}\right)$ satisfies
		\begin{align*}
			v(t)= P_{t-\frac{t_0}{2}}v\Big(\frac{t_0}{2}\Big) +\int_{\frac{t_0}{2}}^{t} P_{t-s} \Psi\left(v(s),\underline{Z}(s)\right)\mathrm{d}s, \quad t_0\leq t\leq T.
		\end{align*}

		We first estimate the $L^{3}([t_0,T]; \mathcal{B}_{6p,\infty}^{\epsilon+\epsilon_1})$-norm. For the initial data part, by Proposition \ref{regularity embedding} and \ref{heat kernel smoothing}, and assumption \eqref{eq: bootstrap},
		\begin{align*}
			\left\| P_{t-\frac{t_0}{2}}v\Big(\frac{t_0}{2}\Big)\right\| _{\mathcal{B}_{6p,\infty}^{\epsilon+\epsilon_1}} 
			&\lesssim 	\left\| P_{t-\frac{t_0}{2}}v\Big(\frac{t_0}{2}\Big)\right\| _{\mathcal{B}_{2p,\infty}^{\epsilon+\epsilon_1+\frac{2}{3p}}}
			\lesssim  \Big(t-\frac{t_0}{2}\Big)^{-(\frac{\epsilon_1}{2}+\frac{1}{3p})}\left\| v\Big(\frac{t_0}{2}\Big)\right\| _{\mathcal{B}_{2p,\infty}^{\epsilon}}\\
			&\lesssim \left\| \underline{Z}\right\|_{\alpha,\alpha^{\prime},T}^{p_0} t_0^{-\kappa} \Big(t-\frac{t_0}{2}\Big)^{-(\frac{\epsilon_1}{2}+\frac{1}{3p})},
		\end{align*}
		where the exponent $p_0>0$ will vary from line to line in this proof. Then
		\begin{align*}
			\int_{t_0}^{T} \left\| P_{t-\frac{t_0}{2}}v\Big(\frac{t_0}{2}\Big)\right\|^3 _{\mathcal{B}_{6p,\infty}^{\epsilon+\epsilon_1}} \mathrm{d}t
			& \lesssim \left\| \underline{Z}\right\|_{\alpha,\alpha^{\prime},T}^{p_0} t_0^{-3(\frac{\epsilon_1}{2}+\frac{1}{3p}+\kappa)}.
		\end{align*} 
		
		Recall that $\Psi(v,\underline{Z})=-\nu v^{2} \overline{v}-\nu\sum_{0\leq i+j<3}\binom{2}{i}\binom{1}{j}v^i\overline{v}^jZ^{:2-i,1-j:} +\tau(v+Z)$ in \eqref{Psi}. Applying Proposition \ref{regularity embedding}, \ref{heat kernel smoothing} and \ref{multiplicative structure},
		\begin{align*}
			&\left\| P_{t-s} v(s)^{2}\overline{v(s)}\right\| _{\mathcal{B}_{6p,\infty}^{\epsilon+ \epsilon_1}}
			\lesssim \left\| P_{t-s} v(s)^{2}\overline{v(s)}\right\|_{\mathcal{B}_{p,\infty}^{\epsilon+ \epsilon_1+\frac{5}{3p}}}\lesssim (t-s)^{-(\frac{\epsilon_1}{2}+\frac{5}{6p} )}  \left\|  v(s)^{2}\overline{v(s)}\right\| _{\mathcal{B}_{p,\infty}^{\epsilon }}\\
			\lesssim &\, (t-s)^{-(\frac{\epsilon_1}{2}+\frac{5}{6p} )}  \left\|  {v(s)}^{2}\right\| _{\mathcal{B}_{2p,\infty}^{\epsilon }} \left\|  {v(s)}\right\| _{\mathcal{B}_{2p,\infty}^{\epsilon }}.
		\end{align*}
		By Proposition \ref{multiplicative structure} and \ref{interpolation},
		\begin{align}\label{eq: v2 est}
			\left\|  {v(s)}^{2}\right\| _{\mathcal{B}_{2p,\infty}^{\epsilon }} \lesssim \left\|  {v(s)}\right\| _{\mathcal{B}_{4p,\infty}^{\epsilon }}^{2} \lesssim \left\|  {v(s)}\right\| _{\mathcal{B}_{6p,\infty}^{\epsilon }}^{\frac32}  \left\|  {v(s)}\right\| _{\mathcal{B}_{2p,\infty}^{\epsilon }}^{\frac12}.
		\end{align}
		Then using assumption \eqref{eq: bootstrap},
		\begin{align*}
			\left\| P_{t-s} v(s)^{2}\overline{v(s)}\right\| _{\mathcal{B}_{6p,\infty}^{\epsilon+ \epsilon_1}}
			&\lesssim (t-s)^{-(\frac{\epsilon_1}{2}+\frac{5}{6p} )}   \left\|  {v(s)}\right\| _{\mathcal{B}_{6p,\infty}^{\epsilon }}^{\frac32}  \left\|  {v(s)}\right\| _{\mathcal{B}_{2p,\infty}^{\epsilon }}^{\frac32}\\
			&\lesssim s^{-\frac{3\kappa}{2}} (t-s)^{-(\frac{\epsilon_1}{2}+\frac{5}{6p} )}   \left\|  {v(s)}\right\| _{\mathcal{B}_{6p,\infty}^{\epsilon }}^{\frac32}.
		\end{align*}
		According to Young's convolution inequality and assumption \eqref{eq: bootstrap}, 
		\begin{align*}
			\int_{t_0}^{T} \left( \int_{\frac{t_0}{2}}^{t} 	\left\| P_{t-s} v(s)^{2}\overline{v(s)}\right\| _{\mathcal{B}_{6p,\infty}^{\epsilon+\epsilon_1}} \mathrm{d} s\right)^3 \mathrm{d} t
			&\lesssim t_0^{-\frac{9\kappa}{2}} 		\int_{t_0}^{T} \left( \int_{\frac{t_0}{2}}^{t} (t-s)^{-(\frac{\epsilon_1}{2}+\frac{5}{6p} )}  \left\|  {v(s)}\right\| _{\mathcal{B}_{6p,\infty}^{\epsilon }}^{\frac32} \mathrm{d} s\right)^3 \mathrm{d} t \\
			&\lesssim t_0^{-\frac{9\kappa}{2}}   \left( \int_{0}^{T} s^{-\frac65(\frac{\epsilon_1}{2}+\frac{5}{6p} )} \mathrm{d} s\right) ^{\frac{5}{2}} \left( \int_{\frac{t_0}{2}}^{T}  \left\|  {v(s)}\right\| _{\mathcal{B}_{6p,\infty}^{\epsilon }}^3 \mathrm{d} s\right) ^{\frac{3}{2}}\\
			&\lesssim  \left\| \underline{Z}\right\|_{\alpha,\alpha^{\prime},T}^{p_0} t_0 ^{-6\kappa } ,
		\end{align*}
		where we use the fact that $0< \frac65(\frac{\epsilon_1}{2}+\frac{5}{6p} ) = \frac35 \epsilon_1 + \frac{1}{p} <1 $ since 
		$\epsilon_1 < \frac53(1-\frac{1}{p})$. Using the similar argument, we know that 
		\begin{align*}
			\int_{t_0}^{T} \left( \int_{\frac{t_0}{2}}^{t} 	\left\| P_{t-s} v(s)\right\| _{\mathcal{B}_{6p,\infty}^{\epsilon+\epsilon_1}} \mathrm{d} s \right)^3 \mathrm{d} t &\lesssim \left(\int_{0}^{T} s^{-\frac{\epsilon_1}{2}} \mathrm{d} s\right)^3  \int_{\frac{t_0}{2}}^{T}  \left\|v(t)\right\|_{\mathcal{B}_{6p,\infty}^{\epsilon}}^3 \mathrm{d} t \\
			&\lesssim   \left\| \underline{Z}\right\|^{p_0}_{\alpha,\alpha^{\prime},T} t_0^{-\kappa},
		\end{align*}
		where we use the fact that $0<\epsilon_1
		< 2 $. We take $\alpha$ small enough such that $\alpha+\epsilon+\epsilon_1<2$, then
		\begin{align*}
			\int_{t_0}^{T} \left( \int_{\frac{t_0}{2}}^{t} \left\| P_{t-s} Z(s)\right\| _{\mathcal{B}_{6p,\infty}^{\epsilon+\epsilon_1}} \mathrm{d} s \right)^3 \mathrm{d} t  
			&\lesssim 	\int_{t_0}^{T} \left( \int_{\frac{t_0}{2}}^{t} (t-s)^{-\frac{\alpha+\epsilon+\epsilon_1}{2}}\left\| Z(s)\right\| _{\mathcal{C}^{-\alpha}} \mathrm{d} s \right)^3 \mathrm{d} t \\
			&\lesssim 	\left( \int_{0}^{T} s^{-\frac{\alpha+\epsilon+\epsilon_1}{2}}  \mathrm{d} s\right)^3 \left\| \underline{Z}\right\|_{\alpha,\alpha^{\prime},T}^{p_0} \lesssim \left\| \underline{Z}\right\|_{\alpha,\alpha^{\prime},T}^{p_0}.
		\end{align*}
		For $i+j=0$, by Proposition \ref{embedding} and \ref{heat kernel smoothing}, Young's convolution inequality and \eqref{norm of Z},
		\begin{align*}
			\int_{t_0}^{T} \left(\int_{\frac{t_0}{2}}^{t}	\left\| P_{t-s} Z^{:2,1:}(s) \right\| _{\mathcal{B}_{6p,\infty}^{\epsilon+\epsilon_1}} \mathrm{d} s \right)^3 \mathrm{d} t
			&\lesssim \left\| \underline{Z}\right\|_{\alpha,\alpha^{\prime},T}^{p_0}	\int_{t_0}^{T} \left( \int_{\frac{t_0}{2}}^{t} (t-s)^{-\frac{\alpha+\epsilon+\epsilon_1}{2}} s^{-2 \alpha^{'}} \mathrm{d} s \right)^3 \mathrm{d} t \\
			&\lesssim \left\| \underline{Z}\right\|_{\alpha,\alpha^{\prime},T}^{p_0} 	\left( \int_{0}^{T} s^{-\frac{\alpha+\epsilon+\epsilon_1}{2}}  \mathrm{d} s\right)^3 \int_{\frac{t_0}{2}}^{T} s^{-6\alpha^{\prime}}  \mathrm{d} s \\
			&\lesssim \left\| \underline{Z}\right\|_{\alpha,\alpha^{\prime},T}^{p_0}  {t_0}^{-6 \alpha^{'}} .
		\end{align*}
		We take $\alpha$ small enough such that $\alpha<\epsilon$. Then for $i+j=1$,  by Proposition \ref{regularity embedding}, \ref{heat kernel smoothing} and \ref{multiplicative structure}, \eqref{norm of Z} and assumption \eqref{eq: bootstrap},
		\begin{align*}
			\left\| P_{t-s} v(s)^i\overline{v(s)}^jZ^{:2-i,1-j:}(s) \right\| _{\mathcal{B}_{6p,\infty}^{\epsilon+\epsilon_1}} &\lesssim  	\left\| P_{t-s} v(s)^i\overline{v(s)}^jZ^{:2-i,1-j:}(s) \right\| _{\mathcal{B}_{2p,\infty}^{\epsilon+\epsilon_1+\frac{2}{3p}}}\\
			&\lesssim (t-s)^{-(\frac{\alpha+\epsilon+\epsilon_1}{2}+\frac{1}{3p})}\left\| v(s)^i\overline{v(s)}^jZ^{:2-i,1-j:}(s) \right\| _{\mathcal{B}_{2p,\infty}^{-\alpha}}
			\\
			&\lesssim (t-s)^{-(\frac{\alpha+\epsilon+\epsilon_1}{2}+\frac{1}{3p})}\left\| v(s)^i\overline{v(s)}^j \right\| _{\mathcal{B}_{2p,\infty}^{\epsilon}} \left\| Z^{:2-i,1-j:}(s) \right\| _{\mathcal{C}^{-\alpha}}\\
			&\lesssim(t-s)^{-(\frac{\alpha+\epsilon+\epsilon_1}{2}+\frac{1}{3p})} s^{-\alpha^{\prime}} \left\| v(s) \right\| _{\mathcal{B}_{2p,\infty}^{\epsilon}} \left\|\underline{Z}\right\|_{\alpha,\alpha^{\prime},T}^{p_0}\\
			&\lesssim (t-s)^{-(\frac{\alpha+\epsilon+\epsilon_1}{2}+\frac{1}{3p})} s^{-(\alpha^{\prime}+\kappa)}  \left\|\underline{Z}\right\|_{\alpha,\alpha^{\prime},T}^{p_0}.
		\end{align*} 
		Then by Young's convolution inequality,
		\begin{align*}
			&	\int_{t_0}^{T} \left(\int_{\frac{t_0}{2}}^{t}	\left\| P_{t-s} v(s)^i\overline{v(s)}^jZ^{:2-i,1-j:}(s) \right\| _{\mathcal{B}_{2p,\infty}^{\epsilon+\epsilon_1}}  \mathrm{d} s \right)^3 \mathrm{d} t \\
			\lesssim &\, \left\|\underline{Z}\right\|_{\alpha,\alpha^{\prime},T}^{p_0} {t_0}^{-3(\alpha^{\prime}+ \kappa)} \left( \int_{0}^{T} s^{-(\frac{\alpha+\epsilon+\epsilon_1}{2}+\frac{1}{3p})} \mathrm{d} s \right) ^3 
			\lesssim \left\|\underline{Z}\right\|_{\alpha,\alpha^{\prime},T}^{p_0} {t_0}^{-3(\alpha^{\prime}+ \kappa)}   ,
		\end{align*}
		where we use the fact that $0<\frac{\alpha+\epsilon+\epsilon_1}{2}+\frac{1}{3p}<1$ since $\alpha$ is sufficiently small and 
		$\epsilon+\epsilon_1<\frac53-\frac{2}{3p}< 2-\frac{2}{3p}$. For $i+j=2$, by Proposition \ref{regularity embedding}, \ref{heat kernel smoothing} and \ref{multiplicative structure}, \eqref{norm of Z}, estimate \eqref{eq: v2 est} and assumption \eqref{eq: bootstrap},
		\begin{align*}
			\left\| P_{t-s} v(s)^i\overline{v(s)}^jZ^{:2-i,1-j:}(s) \right\| _{\mathcal{B}_{6p,\infty}^{\epsilon+\epsilon_1}} 
			&\lesssim 	\left\| P_{t-s} v(s)^i\overline{v(s)}^jZ^{:2-i,1-j:}(s) \right\| _{\mathcal{B}_{2p,\infty}^{\epsilon+\epsilon_1+\frac{2}{3p}}} \\
			&\lesssim (t-s)^{-(\frac{\alpha+ \epsilon+\epsilon_1}{2}+\frac{1}{3p})}\left\| v(s)^i\overline{v(s)}^jZ^{:2-i,1-j:}(s) \right\| _{\mathcal{B}_{2p,\infty}^{-\alpha}}
			\\
			&\lesssim (t-s)^{-(\frac{\alpha+ \epsilon+\epsilon_1}{2}+\frac{1}{3p})} \left\| v(s)^i\overline{v(s)}^j \right\| _{\mathcal{B}_{2p,\infty}^{\epsilon}} \left\| Z^{:2-i,1-j:}(s) \right\| _{\mathcal{C}^{-\alpha}}\\
			&\lesssim (t-s)^{-(\frac{\alpha+ \epsilon+\epsilon_1}{2}+\frac{1}{3p})} \left\|  {v(s)}\right\| _{\mathcal{B}_{6p,\infty}^{\epsilon }}^{\frac32}  \left\|  {v(s)}\right\| _{\mathcal{B}_{2p,\infty}^{\epsilon }}^{\frac12} \left\|\underline{Z}\right\|_{\alpha,\alpha^{\prime},T}^{p_0}\\
			&\lesssim (t-s)^{-(\frac{\alpha+ \epsilon+\epsilon_1}{2}+\frac{1}{3p})} s^{-\frac{\kappa}{2}} \left\|  {v(s)}\right\| _{\mathcal{B}_{6p,\infty}^{\epsilon }}^{\frac32} \left\|\underline{Z}\right\|_{\alpha,\alpha^{\prime},T}^{p_0}.
		\end{align*} 
		Then by Young's convolution inequality and assumption \eqref{eq: bootstrap},
		\begin{align*}
			&	\int_{t_0}^{T} \left( \int_{\frac{t_0}{2}}^{t} 	\left\| P_{t-s} v(s)^i\overline{v(s)}^jZ^{:2-i,1-j:}(s) \right\| _{\mathcal{B}_{2p,\infty}^{\epsilon}}  \mathrm{d} s \right)^3 \mathrm{d} t\\
			\lesssim &\, \left\|\underline{Z}\right\|_{\alpha,\alpha^{\prime},T}^{p_0} {t_0}^{- \frac{3\kappa}{2}} \left(\int_{0}^{T} s^{-\frac65(\frac{\alpha+\epsilon+\epsilon_1}{2}+\frac{1}{3p})} \mathrm{d} s   \right)^{\frac52} \left(\int_{\frac{t_0}{2}}^{T}  \left\|  {v(s)}\right\| _{\mathcal{B}_{6p,\infty}^{\epsilon }}^3 \mathrm{d} s \right)^{\frac32}
			\lesssim \left\|\underline{Z}\right\|_{\alpha,\alpha^{\prime},T}^{p_0} {t_0}^{-3\kappa}   ,
		\end{align*}
		where we use the fact that $0<\frac65(\frac{\alpha+\epsilon+\epsilon_1}{2}+\frac{1}{3p})<1$ since $\alpha$ is sufficiently small and 
		$\epsilon+\epsilon_1<\frac53-\frac{2}{3p}$.

		We next estimate the $L^{\infty}([t_0,T]; \mathcal{B}_{2p,\infty}^{\epsilon+\epsilon_1})$-norm. For the initial data part, by Proposition \ref{heat kernel smoothing} and assumption \eqref{eq: bootstrap},
		\begin{align*}
			\left\| P_{t-\frac{t_0}{2}}v\Big(\frac{t_0}{2}\Big)\right\| _{\mathcal{B}_{2p,\infty}^{\epsilon+\epsilon_1}}
			\lesssim  \Big(t-\frac{t_0}{2}\Big)^{-\frac{\epsilon_1}{2}}\left\| v\Big(\frac{t_0}{2}\Big)\right\| _{\mathcal{B}_{2p,\infty}^{\epsilon}}
			\lesssim \left\| \underline{Z}\right\|_{\alpha,\alpha^{\prime},T}^{p_0} t_0^{-\kappa} \Big(t-\frac{t_0}{2}\Big)^{-\frac{\epsilon_1}{2}},
		\end{align*}
		where the exponent $p_0>0$ will vary from line to line in this proof. Then
		\begin{align*}
			\sup_{t_0\leq t\leq T} 	 \left\| P_{t-\frac{t_0}{2}}v\Big(\frac{t_0}{2}\Big)\right\| _{\mathcal{B}_{2p,\infty}^{\epsilon+\epsilon_1}}
			& \lesssim \left\| \underline{Z}\right\|_{\alpha,\alpha^{\prime},T}^{p_0} t_0^{-\frac{\epsilon_1}{2}-\kappa}.
		\end{align*} 
		
		Recall that $\Psi(v,\underline{Z})=-\nu v^{2} \overline{v}-\nu\sum_{0\leq i+j<3}\binom{2}{i}\binom{1}{j}v^i\overline{v}^jZ^{:2-i,1-j:} +\tau(v+Z)$ in \eqref{Psi}. By Proposition \ref{heat kernel smoothing} and \ref{multiplicative structure}, 
		\begin{align*}
			&\left\| P_{t-s} v(s)^{2}\overline{v(s)}\right\| _{\mathcal{B}_{2p,\infty}^{\epsilon+ \epsilon_1}}
			\lesssim \left\|  v(s)^{2}\overline{v(s)}\right\| _{\mathcal{B}_{2p,\infty}^{\epsilon+ \epsilon_1}} \lesssim \left\|  v(s)\right\| _{\mathcal{B}_{6p,\infty}^{\epsilon+ \epsilon_1}}^3.
		\end{align*}
		Using the $L^{3}([t_0,T]; \mathcal{B}_{6p,\infty}^{\epsilon+\epsilon_1})$-norm estimate above, there exists $\kappa_1>0$ such that
		\begin{align*}
			\sup_{t_0\leq t\leq T} \int_{\frac{t_0}{2}}^{t} 	\left\| P_{t-s} v(s)^{2}\overline{v(s)}\right\| _{\mathcal{B}_{2p,\infty}^{\epsilon}} \mathrm{d} s
			\lesssim \int_{\frac{t_0}{2}}^{T} \left\|  v(s)\right\| _{\mathcal{B}_{6p,\infty}^{\epsilon+ \epsilon_1}}^3 \mathrm{d} s \lesssim  \left\| \underline{Z}\right\|_{\alpha,\alpha^{\prime},T}^{p_0} t_0 ^{-\kappa_1} .
		\end{align*}
		Applying Proposition \ref{heat kernel smoothing}, Young's convolution inequality and assumption \eqref{eq: bootstrap}, we know that 
		\begin{align*}
			\sup_{t_0\leq t\leq T} \int_{\frac{t_0}{2}}^{t} 	\left\| P_{t-s} v(s)\right\| _{\mathcal{B}_{2p,\infty}^{\epsilon+\epsilon_1}} \mathrm{d} s &\lesssim  \int_{0}^{T} s^{-\frac{\epsilon_1}{2}} \mathrm{d} s \sup_{\frac{t_0}{2}\leq t\leq T}   \left\|v(t)\right\|_{\mathcal{B}_{2p,\infty}^{\epsilon}}  \\
			&\lesssim   \left\| \underline{Z}\right\|^{p_0}_{\alpha,\alpha^{\prime},T} t_0^{-\kappa},
		\end{align*}
		where we use the fact that $0<\epsilon_1<\frac53
		<2$. We take $\alpha$ small enough such that $\alpha+\epsilon+\epsilon_1<2$, then
		\begin{align*}
			\sup_{t_0\leq t\leq T} \int_{\frac{t_0}{2}}^{t} \left\| P_{t-s} Z(s)\right\| _{\mathcal{B}_{2p,\infty}^{\epsilon+\epsilon_1}} \mathrm{d} s \lesssim \sup_{t_0\leq t\leq T} \int_{\frac{t_0}{2}}^{t} (t-s)^{-\frac{\alpha+\epsilon+\epsilon_1}{2}}\left\| Z(s)\right\| _{\mathcal{C}^{-\alpha}} \mathrm{d} s \lesssim \left\| \underline{Z}\right\|_{\alpha,\alpha^{\prime},T}^{p_0} .
		\end{align*}
		For $i+j=0$, by Proposition \ref{heat kernel smoothing}, \eqref{norm of Z} and Young's convolution inequality,
		\begin{align*}
			\sup_{t_0\leq t\leq T} \int_{\frac{t_0}{2}}^{t}	\left\| P_{t-s} Z^{:2,1:}(s) \right\| _{\mathcal{B}_{2p,\infty}^{\epsilon+\epsilon_1}} \mathrm{d} s 
			&\lesssim \left\| \underline{Z}\right\|_{\alpha,\alpha^{\prime},T}^{p_0}\sup_{t_0\leq t\leq T} \int_{\frac{t_0}{2}}^{t} (t-s)^{-\frac{\alpha+\epsilon+\epsilon_1}{2}} s^{-2 \alpha^{'}} \mathrm{d} s  \\
			&\lesssim \left\| \underline{Z}\right\|_{\alpha,\alpha^{\prime},T}^{p_0} {t_0}^{-2 \alpha^{'}} .
		\end{align*}
		We take $\alpha$ small enough such that $\alpha<\epsilon$. Then for $i+j=1$, according to Proposition \ref{heat kernel smoothing} and \ref{multiplicative structure}, \eqref{norm of Z} and assumption \eqref{eq: bootstrap},
		\begin{align*}
			\left\| P_{t-s} v(s)^i\overline{v(s)}^jZ^{:2-i,1-j:}(s) \right\| _{\mathcal{B}_{2p,\infty}^{\epsilon+\epsilon_1}} 
			&\lesssim (t-s)^{-\frac{\alpha+ \epsilon+\epsilon_1}{2}}\left\| v(s)^i\overline{v(s)}^jZ^{:2-i,1-j:}(s) \right\| _{\mathcal{B}_{2p,\infty}^{-\alpha}}
			\\
			&\lesssim (t-s)^{-\frac{\alpha+ \epsilon+\epsilon_1}{2}}\left\| v(s)^i\overline{v(s)}^j \right\| _{\mathcal{B}_{2p,\infty}^{\epsilon}} \left\| Z^{:2-i,1-j:}(s) \right\| _{\mathcal{C}^{-\alpha}}\\
			&\lesssim (t-s)^{-\frac{\alpha+ \epsilon+\epsilon_1}{2}} s^{-\alpha^{\prime}} \left\| v(s) \right\| _{\mathcal{B}_{2p,\infty}^{\epsilon}} \left\|\underline{Z}\right\|_{\alpha,\alpha^{\prime},T}^{p_0}\\
			&\lesssim (t-s)^{-\frac{\alpha+ \epsilon+\epsilon_1}{2}} s^{-(\alpha^{\prime}+\kappa)}  \left\|\underline{Z}\right\|_{\alpha,\alpha^{\prime},T}^{p_0}.
		\end{align*} 
		Then by Young's convolution inequality,
		\begin{align*}
			&\sup_{t_0\leq t\leq T} \int_{\frac{t_0}{2}}^{t} 	\left\| P_{t-s} v(s)^i\overline{v(s)}^jZ^{:2-i,1-j:}(s) \right\| _{\mathcal{B}_{2p,\infty}^{\epsilon+\epsilon_1}}  \mathrm{d} s\\
			\lesssim &\, \left\|\underline{Z}\right\|_{\alpha,\alpha^{\prime},T}^{p_0} {t_0}^{-(\alpha^{\prime}+ \kappa)} \int_{0}^{T} s^{-\frac{\alpha+ \epsilon+\epsilon_1}{2}} \mathrm{d} s  
			\lesssim \left\|\underline{Z}\right\|_{\alpha,\alpha^{\prime},T}^{p_0} {t_0}^{-(\alpha^{\prime}+ \kappa)}   ,
		\end{align*}
		where we use the fact that $0<\alpha+\epsilon+\epsilon_1<2$. For $i+j=2$, using Proposition \ref{regularity embedding}, \ref{heat kernel smoothing}, \ref{multiplicative structure}, \eqref{norm of Z} and assumption \eqref{eq: bootstrap},
		\begin{align*}
			\left\| P_{t-s} v(s)^i\overline{v(s)}^jZ^{:2-i,1-j:}(s) \right\| _{\mathcal{B}_{2p,\infty}^{\epsilon+\epsilon_1}} 
			&\lesssim 	\left\| P_{t-s} v(s)^i\overline{v(s)}^jZ^{:2-i,1-j:}(s) \right\| _{\mathcal{B}_{p,\infty}^{\epsilon+\epsilon_1+\frac{1}{p}}} \\
			&\lesssim (t-s)^{-(\frac{\alpha+ \epsilon+\epsilon_1}{2}+\frac{1}{2p})}\left\| v(s)^i\overline{v(s)}^jZ^{:2-i,1-j:}(s) \right\| _{\mathcal{B}_{p,\infty}^{-\alpha}}
			\\
			&\lesssim (t-s)^{-(\frac{\alpha+ \epsilon+\epsilon_1}{2}+\frac{1}{2p})}\left\| v(s)^i\overline{v(s)}^j \right\| _{\mathcal{B}_{p,\infty}^{\epsilon}} \left\| Z^{:2-i,1-j:}(s) \right\| _{\mathcal{C}^{-\alpha}}\\
			&\lesssim (t-s)^{-(\frac{\alpha+ \epsilon+\epsilon_1}{2}+\frac{1}{2p})}  \left\| v(s) \right\| _{\mathcal{B}_{2p,\infty}^{\epsilon}}^2 \left\|\underline{Z}\right\|_{\alpha,\alpha^{\prime},T}^{p_0}\\
			&\lesssim (t-s)^{-(\frac{\alpha+ \epsilon+\epsilon_1}{2}+\frac{1}{2p})} s^{-2\kappa}  \left\|\underline{Z}\right\|_{\alpha,\alpha^{\prime},T}^{p_0}.
		\end{align*} 
		Then by Young's convolution inequality,
		\begin{align*}
			&\sup_{t_0\leq t\leq T} \int_{\frac{t_0}{2}}^{t} 	\left\| P_{t-s} v(s)^i\overline{v(s)}^jZ^{:2-i,1-j:}(s) \right\| _{\mathcal{B}_{2p,\infty}^{\epsilon}}  \mathrm{d} s\\
			\lesssim &\, \left\|\underline{Z}\right\|_{\alpha,\alpha^{\prime},T}^{p_0} {t_0}^{- 2\kappa} \int_{0}^{T} s^{-(\frac{\alpha+ \epsilon+\epsilon_1}{2}+\frac{1}{2p})} \mathrm{d} s  
			\lesssim \left\|\underline{Z}\right\|_{\alpha,\alpha^{\prime},T}^{p_0} {t_0}^{-2\kappa}   ,
		\end{align*}
		where we use the fact that $0<\frac{\alpha+ \epsilon+\epsilon_1}{2}+\frac{1}{2p}<1$ since $\alpha$ is sufficiently small and $\epsilon+\epsilon_1<\frac53-\frac{2}{3p} < 2-\frac{1}{p}$.
	\end{proof}

	\begin{proposition}\label{prop: B1+ est}
		Let $T>0$ and $v\in C\left( \left( 0,T\right]; \mathcal{C}^{\beta}\right) $ be a mild solution to \eqref{remainder} with $m=1$. Then there exist $\epsilon>1$, $p>1$, $p_0>0$ and $\kappa>0$, such that for any $0<t_0\leq T$,   
		\begin{equation*}
			\sup_{t_0\leq t\leq T}\left\| v(t)\right\| _{\mathcal{B}_{2p,\infty}^{\epsilon}} \lesssim \left\| \underline{Z}\right\|_{\alpha,\alpha^{\prime},T}^{p_0}  t_0^{-\kappa}.
		\end{equation*}
	\end{proposition}
	
	\begin{proof}
		By Proposition \ref{prop: B-epsilon-6p est} and Proposition \ref{prop: B-epsilon est}, we know that \eqref{eq: bootstrap} holds. Then using Proposition \ref{prop: B-2epsilon est} repeatedly and noting that $\frac53-\frac{2}{3p}>1$ for $p>1$, we get the conclusion. 
	\end{proof}

		\subsection{Proof of Theorem \ref{Global existence and uniqueness} and Theorem \ref{Global existence and uniqueness for m=1}: Global well-posedness}

		\begin{proof}[Proof of Theorem \ref{Global existence and uniqueness}]
		Let $T>0$. By Theorem \ref{Local existence and uniqueness}, there exists $0<T^{*}\leq T$ such that \eqref{remainder} has a unique mild solution on $[0,T^{*}]$. Under the assumptions $\mu>\frac{2m-1}{2\sqrt{2m}}$, which is equivalent to $2m+1<2+2\mu(\mu+\sqrt{1+\mu^2})$,  and $\alpha_0\in [1/( 1+\mu(\mu+\sqrt{1+\mu^2})) ,2/( 2m+1) ]$, by Proposition \ref{L^p embedding} and \ref{regularity embedding}, we can find $1\leq p\leq 1+\mu(\mu+\sqrt{1+\mu^2}) $ such that $L^{2p}\hookrightarrow\mathcal{C}^{-\alpha_0}$. Then by Proposition \ref{priori estimate}, there exists $p_0>1$ such that 
		\begin{align*}
			\sup_{T^{*} \leq t\leq T} \left\|v(t)\right\|_{\mathcal{C}^{-\alpha_0}} &\leq  \sup_{T^{*} \leq t\leq T} \left\|v(t)\right\|_{L^{2p}}\leq C \left\| \underline{Z}\right\|_{\alpha,\alpha^{\prime},T}^{p_0} \sup_{T^{*} \leq t\leq T}  t^{-\frac{1}{2m}}\\&\leq C \left\| \underline{Z}\right\|_{\alpha,\alpha^{\prime},T}^{p_0} {(T^{*})}^{-\frac{1}{2m}}.
		\end{align*}
		Using Theorem \ref{Local existence and uniqueness} again, there exists $T^{**}$ bounded from below and satisfying
		\begin{equation*}
			T^{**}\geq  C\left(\left\| \underline{Z}\right\|_{\alpha,\alpha^{\prime},T}^{p_0} {(T^{*})}^{-\frac{1}{2m}}\right)^{-\frac{2m+1}{\gamma-\frac{\alpha_0+\beta}{2}}}\wedge T,
		\end{equation*} 
		such that \eqref{remainder} has a unique mild solution on $[T^*,(T^*+T^{**})\wedge T]$ with initial condition $v(T^*)$. We then continue this process with the same time $T^{**}$ until the whole interval $[0, T]$ is covered.
	\end{proof}

		\begin{proof}[Proof of Theorem \ref{Global existence and uniqueness for m=1}]
			Using the embedding $\mathcal{B}_{2p, \infty}^{\epsilon} \hookrightarrow \mathcal{C}^{-\alpha_0} $ for any $p,\epsilon>1$, Proposition \ref{prop: B1+ est} and a similar argument in the proof of Theorem \ref{Global existence and uniqueness}, we get the conclusion.
		\end{proof}
		
		\section{Ergodicity}\label{Ergodicity}
		
	Let $u(\cdot;u_0)=Z(\cdot)+v(\cdot;u_0)$, where $Z$ as defined in \eqref{solution of heat equation} is the solution to \eqref{heat equation} with zero initial data at time $a=0$  and $v(\cdot;u_0)$ is the mild solution to \eqref{remainder} with initial data $u_0\in \mathcal{C}^{-\alpha_0}$. For every $t\geq0$ and $\Phi\in \mathcal{B}_b(\mathcal{C}^{-\alpha_0})$, where $\mathcal{B}_b(\mathcal{C}^{-\alpha_0})$ denotes the space of bounded measurable functions from $\mathcal{C}^{-\alpha_0}$ to $\mathbb{R}$, we define the map $S_t:\Phi\rightarrow S_t\Phi$ by 
	\begin{equation}\label{semigroup}
		S_t\Phi(x)=\mathrm{E}\left[ \Phi(u(t;x)\right] ,\quad x\in \mathcal{C}^{-\alpha_0}.
	\end{equation}
	Let $\left( S_t^{*}\right)_{t\geq0} $ be the dual semigroup of $\left( S_t\right)_{t\geq0} $, that is,
	\begin{equation*}
		S_t^{*}\eta(A)=\int_{\mathcal{C}^{-\alpha_0}}\mathrm{P}(u(t;x)\in A)\eta(\mathrm{d} x),
	\end{equation*}
	where $A$ is a Borel subset of $\mathcal{C}^{-\alpha_0}$ and $\eta$ is a probability measure on $\mathcal{C}^{-\alpha_0}$. We set 
	\begin{equation*}
		\mathcal{G}_t=\sigma\left( \left\lbrace \xi(\varphi): \varphi|_{(t,\infty)\times \mathbb{T}^2} \equiv0,\,\varphi\in L^2(\mathbb{R}\times\mathbb{T}^2)\right\rbrace \right) 
	\end{equation*}
	for $t>-\infty$ and let $\left( \mathcal{F}_t\right)_{t>-\infty} $ be the usual augmentation of the filtration $\left( \mathcal{G}_t\right)_{t>-\infty} $. By a similar argument in \cite[Section 4]{Matsuda2020}, we can show that $\left(u(t;u_0) \right)_{t\geq0} $ is a Feller Markov process  on $\mathcal{C}^{-\alpha_0}$ with transition semigroup $\left( S_t\right)  _{t\geq0}$ with respect to the filtration $\left( \mathcal{F}_t\right)_{t\geq0} $ by using Lemma \ref{lem: unique for u} and continuous dependence of $v(\cdot;u_0)$ on the initial data $u_0$. Then we prove that there exists a unique invariant measure of $\left( S_t\right)_{t\geq0}$ by the Krylov-Bogoliubov theorem and an asymptotic coupling argument.

		\subsection{Existence of invariant measure for $m>1$}

		In this section, we fix $m>1$, $\mu>\frac{2m-1}{2\sqrt{2m}}$ (which is equivalent to $2+2\mu(\mu+\sqrt{1+\mu^2})>2m+1$) and $\alpha_0,\beta$, $\gamma$ satisfying $$ \frac{1}{1+\mu(\mu+\sqrt{1+\mu^2})}\leq \alpha_0\leq \frac{2}{2m+1}, \quad \beta>0,\quad \frac{\alpha_0+\beta}{2}< \gamma<\frac{1}{2m+1},$$
		such that there exists a unique mild solution $v(t;u_0)\in C\left( \left( 0,\infty\right); \mathcal{C}^{\beta}\right)$ to \eqref{remainder} for every initial data $u_0\in \mathcal{C}^{-\alpha_0}$ by Theorem \ref{Global existence and uniqueness}. To derive a priori estimate of $u(\cdot;u_0)$, we first show the following lemma, which is proved in Appendix \ref{Some technical estimation}.
		\begin{lemma}\label{lem: unique for u}
			Let $u(t;u_0)=Z(t)+v(t;u_0)$ for $t\geq0$. Then for every $h\geq0$,
			\begin{equation*}
				u(t+h;u_0)= Z_{t}(t+h)+ v_t(t+h),
			\end{equation*}
			where $v_t(t+\cdot)$ solves \eqref{remainder} driven by $\underline{Z}_{t}(t+\cdot)$ with initial data $u(t;u_0)$, that is,
			\begin{equation*}
				v_t(t+h)= P_hu(t;u_0)+\int_{0}^{h} P_{h-s} \Psi\left(v_t(t+s),\underline{Z}_t(t+s)\right)\mathrm{d}s,\quad h\geq0.
			\end{equation*}
		\end{lemma}
		
		Based on the estimates for Wick products (Theorem \ref{regularity of Z}) and for the remainder $v$ (Proposition \ref{priori estimate}), we give a control of $u(t;u_0)$, which is useful for proving the existence of invariant probability measures via the Krylov-Bogoliubov theorem.
		
		\begin{corollary}\label{priori extimate for u}
			Let $\mu>\frac{2m-1}{2\sqrt{2m}}$, $1\leq p\leq 1+\mu(\mu+\sqrt{1+\mu^2})$ and $\alpha_0,\alpha_1\in[1/( 1+\mu(\mu+\sqrt{1+\mu^2})) ,2/\left( 2m+1\right) ]$. Then
			\begin{equation*}
				\sup_{u_0\in \mathcal{C}^{-\alpha_0}}\sup_{t\geq0}\left( t^{\frac{p}{m}}\wedge1 \right)\mathrm{E}\left[ \left\| u(t;u_0) \right\| _{\mathcal{C}^{-\alpha_1}}^{2p}\right]<\infty.   
			\end{equation*}
		\end{corollary}
		\begin{proof}
			Under the assumptions that $$2m+1<2+2\mu(\mu+\sqrt{1+\mu^2}),\quad \alpha_1\in \left[1/\left( 1+\mu(\mu+\sqrt{1+\mu^2})\right) ,2/\left( 2m+1\right) \right],$$ by Proposition \ref{L^p embedding} and \ref{regularity embedding}, we can find $1\leq p\leq 1+\mu(\mu+\sqrt{1+\mu^2}) $ such that $L^{2p}\hookrightarrow\mathcal{C}^{-\alpha_1}$. For $t>1$, by Lemma \ref{lem: unique for u}, we have that
			\begin{equation*}
				u(t;u_0)= Z_{t-1}(t)+ v_{t-1}(t),
			\end{equation*}
			where $v_{t-1}(t-1+\cdot)$ solves \eqref{remainder} driven by $\underline{Z}_{t-1}(t-1+\cdot)$ with initial data $u(t-1;u_0)$, that is,
			\begin{equation*}
				v_{t-1}(t)= P_1u(t-1;u_0)+\int_{0}^{1} P_{1-s} \Psi\left(	v_{t-1}(t-1+s),\underline{Z}_{t-1}(t-1+s)\right)\mathrm{d}s.
			\end{equation*}
			Using \eqref{variance of wick product}, we know that the law of Gaussian variable $Z_{t-1}(t-1+\cdot)$ is independent of $t$ and thus the law of $\underline{Z}_{t-1}(t-1+\cdot)$ is independent of $t$. Combining this fact with Theorem \ref{regularity of Z} and Proposition \ref{priori estimate}, we derive that there exists $p_0>0$ such that
			\begin{align*}
				\sup_{u_0\in \mathcal{C}^{-\alpha_0}}\sup_{t>1}\mathrm{E}\left[ \left\| u(t;u_0) \right\| _{\mathcal{C}^{-\alpha_1}}^{2p}\right] &\lesssim \sup_{u_0\in \mathcal{C}^{-\alpha_0}}\sup_{t>1}\mathrm{E}\left[ \left\|v_{t-1}(t) \right\| _{L^{2p}}^{2p}\right]  + \sup_{t>1}\mathrm{E}\left[ \left\| Z_{t-1}(t) \right\| _{\mathcal{C}^{-\alpha_1}}^{2p}\right]\\
				& \lesssim \sup_{t>1}\mathrm{E}\left[ \left\| \underline{Z}_{t-1}(t-1+\cdot)\right\|_{\alpha,\alpha^{\prime},1}^{p_0}\right]  +  \sup_{t>1}\mathrm{E}\left[ \left\| Z_{t-1}(t) \right\| _{\mathcal{C}^{-\alpha_1}}^{2p}\right]\lesssim1.
			\end{align*}
			For $0<t\leq 1$, by Theorem \ref{regularity of Z} and Proposition \ref{priori estimate}, there exists $p_0>0$ such that
			\begin{align*}
				\sup_{u_0\in \mathcal{C}^{-\alpha_0}}\sup_{0<t \leq 1}  t^{\frac{p}{m}}\mathrm{E}\left[ \left\| u(t;u_0) \right\| _{\mathcal{C}^{-\alpha_1}}^{2p}\right]
				& \lesssim 	 	\sup_{u_0\in \mathcal{C}^{-\alpha_0}}\sup_{0<t \leq 1}  t^{\frac{p}{m}}\mathrm{E}\left[ \left\|v(t;u_0) \right\| _{L^{2p}}^{2p}\right]  +  \sup_{0<t \leq 1} \mathrm{E}\left[ \left\| Z(t) \right\| _{\mathcal{C}^{-\alpha_1}}^{2p}\right]\\
				&\lesssim   	\mathrm{E}\left[\left\| \underline{Z}\right\|_{\alpha,\alpha^{\prime},1}^{p_0}\right]   +  \sup_{0<t \leq 1}  \mathrm{E}\left[ \left\| Z(t) \right\| _{\mathcal{C}^{-\alpha_1}}^{2p}\right]
				\lesssim   1.
			\end{align*}
		\end{proof}
		
		\begin{proposition}\label{Prop: existence for m>1}
			 Let $m>1$ and $\mu>\frac{2m-1}{2\sqrt{2m}}$. Then there exists a probability measure $\eta$ on $ \mathcal{C}^{-\alpha_0}$ such that $ S_t^{*}\eta=\eta$ for all $t\geq0$.
		\end{proposition}
		
		\begin{proof}
			Let $t>0$, $\alpha_1\in \left[1/\left( 1+\mu(\mu+\sqrt{1+\mu^2})\right) ,2/\left( 2m+1\right) \right]$ satisfy $\alpha_1<\alpha_0$ and $p\in[1,1+\mu(\mu+\sqrt{1+\mu^2})]$. Let $x\in \mathcal{C}^{-\alpha_0}$. By Markov's and Jensen's inequality, for every $a>0$,
			\begin{equation*}
				\mathrm{P}\left(\|u(t ; x)\|_{\mathcal{C}^{-\alpha_1}}>a\right) \leq \frac{1}{a}\mathrm{E}\left[ \|u(t ; x)\|_{\mathcal{C}^{-\alpha_1}}\right]  \leq \frac{1}{a}\left(\mathrm{E}\left[ \|u(t ; x)\|_{\mathcal{C}^{-\alpha_1}}^{2p}\right] \right)^{\frac{1}{2p}}.
			\end{equation*}
			Therefore, by Corollary \ref{priori extimate for u}, there exists a constant $C>0$ such that
			\begin{align}
				\int_0^t 	\mathrm{P}\left(\|u(s ; x)\|_{\mathcal{C}^{-\alpha_1}}>a\right) \mathrm{d} s & \leq \frac{1}{a} \int_0^t\left(\mathrm{E}\left[ \|u(s ; x)\|_{\mathcal{C}^{-\alpha_1}}^{2p}\right] \right)^{\frac{1}{2p}} \mathrm{d} s \label{tightness} \\
				& \leq \frac{C}{a}\left[\int_0^1 s^{-\frac{1}{2m}} \mathrm{d} s+\int_1^t 1 \mathrm{d} s\right] =\frac{C}{a} t.\nonumber
			\end{align}
			Let $Q_t(\cdot)=\frac{1}{t} \int_0^t S_s^* \mathbf{1}_{\{x\}}(\cdot) \mathrm{d} s$. For any $\epsilon>0$, we consider the set $K_{\varepsilon}:=\left\{y \in \mathcal{C}^{-\alpha_1}:\|y\|_{\mathcal{C}^{-\alpha_1}} \right. $ $\left. \leq  C/\varepsilon\right\}$. By \eqref{tightness}, we have that for every $t>0$,
			\begin{align*}
				Q_t\left(K_{\varepsilon} \right)&=1- Q_t\left(\left\{y \in \mathcal{C}^{-\alpha_1}:\|y\|_{\mathcal{C}^{-\alpha_1}}>\frac{C}{\varepsilon}\right\}\right)\\&=1-	\frac{1}{t}\int_0^t 	\mathrm{P}\left(\|u(s ; x)\|_{\mathcal{C}^{-\alpha_1}}>\frac{C}{\varepsilon}\right) \mathrm{d} s \geq 1-\varepsilon .
			\end{align*}
			Since the embedding $\mathcal{C}^{-\alpha_1} \hookrightarrow \mathcal{C}^{-\alpha_0}$ is compact for $\alpha_1<\alpha_0$, $K_{\varepsilon}$ is a compact subset of $\mathcal{C}^{-\alpha_0}$. This implies that $\left\{Q_t\right\}_{t >0}$ is tight in $\mathcal{C}^{-\alpha_0}$. Then by the Krylov-Bogoliubov theorem (see \cite[Corollary 3.1.2]{DZ1996}), there exists a probability measure $\eta$ on $ \mathcal{C}^{-\alpha_0}$ such that $ S_t^{*}\eta=\eta$ for all $t\geq0$.
		\end{proof}
		
		\subsection{Existence of invariant measure for $m=1$}
		
		In this section, we fix $m=1$, $\mu>0$ and $\alpha_0>0$, $\beta>0$, $\gamma$ satisfying $ \frac{\alpha_0+\beta}{2}< \gamma<\frac{1}{2m+1}=\frac13$,
		such that there exists a unique mild solution $v(t;u_0)\in C\left( \left( 0,\infty\right); \mathcal{C}^{\beta}\right)$ to \eqref{remainder} for every initial data $u_0\in \mathcal{C}^{-\alpha_0}$ by Theorem \ref{Global existence and uniqueness for m=1}. Before deriving a compact estimate for $u(\cdot;u_0)=Z(\cdot)+v(\cdot;u_0)$, we estimate $	\left\| v(t;u_0)\right\| _{\mathcal{C}^{-\alpha_1}}$ for any $0<t\leq 1$ and any $0<\alpha_1<1$.
		
		\begin{proposition}\label{prop: C-alpha est on 01}
			Let $v\in C\left( \left( 0,1\right]; \mathcal{C}^{\beta}\right) $ be a mild solution to \eqref{remainder} with $m=1$. For any $0<t\leq 1$ and any $0<\alpha_1<1$, there exist some exponents $p_1,q_1>0$ and $0<\kappa_1<1$ such that  
			\begin{align*}
				\left\| v(t)\right\| _{\mathcal{C}^{-\alpha_1}} \lesssim \left\| \underline{Z}\right\|_{\alpha,\alpha^{\prime},1}^{p_1} (1+ \left\|u_0 \right\|_{\mathcal{C}^{-\alpha_0}} )^{q_1} t^{-\kappa_1}.
			\end{align*}
		\end{proposition}
		
		\begin{proof}
			Let $T^*$ be the time obtained by the fixed point argument in \eqref{eq: cond for R and local T}, namely			
			\begin{equation*}
				T^*= C (1+ \left\|u_0 \right\|_{\mathcal{C}^{-\alpha_0}}+ \left\| \underline{Z}\right\|_{\alpha,\alpha^{\prime},1} )^{-\frac{3}{\gamma-\frac{\alpha_0+\beta}{2}}} \wedge 1.
			\end{equation*}
			Then for $0<t<T^{*}$,
			\begin{align*}
				\left\| v(t)\right\| _{\mathcal{C}^{-\alpha_1}} \lesssim 	\left\| v(t)\right\| _{\mathcal{C}^{\beta}} \lesssim t^{-\gamma} (1+ \left\|u_0 \right\|_{\mathcal{C}^{-\alpha_0}}+ \left\| \underline{Z}\right\|_{\alpha,\alpha^{\prime},1}).
			\end{align*}
			For $T^{*} \leq t\leq 1$, by Proposition \ref{regularity embedding} and \ref{embedding}, and Proposition \ref{prop: B1+ est}, there exist $p>1$, $\epsilon>1$ and $\kappa,p_0>0$ such that
			\begin{align*}
				\left\| v(t)\right\| _{\mathcal{C}^{-\alpha_1}}  
				&\lesssim \left\| v(t)\right\| _{\mathcal{B}_{2p,\infty}^{-\alpha_1+\frac1p}}  \lesssim \left\| v(t)\right\| _{\mathcal{B}_{2p,\infty}^{\epsilon}}\lesssim \left\| \underline{Z}\right\|_{\alpha,\alpha^{\prime},1}^{p_0}  t^{-\kappa} 
				 \lesssim \left\| \underline{Z}\right\|_{\alpha,\alpha^{\prime},1}^{p_0}  (T^{*})^{-\kappa}\\ 
				&\lesssim \left\| \underline{Z}\right\|_{\alpha,\alpha^{\prime},1}^{p_0} (1+ \left\|u_0 \right\|_{\mathcal{C}^{-\alpha_0}}+ \left\| \underline{Z}\right\|_{\alpha,\alpha^{\prime},1} )^{\kappa\frac{3}{\gamma-\frac{\alpha_0+\beta}{2}}}\\
				&\lesssim \left\| \underline{Z}\right\|_{\alpha,\alpha^{\prime},1}^{p_0+\kappa\frac{3}{\gamma-\frac{\alpha_0+\beta}{2}}} (1+ \left\|u_0 \right\|_{\mathcal{C}^{-\alpha_0}} )^{\kappa\frac{3}{\gamma-\frac{\alpha_0+\beta}{2}}}.
			\end{align*}
		\end{proof}

	Combining Proposition \ref{prop: C-alpha est on 01} with Proposition \ref{prop: B1+ est}, we give an estimate on $u(t;u_0)$ when $m=1$.
		
		\begin{corollary}\label{priori extimate for u with m=1}
			Let $m=1$ and $0<\alpha_1<1$. Then there exist $0<\kappa_1<1$ and $q_1>0$ such that 
			\begin{equation*}
				\sup_{t\geq 0}\left( t^{\kappa_1}\wedge1 \right)\mathrm{E}\left[ \left\| u(t;u_0) \right\| _{\mathcal{C}^{-\alpha_1}}\right]\lesssim (1+ \left\|u_0 \right\|_{\mathcal{C}^{-\alpha_0}} )^{q_1}.   
			\end{equation*}
		\end{corollary}
		\begin{proof}
			Following the similar argument in the proof of Corollary \ref{priori extimate for u}, using the embedding $\mathcal{B}_{2p,\infty}^{\epsilon}\hookrightarrow\mathcal{C}^{-\alpha}$ for any $p>1$, $\epsilon>1$, and applying Proposition \ref{prop: B1+ est} for $t>1$, Proposition \ref{prop: C-alpha est on 01} for $0<t\leq 1$, we finish the proof.
		\end{proof}

		\begin{proposition}\label{existence of IM for m=1}
			Let $m=1$ and $\mu>0$. Then there exists a probability measure $\eta$ on $ \mathcal{C}^{-\alpha_0}$ such that $ S_t^{*}\eta=\eta$ for all $t\geq0$.
		\end{proposition}
		
		\begin{proof}
			By a similar argument in the proof of Proposition \ref{Prop: existence for m>1}, we get the conclusion using Corollary \ref{priori extimate for u with m=1}.
		\end{proof}

		\subsection{A priori estimates of an auxiliary system}
		
		We first introduce some notations and recall an abstract result in \cite{HairerMattinglyScheutzow2011}. Let $\mathcal{P}$ be a Markov transition kernel on a Polish space $(\mathbb{X}, \rho)$ and let $\mathbb{X}_{\infty}=\mathbb{X}^{\mathbb{N}}$ be the space of one-sided infinite sequences with product topology. Denote by $\mathcal{M}(\mathbb{X})$, $\mathcal{M}\left(\mathbb{X}_{\infty}\right)$ and $\mathcal{M}\left(\mathbb{X}_{\infty} \times \mathbb{X}_{\infty}\right)$ the space of all Borel probability measures on $\mathbb{X}$, $\mathbb{X}_{\infty}$ and $\mathbb{X}_{\infty} \times \mathbb{X}_{\infty}$ respectively. Let $\mathcal{P}_{\infty}: \mathbb{X} \rightarrow \mathcal{M}\left(\mathbb{X}_{\infty}\right)$ be the probability kernel defined by stepping with the Markov kernel $\mathcal{P}$. For $\eta \in \mathcal{M}(\mathbb{X})$, we define $\eta \mathcal{P}_{\infty}\in \mathcal{M}\left(\mathbb{X}_{\infty}\right)$ as $\eta \mathcal{P}_{\infty}= \int_{\mathbb{X}} \mathcal{P}_{\infty}(x, \cdot) \mathrm{d} \eta(x)$.  Given $\eta_1, \eta_2 \in \mathcal{M}(\mathbb{X})$, we define the set of all generalized couplings as
		$$
		\tilde{\mathcal{C}}\left(\eta_1 \mathcal{P}_{\infty}, \eta_2 \mathcal{P}_{\infty}\right):=\left\{\Gamma \in \mathcal{M}\left(\mathbb{X}_{\infty} \times \mathbb{X}_{\infty}\right): \Gamma \circ \Pi_i^{-1} \ll \eta_i \mathcal{P}_{\infty}, i=1,2\right\},
		$$
		where $\Pi_i$ is the projection onto the $i$-th coordinate for $i=1,2$ and for two measures $M_1$ and $M_2$, $M_1\ll M_2$ means $M_1$ is absolutely continuous with respect to $M_2$. We define the diagonal at infinity as
		\begin{equation}\label{diagonal at infinity}
			D:=\left\{(x, y) \in \mathbb{X}_{\infty} \times \mathbb{X}_{\infty}: \lim _{n \rightarrow \infty} \rho\left(x_n, y_n\right)=0\right\} .
		\end{equation}
		We recall an abstract result based on asymptotic argument in \cite{HairerMattinglyScheutzow2011}.
		\begin{theorem}{\cite[Corollary 2.2]{HairerMattinglyScheutzow2011}}\label{HM Corollary 2.2}
			If there exists a Borel measurable set $B \subset \mathbb{X}$ such that
			\begin{enumerate}
				\item $\eta(B)>0$ for any invariant probability measure $\eta$ of $\mathcal{P}$;
				\item there exists a measurable map $B \times B \ni(x, y) \mapsto \Gamma_{x, y} \in \mathcal{M}\left(\mathbb{X}_{\infty} \times \mathbb{X}_{\infty}\right)$ such that $\Gamma_{x, y} \in \tilde{\mathcal{C}}\left(\delta_x \mathcal{P}_{\infty}, \delta_y \mathcal{P}_{\infty}\right)$ and $\Gamma_{x, y}(D)>0$  for every $x, y \in B$.
			\end{enumerate}
			Then there exists at most one invariant probability measure for $\mathcal{P}$.
		\end{theorem}
		
		In the following, we prove the uniqueness of invariant measure of $\left( S_t\right) _{t\geq 0}$ defined in \eqref{semigroup} using Theorem \ref{HM Corollary 2.2}. In order to construct an asymptotic coupling, we first consider the following auxiliary system
		\begin{equation}\label{auxiliary system}
			\begin{cases}
				\partial_t \tilde{v}=\left[ (\i+\mu) \Delta-1\right]  \tilde{v}+\lambda\left(v-\tilde{v} \right) +	\Psi(\tilde{v},\underline{Z}) , & t>0, x \in \mathbb{T}^2, \\
				\tilde{v}(0, \cdot)=u_1\in \mathcal{C}^{-\alpha_0},
			\end{cases}
		\end{equation}
		where $\lambda>1$ will be taken sufficiently large later. Compared to \eqref{remainder}, there is a new dissipation term $\lambda\left(v-\tilde{v} \right)$ in \eqref{auxiliary system} and the initial data $u_1\in \mathcal{C}^{-\alpha_0}$ is different. By a similar argument in the proofs of Theorem \ref{Global existence and uniqueness} and Theorem \ref{Global existence and uniqueness for m=1}, we can show that \eqref{auxiliary system} is globally well-posed.
		\begin{theorem}\label{Global existence and uniqueness for auxiliary system}
		Let $m>1$, $\mu>\frac{2m-1}{2\sqrt{2m}}$,  $\frac{1}{1+\mu(\mu+\sqrt{1+\mu^2})}\leq \alpha_0\leq \frac{2}{2m+1}$,
		or 
		$m=1$, $\mu>0$, $\alpha_0>0$. Let $\beta>0$ and $  \frac{\alpha_0+\beta}{2}< \gamma<\frac{1}{2m+1}$. 
		Then for every $u_1\in \mathcal{C}^{-\alpha_0}$, there exists a unique mild solution $\tilde{v}\in C\left( \left( 0,\infty\right); \mathcal{C}^{\beta}\right)$ of \eqref{auxiliary system}.  
	\end{theorem}
	Analogue to Proposition \ref{priori estimate}, we can prove the following proposition.

	\begin{proposition}\label{priori estimate for auxiliary system}
		Let $T\geq1 $ and $\tilde{v}\in C\left( \left( 0,T\right]; \mathcal{C}^{\beta}\right) $ be a mild solution to \eqref{auxiliary system}. Let $1\leq p\leq 1+\mu(\mu+\sqrt{1+\mu^2})$ and $0\leq \delta<1$ satisfy \eqref{delta condition}. Then there exist $C>0$ and $p_0>0$, such that for every $1\leq t\leq T$,  
		\begin{align}
			&\frac{1}{2 p}\|\tilde{v}(t)\|_{L^{2 p}}^{2 p}+\frac{\delta \mu}{2} \int_{1}^t\left\||\nabla \tilde{v}(s)|^2|\tilde{v}(s)|^{2 p-2}\right\|_{L^1}\mathrm{d} s +\frac{\mathrm{Re}\nu}{2} \int_{1}^t\|\tilde{v}(s)\|_{L^{2 p+2m}}^{2 p+2m}\mathrm{d} s \label{est for tildev}\\ \leq&\,    C    \left\| \underline{Z}\right\|_{\alpha,\alpha^{\prime},1}^{p_0}+ C  \sum_{i=0}^{m+1}\sum_{j=0}^{m} \int_{1}^t \left\| Z^{:m+1-i,m-j:}(s) \right\| _{\mathcal{C}^{-\alpha}}^{p_0} \mathrm{d} s .\nonumber
		\end{align}
	\end{proposition}
	\begin{proof}
		Similar to Proposition \ref{L^p of v(t)} and \ref{expression of v_2p}, we have that for every $0<s\leq t\leq T$,
		\begin{align}
			& \frac{1}{2 p}\partial_s\|\tilde{v}(s)\|_{L^{2 p}}^{2 p}+\delta \mu \left\||\nabla \tilde{v}(s)|^2|\tilde{v}(s)|^{2 p-2}\right\|_{L^1}+ \mathrm{Re}\nu \|\tilde{v}(s)\|_{L^{2 p+2m}}^{2 p+2m}+(\lambda+1)\|\tilde{v}(s)\|_{L^{2 p}}^{2 p} \label{experssion of tilde v 2p}\\
			\leq&\,\left\langle |\tilde{v}(s)|^{2 p-2},\mathrm{Re}\left(\overline{\tilde{v}(s)}\Psi^{\prime}\left( \tilde{v}(s),\underline{Z}(s)\right)\right) \right\rangle+\lambda \left\langle |\tilde{v}(s)|^{2 p-2},\mathrm{Re}\left(\overline{\tilde{v}(s)}v(s)\right) \right\rangle.\nonumber
		\end{align}
		Using the similar argument to \eqref{vRev}, we know that there exist $\tilde{q}\geq1$ and $C>0$ such that 
		\begin{align}
			& \left| 	\left\langle |\tilde{v}(s)|^{2 p-2},\mathrm{Re}\left(\overline{\tilde{v}(s)}\Psi^{\prime}\left( \tilde{v}(s),\underline{Z}(s)\right)\right) \right\rangle \right| \label{Re tilde v}\\
			\leq&\,  \frac{\delta\mu}{2}  \left\||\nabla \tilde{v}(s)|^2|\tilde{v}(s)|^{2 p-2}\right\|_{L^1}+ \frac{\mathrm{Re}\nu}{2} \|\tilde{v}(s)\|_{L^{2 p+2m}}^{2 p+2m}+ C  \sum_{i=0}^{m+1}\sum_{j=0}^{m} \left\| Z^{:m+1-i,m-j:}(s) \right\| _{\mathcal{C}^{-\alpha}}^{\tilde{q}}.\nonumber
		\end{align}
		For the term $ \left\langle |\tilde{v}(s)|^{2 p-2},\mathrm{Re}\left(\overline{\tilde{v}(s)}v(s)\right) \right\rangle  $, by H\"older's and Young's inequality,
		\begin{align}
			\left| \left\langle |\tilde{v}(s)|^{2 p-2},\mathrm{Re}\left(\overline{\tilde{v}(s)}v(s)\right) \right\rangle \right|&\leq \left\langle |\tilde{v}(s)|^{2p-1},|v(s)|\right\rangle
			\leq \left\| v(s)\right\| _{L^{2p}}  \left\| \tilde{v}(s)\right\| _{L^{2p}}^{2p-1} \label{estimate of new term}\\
			&\leq \frac{1}{2p} \left\| v(s)\right\| _{L^{2p}} ^{2p}+\frac{2p-1}{2p} \left\| \tilde{v}(s)\right\| _{L^{2p}}^{2p} \nonumber\\
			&\leq \frac{1}{2p} \left\| v(s)\right\| _{L^{2p}} ^{2p}+ \left\| \tilde{v}(s)\right\| _{L^{2p}}^{2p}.\nonumber
		\end{align}
		Combining \eqref{experssion of tilde v 2p}, \eqref{Re tilde v} with \eqref{estimate of new term}, we get that for $1\leq s\leq t$, 
		\begin{align}
			&\frac{1}{2 p}\partial_s\|\tilde{v}(s)\|_{L^{2 p}}^{2 p}+\frac{\delta \mu}{2} \left\||\nabla \tilde{v}(s)|^2|\tilde{v}(s)|^{2 p-2}\right\|_{L^1}+ \frac{\mathrm{Re}\nu}{2} \|\tilde{v}(s)\|_{L^{2 p+2m}}^{2 p+2m}+\|\tilde{v}(s)\|_{L^{2 p}}^{2 p} \label{differentiate tilde v}\\
			\leq&\,  \frac{\lambda}{2p} \left\| v(s)\right\| _{L^{2p}} ^{2p}+ C  \sum_{i=0}^{m+1}\sum_{j=0}^{m} \left\| Z^{:m+1-i,m-j:}(s) \right\| _{\mathcal{C}^{-\alpha}}^{\tilde{q}}.\nonumber
		\end{align}
		By a similar argument to Proposition \ref{priori estimate}, we derive that there exists $p^{\prime}\geq 0$ such that
		\begin{equation}\label{tilde v1}
			\|\tilde{v}(1)\|_{L^{2 p}}^{2 p}\leq C\left\| \underline{Z}\right\|_{\alpha,\alpha^{\prime},1}^{p^{\prime}}.
		\end{equation}
		Integrating \eqref{differentiate tilde v} with respect $s$ from $1$ to $t$ and using \eqref{tilde v1} and Proposition \ref{priori estimate}, we have that there exist $p_0\geq 1$ and $C>0$ depending on $\lambda$ such that \eqref{est for tildev} holds.
	\end{proof}
		
		\subsection{Proof of Theorem \ref{ergodicity of u}: Ergodicity}\label{sec: Proof of Ergodicity}
		
		Using an asymptotic coupling argument (see \cite[Corollary 2.2]{HairerMattinglyScheutzow2011} or Theorem \ref{HM Corollary 2.2} below), we prove the uniqueness of invariant measure of the transition semigroup $\left( S_t\right) _{t\geq 0}$ in this section.

		\begin{proof}[Proof of Theorem \ref{ergodicity of u}]
			Let $u(t;u_0)=Z(t)+v\left(t;u_0\right)$ and $\tilde{u}(t;u_1)=Z(t)+\tilde{v}\left(t;u_1\right)$, where $v\left(t;u_0\right)$, $\tilde{v}\left(t;u_1\right)$ is the unique mild solution to \eqref{remainder} and \eqref{auxiliary system} with initial data $u_0,u_1 \in \mathcal{C}^{-\alpha_0}$ respectively. We apply Theorem \ref{HM Corollary 2.2} to $B=\mathbb{X}=\mathcal{C}^{-\alpha_0}$ and $\mathcal{P}=P_1(x, \mathrm{d} y)$, where $P_1(x, \cdot)$ denotes the marginal distribution of $u(t;x)$ with initial data $x\in \mathcal{C}^{-\alpha_0}$ at time $t=1$.

			\textbf{Step 1.} We first construct a coupling $\Gamma_{u_0, u_1} \in \tilde{\mathcal{C}}\left(\delta_{u_0} \mathcal{P}_{\infty}, \delta_{u_1} \mathcal{P}_{\infty}\right)$ for any $u_0,u_1 \in \mathcal{C}^{-\alpha_0}$. 
			
			Let $w(t):=\tilde{u}\left(t;u_1\right)-u\left(t;u_0\right)=\tilde{v}\left(t;u_1\right)-v\left(t;u_0\right)$ and $\tilde{W}(t)=W(t)-\int_0^{t \wedge \tau_R} \lambda w(s) \mathrm{d} s$, where $W$ is a cylindrical Wiener process defined as \eqref{xi(t,x)} and
			$$
			\tau_R:=\inf \left\{t>0: \int_0^t\left\|w(s)\right\|_{L^2}^2 \mathrm{d} s \geq R\right\} .
			$$
			Note that
			$$
			\mathrm{E} \left[\exp \left(\frac{1}{2} \int_0^{\tau_R}\|\lambda w(s)\|_{L^2}^2 \mathrm{d} s\right) \right] \leq e^{\frac{1}{2} R \lambda^2},
			$$
			then by Girsanov theorem (see \cite{CM1944,Girsanov1960} or \cite[Theorem 5.1]{Karatzas1991}), there exists a probability measure $Q$ on $\left(\Omega, \mathcal{F},\left(\mathcal{F}_t\right)_{t \geq 0}\right)$ such that under $Q$, $\tilde{W}$ is a cylindrical Wiener process. Moreover, it holds that $\mathrm{P}\sim Q$ on  $\mathcal{F}_{\infty}=\sigma\left(\cup_{t \geq 0} \mathcal{F}_t\right)$.

			Let $\tilde{Z}$ be the solution to the following linear equation
			\begin{equation*}
				\begin{cases}
					\mathrm{d} \tilde{Z}(t)=\left[ (\i+\mu) \Delta-1\right]  \tilde{Z}(t) \mathrm{d} t+\mathrm{d} \tilde{W}(t), & t>0, x \in \mathbb{T}^2, \\
					\tilde{Z}(0)=0.
				\end{cases} 
			\end{equation*}
			Then $\tilde{Z}^{:k,l:}$ can be defined similarly as ${Z}^{:k,l:}$ in Section \ref{Stochastic heat equation with dispersion}. Since $\tilde{Z}(t)=Z(t)+a(t)$, where $a(t)=-\int_0^t e^{(t-s) A} \lambda w(s) \textbf{1}_{\left\lbrace s \leq \tau_R\right\rbrace}  \mathrm{d} s \in C([0, \infty) ; \mathcal{C}^{\alpha^{\prime}})$ for some $\alpha^{\prime}>\alpha$, similar to Corollary \ref{Cor relation of wick product of diff initial moments}, we have that $\mathrm{P}$-almost surely, 
			\begin{equation}\label{hat{Z}_{0,t}}
				\tilde{Z}^{:k,l:}=\sum_{i=0}^{k}\sum_{j=0}^{l}\binom{k}{i}\binom{l}{j} a  ^i \overline{a }^j Z^{:k-i,l-j:}.
			\end{equation}
			Using Theorem \ref{regularity of Z} and $a(t) \in C([0, \infty) ; \mathcal{C}^{\alpha^{\prime}})$ for some $\alpha^{\prime}>\alpha$, we get that for $\mathrm{P}$-almost surely $\omega\in\Omega$,
			\begin{align*}
				\underline{\tilde{Z}}:=\left\lbrace \tilde{Z}^{:i,j:}: (i,j)\in L\right\rbrace \in C([0,\infty);\mathcal{C}^{-\alpha})\times \left( C(( 0,\infty);\mathcal{C}^{-\alpha})\right) ^{M-1},
			\end{align*}
			and for every $T>0$, there exists a constant $C$ depending on $\alpha$, $\alpha^{\prime}$ and $T$ such that
			\begin{align*}
				\| \underline{\tilde{Z}}\|_{\alpha,\alpha^{\prime},T}:=1\vee \sup_{(i,j)\in L} \sup_{0< t\leq T } t^{(i+j-1)\alpha^{\prime}} \|  \tilde{Z}^{:i,j:}(t,\cdot)\| _{\mathcal{C}^{-\alpha}} <C,
			\end{align*} 
			where the set $L$ is defined in \eqref{set L}. 
			Then by Theorem \ref{Global existence and uniqueness}, for $\mathrm{P}$-almost surely $\omega\in\Omega$, there exists a unique mild solution $\hat{v}(\cdot;u_1) \in C\left((0, \infty), \mathcal{C}^\beta\right)$ to the following equation
			\begin{equation*}
				\begin{cases}
					\partial_t \hat{v}=\left[ (\i+\mu) \Delta-1\right]  \hat{v}+	\Psi(\hat{v},\underline{\tilde{Z}}) , & t>0, x \in \mathbb{T}^2, \\
					\hat{v}(0, \cdot)=u_1\in \mathcal{C}^{-\alpha_0}.
				\end{cases}
			\end{equation*}
			Combining Theorem \ref{Global existence and uniqueness} with Yamada-Watanabe Theorem (see \cite{Kurtz2007}), we have that the law of $\hat{v}(t;u_1)$ under $Q$ equals the law of $v(t;u_1)$ under $\mathrm{P}$. Let $\hat{u}(t;u_1)=\hat{v}(t;u_1)+\tilde{Z}(t)$. Then the law of $\hat{u}(t;u_1)$ under $Q$ equals the law of $u(t;u_1)$ under $\mathrm{P}$. Since $\mathrm{P} \sim Q$, the marginal distributions of the pair $\left(u\left(\cdot;u_0\right), \hat{u}\left(\cdot; u_1\right)\right)$ are equivalent to these of $\left(u\left(\cdot;u_0\right), u\left(\cdot; u_1\right)\right)$ under $\mathrm{P}$. For $u_0, u_1 \in \mathcal{C}^{-\alpha_0}$, we set 
			\begin{equation}\label{generalized coupling}
				\Gamma_{u_0, u_1}:= \mbox{law of }\left(u\left(\cdot;u_0\right), \hat{u}\left(\cdot; u_1\right)\right).
			\end{equation}
			Then $\Gamma_{u_0, u_1} \in \tilde{\mathcal{C}}\left(\delta_{u_0} \mathcal{P}_{\infty}, \delta_{u_1} \mathcal{P}_{\infty}\right)$. 
			
			\textbf{Step 2.} We next prove that $\Gamma_{u_0, u_1}(D)>0$.
			
			Since $\hat{u}(t;u_1)=\hat{v}(t;u_1)+\tilde{Z}(t)$, $\hat{u}\left(t;u_1\right)$ satisfies the following equation with initial data $u_1\in \mathcal{C}^{-\alpha_0}$ under $\mathrm{P}$,
			\begin{align*}
				\partial_t \hat{u}&=\left[ (\i+\mu) \Delta-1\right]  \hat{u}+	\Psi(\hat{u}-\tilde{Z},\underline{\tilde{Z}}) +\xi-\lambda w\mathbf{1}_{\left\lbrace t\leq \tau_R  \right\rbrace }\\
				&=\left[ (\i+\mu) \Delta-1\right]  \hat{u}+	\Psi(\hat{u}-Z,\underline{Z}) +\xi-\lambda w\mathbf{1}_{\left\lbrace t\leq \tau_R  \right\rbrace },\quad t>0, x \in \mathbb{T}^2.
			\end{align*}
			where we use \eqref{hat{Z}_{0,t}} to get that $\Psi(\hat{u}-\tilde{Z},\underline{\tilde{Z}})=\Psi(\hat{u}-Z,\underline{Z}) $.
			Then on $\left\{\tau_R=\infty\right\}$, $\hat{u}(t;u_1)-Z(t)$ also satisfies \eqref{auxiliary system}. By Theorem \ref{Global existence and uniqueness for auxiliary system}, we obtain that on $\left\{\tau_R=\infty\right\}$,  $\hat{u}(t;u_1)-Z(t)=\tilde{v}(t;u_1)=\tilde{u}(t;u_1)-Z(t)$, which implies that $\hat{u}(t;u_1)=\tilde{u}(t;u_1)$ on $\left\{\tau_R=\infty\right\}$. To prove $\Gamma_{u_0, u_1}(D)>0$, it is sufficient to estimate
			$w(t)=\tilde{u}\left(t;u_1\right)-u\left(t;u_0\right)=\tilde{v}\left(t;u_1\right)-v\left(t;u_0\right)$ on $\left\{\tau_R=\infty\right\}$.
			
			By the definition of $w(t)$, we know that $w$ is the mild solution to the following equation
			\begin{equation*}
				\begin{cases}
					\partial_t w=\left[ (\i+\mu)\Delta-1\right] w-\lambda w+\Psi\left(\tilde{v}, \underline{Z}\right)-\Psi\left(v, \underline{Z}\right) , & t>0, x \in \mathbb{T}^2, \\
					w(0,\cdot)=u_1-u_0.
				\end{cases}
			\end{equation*}
			Testing against $\overline{w}_t$ and using a similar argument in Proposition \ref{expression of v_2p}, we obtain that for every $0\leq \delta<1$ (see \eqref{delta condition} for $p=1$),
			$$
			\frac{1}{2} \partial_t\|w(t)\|_{L^2}^2+\delta\mu\|\nabla w(t)\|_{L^2}^2+(\lambda+1)\|w(t)\|_{L^2}^2 \leq \mathrm{Re}\left[ \left\langle \overline{w(t)}, \Psi\left(\tilde{v}(t), \underline{Z}(t)\right)-\Psi\left(v(t), \underline{Z}(t)\right)\right\rangle\right] ,$$
			where
			\begin{align*}
				&\Psi\left(\tilde{v}, \underline{Z}\right)-\Psi\left(v, \underline{Z}\right)\\=&\, -\nu\sum_{i=0}^{m+1} \sum_{j=0}^{m}\binom{m+1}{i} \binom{m}{j} \mathbf{1}_{\left\lbrace i+j>0 \right\rbrace }  Z^{:m+1-i,m-j:}\left(\tilde{v}^{i}\overline{\tilde{v}}^j- v^{i}\overline{v}^j\right)+\tau\left(\tilde{v}-v\right) \\
				=&\,-\nu\sum_{i+j>0} \binom{m+1}{i}\binom{m}{j} Z^{:m+1-i,m-j:}\left[  \left(\tilde{v}-v\right)\overline{v}^j\sum_{l=0}^{i-1} v^{i-1-l}\tilde{v}^{l}+  \left(\overline{\tilde{v}}-\overline{v}\right)\tilde{v}^i\sum_{s=0}^{j-1} \overline{v}^{j-1-s}\overline{\tilde{v}}^{s}\right]\\&\,+\tau\left(\tilde{v}-v\right).
			\end{align*}
			Then 
			\begin{align*}
				&\left| \mathrm{Re}\left[ \left\langle \overline{w}, \Psi\left(\tilde{v}, \underline{Z}\right)-\Psi\left(v, \underline{Z}\right)\right\rangle\right]\right| \leq 2 \left\langle \left| w\right|, \left| \Psi\left(\tilde{v}, \underline{Z}\right)-\Psi\left(v, \underline{Z}\right) \right|  \right\rangle 
				\\\leq&\, C  \sum_{i+j>0} \left\langle \left| w\right|^2 ,   \left| Z^{:m+1-i,m-j:}\right| \left(  \sum_{l=0}^{i-1}\left|v \right|^{i+j-1-l} \left| \tilde{v}\right| ^{l} +   \sum_{s=0}^{j-1}\left|\tilde{v} \right|^{i+s} \left|v\right| ^{j-1-s}\right)  \right\rangle +|\tau|\|w\|_{L^2}^2.
			\end{align*}
			We now estimate each term $ \left\langle \left| w\right|^2 ,    \left|v \right|^{k} \left| \tilde{v}\right| ^{i+j-1-k} \left| Z^{:m+1-i,m-j:}\right| \right\rangle$ for $(i,j,k)$ belongs to the set
			\begin{equation}\label{set I}
				I:=\left\lbrace (i,j,k)\in \mathbb{N}^3: 0\leq i\leq m+1, 0\leq j\leq m, (i,j)\neq(0,0), 0\leq k\leq i+j-1  \right\rbrace.
			\end{equation}
			By Proposition \ref{Sobolev embedding} and \ref{interpolation and multiplicative for Sobolev}, we have that for $0<\alpha<s<\frac12$,
			$$
			\left\||w|^2\right\|_{\mathcal{B}_{4/3, 1}^{\alpha}} \lesssim \left\||w|^2\right\|_{H_{4/3}^s}  \lesssim\left\| w\right\|_{H_2^s}\|w\|_{L^4}  \lesssim\left\| w\right\|_{H_2^{1/2}}^2 \lesssim\|w\|_{L^2}\left(\|\nabla w\|_{L^2}+\|w\|_{L^2}\right).
			$$
			Then by Proposition \ref{duality property} and Young's inequality,
			\begin{align*}
				&\left\langle \left| w\right|^2 ,    \left|v \right|^{k} \left| \tilde{v}\right| ^{i+j-1-k} \left| Z^{:m+1-i,m-j:}\right| \right\rangle
				\\\leq&\, C\left\||w|^2\right\|_{\mathcal{B}_{4/3, 1}^{\alpha}}\left\|v^{k}  \tilde{v} ^{i+j-1-k}  Z^{:m+1-i,m-j:} \right\|_{\mathcal{B}_{4, \infty}^{-\alpha}}
				\\\leq&\, C\|w\|_{L^2}\left(\|\nabla w\|_{L^2}+\|w\|_{L^2}\right)\left\|v ^{k}  \tilde{v} ^{i+j-1-k}  Z^{:m+1-i,m-j:} \right\|_{\mathcal{B}_{4, \infty}^{-\alpha}}
				\\\leq&\, C\|w\|_{L^2}^2\left(1+\left\|v ^{k}  \tilde{v} ^{i+j-1-k}  Z^{:m+1-i,m-j:} \right\|_{\mathcal{B}_{4, \infty}^{-\alpha}}^2\right)  +\frac{2\delta\mu}{(m+1)(m+2)(2m+1)}\|\nabla w\|_{L^2}^2,
			\end{align*}
			where $(m+1)(m+2)(2m+1)/2$ is the cardinality of the set $I$ (see \eqref{set I}). 
			Therefore, 
			\begin{align*}
				&\frac{1}{2} \partial_t\|w(t)\|_{L^2}^2+\lambda\|w(t)\|_{L^2}^2 \\\leq&\, C \|w(t)\|_{L^2}^2 \Bigg(1+ \sum_{(i,j,k)\in I} \left\|v(t) ^{k}  \tilde{v}(t) ^{i+j-1-k}  Z^{:m+1-i,m-j:}(t) \right\|_{\mathcal{B}_{4, \infty}^{-\alpha}}^2 \Bigg)\\ =&:\|w(t)\|_{L^2}^2 L(t).
			\end{align*}
			By Gronwall's inequality, for $t\geq1$,
			\begin{equation}\label{Gronwall}
				\|w(t)\|_{L^2}^2\leq \|w(1)\|_{L^2}^2 \exp\left[ \int_{1}^{t}2(-\lambda+L(s)) \mathrm{d}s\right] .
			\end{equation}

			Recall that by Lemma \ref{lem setE}, for every $\gamma>0$, there exists $K > 0$ such that $\mathrm{P}(E_{K,\gamma})>0$, where $E_{K,\gamma}$ is defined in \eqref{E_{K,gamma}}. We estimate each term of $\int_{1}^{t}L(s) \mathrm{d}s$ on $E_{K,\gamma}$ with $\gamma> 0$ to be determined later. For any $(i,j,k)\in I$ (see \eqref{set I} for the definition of the set $I$), let $l=i+j-1-k$.
			Note that for $1\leq s\leq t$, by Proposition \ref{embedding}, \ref{regularity embedding}, \ref{Sobolev embedding} and \ref{interpolation and multiplicative for Sobolev}, for any $\beta^{\prime}>\alpha$,
			\begin{align*}
				\left\|v(s)^{k}  \tilde{v}(s) ^{l}  \right\|_{\mathcal{B}_{4, \infty}^{\beta^{\prime}}}^2
				&\lesssim \left\|v(s)^{k} \tilde{v}(s) ^{l}  \right\|_{\mathcal{B}_{p, \infty}^{\beta^{\prime}+\frac2p-\frac12}}^2\lesssim \left\|v(s)^{k} \tilde{v}(s) ^{l} \right\|_{\mathcal{B}_{p, 1}^{\beta^{\prime}+\frac2p-\frac12}}^2\\
				&\lesssim \left\|v(s)^{k} \tilde{v}(s) ^{l}  \right\|_{H_{p}^{\beta_0}}^2\lesssim \left\|v(s)^{k} \tilde{v}(s) ^{l}\right\|_{L^p}^{2(1-\beta_0)}   \left\|v(s)^{k}\tilde{v}(s) ^{l} \right\|_{H^1_p}^{2\beta_0}\\
				&\lesssim \left\|v(s)^{k} \tilde{v}(s) ^{l} \right\|_{L^p}^{2(1-\beta_0)}  \left(  \left\|v(s)^{k} \tilde{v}(s) ^{l} \right\|_{L^p}^{2\beta_0}+ \left\| v(s)^{k} \tilde{v}(s) ^{l-1} \nabla \tilde{v}(s) \right\|_{L^p}^{2\beta_0} \right. \\ & \left. \quad+ \left\|\tilde{v}(s) ^{l}  v(s)^{k-1}\nabla v(s) \right\|_{L^p}^{2\beta_0} \right)\\
				&\lesssim \left\|v(s)^{k} \tilde{v}(s) ^{l} \right\|_{L^{p}}^{2}+  \left\|v(s)^{k} \tilde{v}(s) ^{l} \right\|_{L^p}^{2(1-\beta_0)} \left\| v(s)^{k} \tilde{v}(s) ^{l-1} \nabla \tilde{v}(s) \right\|_{L^p}^{2\beta_0} \\&\quad + \left\|v(s)^{k} \tilde{v}(s) ^{l} \right\|_{L^p}^{2(1-\beta_0)}  \left\|\tilde{v}(s) ^{l}  v(s)^{k-1}\nabla v(s) \right\|_{L^p}^{2\beta_0},
			\end{align*}
			where $p<2$ but close to $2$, $\beta_0=\beta^{\prime}+\frac2p-\frac12+\epsilon>\frac12$ and $\epsilon>0$ is arbitrary. 
			By H\"older's inequality, we get that
			\begin{align*}
				\left\|v(s)^{k} \tilde{v}(s) ^{l} \right\|_{L^p} &\leq 	\left\|v(s) \right\|_{L^{2kp}}^{k} \left\|\tilde{v}(s) \right\|_{L^{2lp}}^{l},\\
				\left\| v(s)^{k} \tilde{v}(s) ^{l-1} \nabla \tilde{v}(s) \right\|_{L^p} &\leq  \left\| v(s)\right\|_{L^{\frac{2kp}{2-p}}}^{k}  \left\|  \tilde{v}(s) ^{2(l-1)} |\nabla \tilde{v}(s) |^2 \right\|_{L^1}^{1/2},\\
				\left\|\tilde{v}(s) ^{l}  v(s)^{k-1}\nabla v(s) \right\|_{L^p}&\leq  \left\| \tilde{v}(s)\right\|_{L^{\frac{2lp}{2-p}}}^{l}  \left\|  v(s) ^{2(k-1)} |\nabla v(s) |^2 \right\|_{L^1}^{1/2}.
			\end{align*}
			Thus,
			\begin{align*}
				\left\|v(s)^{k}  \tilde{v}(s) ^{l}  \right\|_{\mathcal{B}_{4, \infty}^{\beta^{\prime}}}^2 &\lesssim \left\|v(s) \right\|_{L^{2kp}}^{2k} \left\|\tilde{v}(s) \right\|_{L^{2lp}}^{2l} \\&\quad + \left\|v(s) \right\|_{L^{2kp}}^{2k(1-\beta_0)} \left\|\tilde{v}(s) \right\|_{L^{2lp}}^{2l(1-\beta_0)}  \left\| v(s)\right\|_{L^{\frac{2kp}{2-p}}}^{2k\beta_0}  \left\|  \tilde{v}(s) ^{2(l-1)} |\nabla \tilde{v}(s) |^2 \right\|_{L^1}^{\beta_0} \\&\quad + \left\|v(s) \right\|_{L^{2kp}}^{2k(1-\beta_0)} \left\|\tilde{v}(s) \right\|_{L^{2lp}}^{2l(1-\beta_0)}  \left\| \tilde{v}(s)\right\|_{L^{\frac{2lp}{2-p}}}^{2l\beta_0}  \left\|  v(s) ^{2(k-1)} |\nabla v(s) |^2 \right\|_{L^1}^{\beta_0}\\
				&\lesssim \left\|v(s) \right\|_{L^{2kp}}^{2k} \left\|\tilde{v}(s) \right\|_{L^{2lp}}^{2l}  +  \left\|\tilde{v}(s) \right\|_{L^{2lp}}^{2l(1-\beta_0)}  \left\| v(s)\right\|_{L^{\frac{2kp}{2-p}}}^{2k}  \left\|  \tilde{v}(s) ^{2(l-1)} |\nabla \tilde{v}(s) |^2 \right\|_{L^1}^{\beta_0} \\&\quad+  \left\|v(s) \right\|_{L^{2kp}}^{2k(1-\beta_0)}  \left\| \tilde{v}(s)\right\|_{L^{\frac{2lp}{2-p}}}^{2l}  \left\|  v(s) ^{2(k-1)} |\nabla v(s) |^2 \right\|_{L^1}^{\beta_0},
			\end{align*}
			where we use the fact that $p>1$, which implies $p<\frac{p}{2-p}$. 
			Then for $(i,j,k)\in I$, we still let $l=i+j-1-k$ and by Proposition \ref{multiplicative structure},
			\begin{align*}
				&\int_{1}^{t}  \left\|v(s) ^{k}  \tilde{v}(s) ^{i+j-1-k}  Z^{:m+1-i,m-j:}(s) \right\|_{\mathcal{B}_{4, \infty}^{-\alpha}}^2 \mathrm{d}s 
				\\\lesssim &\,  \int_{1}^{t} \left\|v(s)^{k}  \tilde{v}(s) ^{l}  \right\|_{\mathcal{B}_{4, \infty}^{\beta^{\prime}}}^{2}  
				\left\| Z^{:m+1-i,m-j:}(s) \right\|_{\mathcal{C}^{-\alpha}}^2 \mathrm{d}s 
				\\\lesssim &\,  \int_{1}^{t}  \left( \left\|v(s) \right\|_{L^{2kp}}^{2k} \left\|\tilde{v}(s) \right\|_{L^{2lp}}^{2l}  +  \left\|\tilde{v}(s) \right\|_{L^{2lp}}^{2l(1-\beta_0)}  \left\| v(s)\right\|_{L^{\frac{2kp}{2-p}}}^{2k}  \left\|  \tilde{v}(s) ^{2(l-1)} |\nabla \tilde{v}(s) |^2 \right\|_{L^1}^{\beta_0} \right.  \\ &\left. \quad+  \left\|v(s) \right\|_{L^{2kp}}^{2k(1-\beta_0)}  \left\| \tilde{v}(s)\right\|_{L^{\frac{2lp}{2-p}}}^{2l}  \left\|  v(s) ^{2(k-1)} |\nabla v(s) |^2 \right\|_{L^1}^{\beta_0}\right) \left\| Z^{:m+1-i,m-j:}(s) \right\|_{\mathcal{C}^{-\alpha}}^2
				\mathrm{d}s.
			\end{align*}
			Taking the constants $p^1_i>1$ for $1\leq i\leq 3$, which satisfy  $1=\sum_{i=1}^{3}\frac{1}{p^1_i}$ and $p^1_1,p^1_2 \geq p$, we have that 
			\begin{align*}
				&\int_{1}^{t}  \left\|v(s) \right\|_{L^{2kp}}^{2k} 
				\left\|\tilde{v}(s) \right\|_{L^{2lp}}^{2l}\left\| Z^{:m+1-i,m-j:}(s) \right\|_{\mathcal{C}^{-\alpha}}^2
				\mathrm{d}s
				\\\lesssim&\, 
				\left( \int_{1}^{t} \left\|v(s) \right\|_{L^{2kp}}^{2kp^1_1} \mathrm{d}s\right) ^{\frac{1}{p^1_1}}   \left( \int_{1}^{t} \left\|\tilde{v}(s) \right\|_{L^{2l p}}^{2lp^1_2} \mathrm{d}s \right) ^{\frac{1}{p^1_2}} \left(\int_{1}^{t} 	\left\| Z^{:m+1-i,m-j:}(s) \right\|_{\mathcal{C}^{-\alpha}}^{2p^1_3} \mathrm{d}s \right) ^{\frac{1}{p^1_3}}	\\
				\lesssim&\, 
				\left( \int_{1}^{t} \left\|v(s) \right\|_{L^{2kp^1_1}}^{2kp^1_1} \mathrm{d}s\right) ^{\frac{1}{p^1_1}}   \left( \int_{1}^{t} \left\|\tilde{v}(s) \right\|_{L^{2l p^1_2}}^{2lp^1_2} \mathrm{d}s \right) ^{\frac{1}{p^1_2}} \left(\int_{1}^{t} 	\left\| Z^{:m+1-i,m-j:}(s) \right\|_{\mathcal{C}^{-\alpha}}^{2p^1_3} \mathrm{d}s \right) ^{\frac{1}{p^1_3}}.
			\end{align*}
			Taking the constants $p^2_i>1$ for $1\leq i\leq 4$ which satisfy  $1=\sum_{i=1}^{4}\frac{1}{p^2_i}$, $p^2_1\geq \frac{p}{2-p}$, $(1-\beta_0)p^2_2\geq p$ and $\beta_0p^2_3\leq 1$, we have that 
			\begin{align*}
				&\int_{1}^{t}    \left\|\tilde{v}(s) \right\|_{L^{2lp}}^{2l(1-\beta_0)}  \left\| v(s)\right\|_{L^{\frac{2kp}{2-p}}}^{2k}  \left\|  \tilde{v}(s) ^{2(l-1)} |\nabla \tilde{v}(s) |^2 \right\|_{L^1}^{\beta_0}  \left\| Z^{:m+1-i,m-j:}(s) \right\|_{\mathcal{C}^{-\alpha}}^2
				\mathrm{d}s
				\\\lesssim&\, 
				\left( \int_{1}^{t} \left\|v(s) \right\|_{L^{\frac{2kp}{2-p}}}^{2kp^2_1} \mathrm{d}s\right) ^{\frac{1}{p^2_1}}  \left( \int_{1}^{t}  \left\|\tilde{v}(s) \right\|_{L^{2lp}}^{2l(1-\beta_0)p^2_2} \mathrm{d}s \right) ^{\frac{1}{p^2_2}} \\&\cdot
				\left( \int_{1}^{t} \left\| \tilde{v}(s)^{2(l-1)}  \left|\nabla \tilde{v}(s)\right| ^2 \right\|^{\beta_0 p^2_3}_{L^1}\mathrm{d}s \right) ^{\frac{1}{p^2_3}}
				\left(\int_{1}^{t} 	\left\| Z^{:m+1-i,m-j:}(s) \right\|_{\mathcal{C}^{-\alpha}}^{2p^2_4} \mathrm{d}s \right) ^{\frac{1}{p^2_4}}	\\\lesssim&\, 
				\left( \int_{1}^{t} \left\|v(s) \right\|_{L^{2kp^2_1}}^{2kp^2_1} \mathrm{d}s\right) ^{\frac{1}{p^2_1}}  \left( \int_{1}^{t}  \left\|\tilde{v}(s) \right\|_{L^{2l(1-\beta_0)p^2_2}}^{2l(1-\beta_0)p^2_2} \mathrm{d}s \right) ^{\frac{1}{p^2_2}} \\&\cdot
				\left( \int_{1}^{t} \left\| \tilde{v}(s)^{2(l-1)}  \left|\nabla \tilde{v}(s)\right| ^2 \right\|_{L^1}\mathrm{d}s +t \right) ^{\frac{1}{p^2_3}}
				\left(\int_{1}^{t} 	\left\| Z^{:m+1-i,m-j:}(s) \right\|_{\mathcal{C}^{-\alpha}}^{2p^2_4} \mathrm{d}s \right) ^{\frac{1}{p^2_4}}.
			\end{align*}
			Similarly, taking the constants $p^3_i>1$ for  $1\leq i\leq 4$, which satisfy  $1=\sum_{i=1}^{4}\frac{1}{p^3_i}$, $(1-\beta_0)p^3_1\geq p$, $p^3_2\geq \frac{p}{2-p}$ and $\beta_0p^3_3\leq 1$, we have that 
			\begin{align*}
				&\int_{1}^{t}    \left\|v(s) \right\|_{L^{2kp}}^{2k(1-\beta_0)}  \left\| \tilde{v}(s)\right\|_{L^{\frac{2lp}{2-p}}}^{2l}  \left\|  v(s) ^{2(k-1)} |\nabla v(s) |^2 \right\|_{L^1}^{\beta_0} \left\| Z^{:m+1-i,m-j:}(s) \right\|_{\mathcal{C}^{-\alpha}}^2
				\mathrm{d}s
				\\\lesssim&\, \left( \int_{1}^{t} \left\|v(s) \right\|_{L^{2k(1-\beta_0)p^3_1}}^{2k(1-\beta_0)p^3_1} \mathrm{d}s\right) ^{\frac{1}{p^3_1}} 	\left( \int_{1}^{t}  \left\|\tilde{v}(s) \right\|_{L^{2l p^3_2}}^{2 l p^3_2}\mathrm{d}s \right) ^{\frac{1}{p^3_2}}  \\&
				\left( \int_{1}^{t}  \left\| v(s)^{2k-2}  \left|\nabla v(s)\right| ^2 \right\|_{L^1} \mathrm{d}s+t \right) ^{\frac{1}{p^3_3}} 
				\left(\int_{1}^{t} 	\left\| Z^{:m+1-i,m-j:}(s) \right\|_{\mathcal{C}^{-\alpha}}^{2p^3_4} \mathrm{d}s \right) ^{\frac{1}{p^3_4}}.
			\end{align*}
			Let $\mu_0>0$ be such that 
			\begin{align}
				&1+\mu_0\left( \mu_0+\sqrt{1+\mu_0^2}\right) \label{condition on mu}\\= &\, 2m \max\left\lbrace p^1_1, p^1_2, p^2_1,(1-\beta_0)p^2_2,(1-\beta_0)p^3_1,p^3_2 \right\rbrace  -m.\nonumber
			\end{align}
			Let $\mu>\mu_0$. This implies that $2+2\mu(\mu+\sqrt{1+\mu^2})>2m+1$.
			Then we can use \eqref{eq: Lp+2m est 2} in Proposition \ref{priori estimate} and Proposition \ref{priori estimate for auxiliary system} to control terms like
			\begin{align*}
				&\int_{1}^{t} \left\| v(s)\right\|^{2q+2m}_{L^{2q+2m}}\mathrm{d}s,\quad \int_{1}^{t} \left\|\tilde{v}(s)\right\|^{2q+2m}_{L^{2q+2m}} \mathrm{d}s,\\&\int_{1}^{t} \left\| v(s)^{2q-2}  \left|\nabla v(s)\right| ^2 \right\|_{L^1} \mathrm{d}s,\quad \int_{1}^{t} \left\| \tilde{v}(s)^{2q-2}  \left|\nabla \tilde{v}(s)\right| ^2 \right\|_{L^1} \mathrm{d}s,
			\end{align*}
			for some $q\in \left[1, 1+\mu\left( \mu+\sqrt{1+\mu^2}\right)\right]$. Finally, we obtain that there exists $\gamma>0$ such that on $E_{K,\gamma}$,
			\begin{align*}
				\int_{1}^{t}  L(s)\mathrm{d}s&\leq C\Bigg( t+ \left\| \underline{Z}\right\|_{\alpha,\alpha^{\prime},1}^{\gamma}+ \sum_{i=0}^{m+1}\sum_{j=0}^{m} \int_{1}^t \left\| Z^{:m+1-i,m-j:}(s) \right\| _{\mathcal{C}^{-\alpha}}^{\gamma} \mathrm{d} s \Bigg)\leq C(1+t).
			\end{align*}
			Combining with \eqref{Gronwall}, we obtain that on $E_{K,\gamma}$, we can take $\lambda$ large enough such that there exist constants $C_1, C_2 > 0$ such that
			\begin{equation*}
				\|w(t)\|_{L^2}^2\leq C_1e^{-C_2 t}\rightarrow 0, \quad t\rightarrow \infty.
			\end{equation*}
			This implies that there exists $R>0$ such that $\tau_R=\infty$ on $E_{K,\gamma}$.
			
			Moreover, under the assumptions that $2m+1<2+2\mu(\mu+\sqrt{1+\mu^2})$ and $\alpha_0\in [1/( 1+\mu(\mu+\sqrt{1+\mu^2})) ,2/\left( 2m+1\right) ]$, we can find $1\leq p\leq 1+\frac12\mu(\mu+\sqrt{1+\mu^2}) $ such that $L^{2p}\hookrightarrow\mathcal{C}^{-\alpha_0}$.
			Thus on $E_{K,\gamma}$, by H\"older's inequality, \eqref{eq: Lp+2m est 2} in Proposition \ref{priori estimate} and  Proposition \ref{priori estimate for auxiliary system}, we know that
			\begin{align*}
				\left\| w(t)\right\|^{2p}_{\mathcal{C}^{-\alpha_0}}& \lesssim \left\| w(t)\right\|_{L^{2p}} ^{2p}\lesssim \left\| w(t)\right\|_{L^{2}} \left(  \left\| v(t)\right\|_{L^{2(2p-1)}} ^{2p-1} + \left\| \tilde{v}(t)\right\|_{L^{2(2p-1)}} ^{2p-1}\right)\\& \lesssim e^{-Ct}(1+t)\rightarrow 0, \quad t\rightarrow \infty.
			\end{align*}
			Therefore, we obtain that 
			\begin{equation*}
				\Gamma_{u_0, u_1}(D)\geq \mathrm{P}(E_{K,\gamma})>0.
			\end{equation*}
			By Theorem \ref{HM Corollary 2.2}, we get the uniqueness of the invariant probability measure.
		\end{proof}

		\appendix
		
		\section{Some calculations and technical estimates}\label{Some technical estimation}

		\begin{proof}[Proof of \eqref{eq: 2nd est with test func}] For $j=(j_1,j_2)\in \mathbb{Z}^2$, we define $|j|=\sqrt{j_1^2+j_2^2}$ and set $\rho_j =1+\mu|j|^2$, $\theta_j=|j|^2$. By isometry property \eqref{eq: general isometry property},
			\begin{align}
				&\mathrm{E} \left[ \left\langle Z_{-\infty}^{:k,l:}(t,\cdot),\phi_1\right\rangle \overline{\left\langle Z_{-\infty}^{:k,l:}(s,\cdot),\phi_2\right\rangle}\right]  \label{eq: 2nd est with test func in App} \\
				=&\, \mathrm{E} \left[ I_{k,l}(f^{k,l}_{-\infty,t,\phi_1}) \overline{I_{k,l}(f^{k,l}_{-\infty,s,\phi_2}) } \right] = k!l!\left\langle f^{k,l}_{-\infty,t,\phi_1}, \overline{f^{k,l}_{-\infty,s,\phi_2} }\right\rangle \nonumber
				\\=&\,k!l!\int_{\left\lbrace (-\infty,t\wedge s]\times \mathbb{T}^2\right\rbrace ^{k+l}} \int_{\mathbb{T}^2} \prod_{j=1}^{k} H(t-s_j,x_1-y_j)\prod_{j=k+1}^{k+l} \overline{H}(t-s_j,x_1-y_j)\phi_1(x_1) \mathrm{d}x_1 \nonumber \\& \cdot \int_{\mathbb{T}^2} \prod_{j=1}^{k} \overline{H}(s-s_j,x_2-y_j)\prod_{j=k+1}^{k+l} H(s-s_j,x_2-y_j) \overline{\phi_2(x_2)}\mathrm{d}x_2 \mathrm{d}s_1\mathrm{d}y_1\cdots\mathrm{d}s_{k+l}\mathrm{d}y_{k+l} \nonumber
				\\=&\,k!l!\int_{ \mathbb{T}^2}\int_{ \mathbb{T}^2}\left( \int_{-\infty}^{t\wedge s}\int_{ \mathbb{T}^2} H(t-s_j,x_1-y_j)\overline{H}(s-s_j,x_2-y_j)\mathrm{d}y_j\mathrm{d}s_j\right)  ^k \nonumber \\&\cdot\left( \int_{-\infty}^{t\wedge s}\int_{ \mathbb{T}^2} \overline{H}(t-s_j,x_1-y_j)H(s-s_j,x_2-y_j)\mathrm{d}y_j\mathrm{d}s_j\right) ^l  \phi_1(x_1) \overline{\phi_2(x_2)}\mathrm{d}x_1\mathrm{d}x_2 \nonumber
				\\=&\,k!l!\int_{ \mathbb{T}^2}\int_{ \mathbb{T}^2}
				\Bigg( \sum_{j \in \mathbb{Z}^2}\frac{e^{-\rho_j|t-s|-\i\theta_j(t-s)}}{2\rho_j}  \frac{1}{4\pi^2} e^{\i j\cdot(x_1-x_2)}\Bigg)^k \nonumber\\& \cdot\Bigg( \sum_{j \in \mathbb{Z}^2}\frac{e^{-\rho_j|t-s|+\i\theta_j(t-s)}}{2\rho_j}  \frac{1}{4\pi^2} e^{\i j\cdot(x_1-x_2)}\Bigg) ^l \phi_1(x_1) \overline{\phi_2(x_2)}\mathrm{d}x_1\mathrm{d}x_2 \nonumber \\
				=&\, \frac{k!l!}{(2\pi)^{2(k+l)}}\sum_{j_1,\ldots, j_{k+l}\in \mathbb{Z}^2}   \prod_{i=1}^{k+l}\frac{ e^{-\rho_{j_i}|t-s| -\i a_i\theta_{j_i} (t-s) }}{2\rho_{j_i}}  \nonumber \\
				& \cdot\int_{ \mathbb{T}^2}\int_{ \mathbb{T}^2} e^{\i (\sum_{i=1}^{k+l} j_i)\cdot x_1} e^{-\i (\sum_{i=1}^{k+l} j_i)\cdot x_2} \phi_1(x_1) \overline{\phi_2(x_2)}\mathrm{d}x_1\mathrm{d}x_2 \nonumber\\
				=&\, \frac{k!l!}{(2\pi)^{2(k+l-1)}}\sum_{j\in \mathbb{Z}^2} \hat{\phi_1}(-j) \overline{\hat{\phi_2}(-j)}   \nonumber\Bigg( \sum_{j_1+\cdots+ j_{k+l}=j}\prod_{i=1}^{k+l}\frac{ e^{ -\rho_{j_i}|t-s|-\i a_i\theta_{j_i} (t-s) }}{2\rho_{j_i}}   \Bigg) \nonumber,
			\end{align}
			where $a_i=1$ for $1\leq i\leq k$ and $a_i=-1$ for $k+1\leq i\leq k+l$.
		\end{proof}
		
		For ease of notations, we write
		\begin{align}\label{A^{k,l}}
			A^{k,l}_{s,t}(j_1,\ldots,j_{k+l}):=&2\prod_{i=1}^{k+l}\frac{1}{2\rho_{j_i}}  - 2\cos\Bigg( (t-s)\sum_{i=1}^{k+l}a_i \theta_{j_i} \Bigg)  \prod_{i=1}^{k+l}\frac{e^{-\rho_{j_i}|t-s|} }{2\rho_{j_i}} ,
		\end{align}
		where $k,l\in\mathbb{N}$ with $k+l>0$, $j_1,\ldots,j_{k+l}\in \mathbb{Z}^2$, $s,t>-\infty$.
		\begin{lemma}\label{lem bound for A^{k,l}}
			For $A^{k,l}_{s,t}(j_1,\ldots,j_{k+l})$ defined in \eqref{A^{k,l}}, we have that for any $0\leq \epsilon\leq 1$,
			\begin{equation}\label{eq bound for A^{k,l}}
				A^{k,l}_{s,t}(j_1,\ldots,j_{k+l})\lesssim |t-s|^{\epsilon}\prod_{i=1}^{k+l}\frac{1}{(1+|j_i|^2)^{1-\epsilon}}.
			\end{equation}
		\end{lemma}
		\begin{proof}
			If $s=t$, $	A^{k,l}_{s,t}(j_1,\ldots,j_{k+l})\equiv0$, then \eqref{eq bound for A^{k,l}} holds. Without loss of generality, we assume that $-\infty< s<t<\infty$.
			Note that 
			\begin{align*}
				&A^{k,l}_{s,t}(j_1,\ldots,j_{k+l})\\
				=&\, 
				\prod_{i=1}^{k+l}\int_{-\infty}^{t} e^{-2\rho_{j_i}(t-v)}\mathrm{d} v+ 	\prod_{i=1}^{k+l}\int_{-\infty}^{s} e^{-2\rho_{j_i}(s-v)}\mathrm{d} v \\&\, - 2\cos\left( (t-s)\sum_{i=1}^{k+l}a_i \theta_{j_i} \right)  	\prod_{i=1}^{k+l}  \int_{-\infty}^{ s} e^{-\rho_{j_i}(t+s-2v)}\mathrm{d} v \\
				=&\, \int_{(-\infty,t]^{k+l}} e^{-\sum_{i=1}^{k+l}2\rho_{j_i}(t-v_i) }\mathrm{d} v_1\cdots \mathrm{d} v_{k+l} + \int_{(-\infty,s]^{k+l}} e^{-\sum_{i=1}^{k+l}2\rho_{j_i}(s-v_i) }\mathrm{d} v_1\cdots \mathrm{d} v_{k+l}\\&\,  - 2\cos\left( (t-s)\sum_{i=1}^{k+l}a_i \theta_{j_i} \right)  \int_{(-\infty,s]^{k+l}} e^{-\sum_{i=1}^{k+l}\rho_{j_i}(t+s-2v_i) }\mathrm{d} v_1\cdots \mathrm{d} v_{k+l}\\
				= &\, \left[ \left( 1- e^{-\sum_{i=1}^{k+l}\rho_{j_i}(t-s) }\right) ^2+ 2e^{-\sum_{i=1}^{k+l}\rho_{j_i}(t-s) }\left( 1-  \cos\left( (t-s)\sum_{i=1}^{k+l}a_i \theta_{j_i} \right)\right) \right]\prod_{i=1}^{k+l}\frac{1}{ 2\rho_{j_i} } \\&\, + \int_{(-\infty,t]^{k+l}\setminus (-\infty,s]^{k+l}} e^{-\sum_{i=1}^{k+l}2\rho_{j_i}(t-v_i) }\mathrm{d} v_1\cdots \mathrm{d} v_{k+l} .
			\end{align*}
			Since for $(v_1,\ldots, v_{k+l})\in (-\infty,t]^{k+l}\setminus (-\infty,s]^{k+l}$, there at least exists some $1\leq i_0\leq k+l$ such that $v_{i_0}\in [s,t]$, then 
			\begin{align*}
				&\int_{(-\infty,t]^{k+l}\setminus (-\infty,s]^{k+l}} e^{-\sum_{i=1}^{k+l}2\rho_{j_i}(t-v_i) }\mathrm{d} v_1\cdots \mathrm{d} v_{k+l} \\
				\leq &\, \int_{s}^{t} e^{-2\rho_{j_{i_0}}(t-v_{i_0}) }\mathrm{d} v_{i_0} \prod_{i\neq i_0} \int_{-\infty}^{t} e^{-2\rho_{j_i}(t-v)}\mathrm{d} v\\
				=&\, \frac{1}{\prod_{i=1}^{k+l} 2\rho_{j_i} } \left(1- e^{-2\rho_{j_{i_0}}(t-s)}\right) \\
				\leq&\,\frac{1}{\prod_{i=1}^{k+l} 2\rho_{j_i} } \left(1- e^{-\sum_{i=1}^{k+l}2\rho_{j_i}(t-s) } \right) .
			\end{align*}
			Using the fact that for any $0\leq \epsilon\leq 1$, $1-e^{-z}\leq z^{ \epsilon}$ and $ 1-\cos z\leq 2z^{ \epsilon}$ for $z\geq0$, we continue the estimate and get that
			\begin{align*}
				A^{k,l}_{s,t}(j_1,\ldots,j_{k+l})
				&\lesssim \left[ 1- e^{-\sum_{i=1}^{k+l}2\rho_{j_i}(t-s) }+ \left( 1-  \cos\left( (t-s)\sum_{i=1}^{k+l}a_i \theta_{j_i} \right)\right) \right]  \prod_{i=1}^{k+l} \frac{1}{2\rho_{j_i} }   \\
				&\lesssim (t-s)^{\epsilon} \left(\sum_{i=1}^{k+l} \left(1+|j_i|^2\right) \right) ^{\epsilon} \prod_{i=1}^{k+l} \frac{1}{1+|j_i|^2 }\\
				&\lesssim (t-s)^{\epsilon} \left(\sum_{i=1}^{k+l} \left(1+|j_i|^2\right)^{\epsilon}\right) \prod_{i=1}^{k+l} \frac{1}{1+|j_i|^2 }\\ 
				&\lesssim (t-s)^{\epsilon} \prod_{i=1}^{k+l}\frac{1}{(1+|j_i|^2)^{1-\epsilon}}.
			\end{align*}
			Then we finish the proof.
		\end{proof}

		\begin{proof}[Proof of Lemma \ref{lem J_{k,l}(z_1+z_2)}]
			Using \eqref{complex Hermite polynomial}, we know that
			\begin{equation}\label{complex Hermite polynomial1}
				\exp\left\{ \lambda(\overline{z_1}+\overline{z_2})+\overline{\lambda}(z_1+z_2)-\rho|\lambda|^2\right\} =\sum_{k=0}^{\infty}\sum_{l=0}^{\infty}\frac{\overline{\lambda}^k\lambda^l}{k!l!}J_{k,l}(z_1+z_2,\rho).
			\end{equation}
			Further expanding the left hand side of \eqref{complex Hermite polynomial1}, we get that
			\begin{align}
				&e^{ \lambda(\overline{z_1}+\overline{z_2})+\overline{\lambda}(z_1+z_2)-\rho|\lambda|^2}= e^{ \lambda\overline{z_1}+\overline{\lambda}z_1-\rho|\lambda|^2} e^{ \lambda\overline{z_2}} e^{\overline{\lambda}z_2} \nonumber\\
				=&\,  \sum_{i=0}^{\infty} \sum_{j=0}^{\infty} \sum_{m=0}^{\infty} \sum_{n=0}^{\infty}  \frac{J_{i,j}(z_1,\rho)}{i!j!m!n!}z_2^n\overline{z_2}^m  \overline{\lambda}^{i+n}\lambda^{j+m} \nonumber\\
				=&\,  \sum_{k=0}^{\infty} \sum_{l=0}^{\infty} \left( \sum_{i=0}^{k} \sum_{j=0}^{l}  \frac{J_{i,j}(z_1,\rho)}{k!l!}\binom{k}{i}\binom{l}{j}z_2^{k-i}\overline{z_2}^{l-j} \right)  \overline{\lambda}^{k}\lambda^{l} \label{complex Hermite polynomial2}.
			\end{align}
			Identifying coefficients of $\overline{\lambda}^{k}\lambda^{l}$ in \eqref{complex Hermite polynomial1} and \eqref{complex Hermite polynomial2}, we derive the conclusion.
		\end{proof}
		
		\begin{proof}[Proof of Lemma \ref{lem: unique for u}]
			By the definition of $u(\cdot;u_0)$, we know that
			\begin{align*}
				u(t+h;u_0)&= Z(t+h)+v(t+h;u_0)\\
				&= Z_t(t+h)+ \left( P_h u(t;u_0)+ \int_{0}^{h} P_{h-r} \Psi\left(v(t+r),\underline{Z}(t+r)\right)\mathrm{d}r\right) \\
				&=:Z_t(t+h)+ v_t(t+h),
			\end{align*}
			which implies 
			\begin{equation}\label{v_t(t+h)}
				v_t(t+h)= Z(t+h)+v(t+h;u_0)- Z_t(t+h)= v(t+h;u_0)+ P_h Z(t).
			\end{equation}
			Applying \eqref{v_t(t+h)} and Corollary \ref{Cor relation of wick product of diff initial moments}, we derive that
			\begin{align*}
				&\sum_{i=0}^{m+1}\sum_{j=0}^{m}\binom{m+1}{i}\binom{m}{j}v(t+r)^i\overline{v(t+r)}^jZ^{:m+1-i,m-j:}(t+r) \label{Psi time shift in App} \\
				=&\,\sum_{i=0}^{m+1}\sum_{j=0}^{m}\binom{m+1}{i}\binom{m}{j}(	v_t(t+r)- P_r Z(t))^i(\overline{	v_t(t+r)}-\overline{ P_r Z(t)})^jZ^{:m+1-i,m-j:}(t+r) \nonumber\\
				=&\,\sum_{i=0}^{m+1}\sum_{j=0}^{m}\binom{m+1}{i}\binom{m}{j} \sum_{k=0}^{i}\sum_{l=0}^{j}\binom{i}{k}\binom{j}{l}  	v_t(t+r)^k(- P_r Z(t))^{i-k}\overline{	v_t(t+r)}^l(-\overline{ P_r Z(t)})^{j-l} \nonumber\\& \, \cdot Z^{:m+1-i,m-j:}(t+r) \nonumber\\
				=&\,\sum_{k=0}^{m+1}\sum_{l=0}^{m}\sum_{s=0}^{m+1-k}\sum_{n=0}^{m-l}\binom{m+1}{s+k}\binom{m}{n+l} \binom{s+k}{k}\binom{n+l}{l}  	v_t(t+r)^k(- P_r Z(t))^{s}\overline{	v_t(t+r)}^l  \nonumber \\& \, \cdot(-\overline{ P_r Z(t)})^{n} Z^{:m+1-s-k,m-n-l:}(t+r) \nonumber\\
				=&\, \sum_{k=0}^{m+1}\sum_{l=0}^{m}\binom{m+1}{k}\binom{m}{l} v_t(t+r)^k \overline{	v_t(t+r)}^l  \\& \,\cdot\left( \sum_{s=0}^{m+1-k}\sum_{n=0}^{m-l}\binom{m+1-k}{s}\binom{m-l}{n}  	(- P_r Z(t))^{s}(-\overline{ P_r Z(t)})^{n} Z^{:m+1-s-k,m-n-l:}(t+r) \right)  \nonumber \\
				=&\, \sum_{k=0}^{m+1}\sum_{l=0}^{m}\binom{m+1}{k}\binom{m}{l} v_t(t+r)^k \overline{	v_t(t+r)}^l Z_t^{:m+1-k,m-l:}(t+r)  \nonumber  .
			\end{align*}
			Combining with $v(t+r)+Z(t+r)= v_t(t+r)+Z_t(t+r)$, we have 
			\begin{equation*}
				\Psi\left(v(t+r),\underline{Z}(t+r)\right)= \Psi\left(v_t(t+r),\underline{Z}_t(t+r)\right).
			\end{equation*}
			Therefore, 
			\begin{align*}
				v_t(t+h)&= P_h u(t;u_0)+ \int_{0}^{h} P_{h-r} \Psi\left(v(t+r),\underline{Z}(t+r)\right)\mathrm{d}r\\
				&= P_h u(t;u_0) + \int_{0}^{h} P_{h-r} \Psi\left(v_t(t+r),\underline{Z}_t(t+r)\right)\mathrm{d}r.
			\end{align*}
			Then we complete the proof.
		\end{proof}



	\end{document}